\documentclass[oupthm,crop,info]{CUP-JNL-CJM}%

\newcommand{\bpf}[1][Proof]{\noindent{\fontseries{sb}\selectfont #1\quad}}

\newcommand{\epf}{\qed\vspace{12pt}}

\usepackage{graphicx}
\usepackage{multicol,multirow}
\usepackage{amsmath,amssymb,amsfonts}
\usepackage{mathrsfs}
\usepackage{rotating}
\usepackage{appendix}

\theoremstyle{oupplain}
\newtheorem{theorem}{Theorem}[section]
\newtheorem{lemma}[theorem]{Lemma}
\newtheorem{corollary}[theorem]{Corollary}

\newtheorem{claim}[theorem]{Claim}
\theoremstyle{oupdefinition}

\theoremstyle{oupremark}

\theoremstyle{oupproof}

\numberwithin{equation}{section}

\usepackage{ifthen,latexsym,bbm}
\usepackage[shortlabels]{enumitem}

\usepackage[nobysame,initials]{amsrefs}

\usepackage{hyperref}
\usepackage{tikz}
\usepackage{dsfont}
\usepackage{anyfontsize}
\usepackage{array}

\AtBeginDocument{%
   \def\MR#1{}
}

\usepackage[normalem]{ulem}

\newcommand{\C}[1]{{\protect\mathcal{#1}}}
\newcommand{\B}[1]{{\bf #1}}
\newcommand{\I}[1]{{\mathbbm #1}}
\renewcommand{\O}[1]{\overline{#1}}

\newcommand{\me}{{\mathrm e}}
\renewcommand{\ldots}{\hspace{0.9pt}.\hspace{0.3pt}.\hspace{0.3pt}.\hspace{1.5pt}}
\renewcommand{\ge}{\geqslant}
\renewcommand{\le}{\leqslant}
\renewcommand{\geq}{\geqslant}
\renewcommand{\leq}{\leqslant}

\newcommand{\cG}[2]{{\mathcal G}^{(#1)}_{#2}}

\newcommand{\GR}[1]{\mathrm{GR}(#1)}

\newif\ifnotesw\noteswtrue

\newcommand{\hide}[1]{}

\newcommand{\beq}[1]{\begin{equation}\label{#1}}
\newcommand{\eeq}{\end{equation}}

\newcommand{\case}[2]{\medskip\noindent{\bf Case #1.} #2\medskip\par\noindent}

\newcommand{\E}[1]{\mathbb{E}\left[\,#1\,\right]}
\newcommand{\p}[2]{P_{#1}(#2)}
\newcommand{\pp}[3]{P_{\,\O{#1}#2}(#3)}
\newcommand{\CI}{C}
\newcommand{\OutIn}[2]{$\ifthenelse{\equal{#1}{}}{}{\O #1}#2$}
\newcommand{\OneClaim}{$1$-claim}
\newcommand{\OMerge}[2]{\ifthenelse{\equal{#1}{}}{(#2)}{(#1|#2)}}
\newcommand{\Merge}[2]{\ifthenelse{\equal{#1}{}}{#2}{#1|#2}}
\newcommand{\cluster}{cluster}
\newcommand{\OneClusters}{\mathcal{M}_{\Merge{}{1}}}
\newcommand{\OneCluster}{$\Merge{}{1}$-\cluster}
\newcommand{\TwoClusters}{\mathcal{M}_{\Merge{}{2}}}
\newcommand{\TwoCluster}{$\Merge{}{2}$-\cluster}
\newcommand{\ThreePClusters}{\mathcal{M}_{3^{+}}}
\newcommand{\ThreePCluster}{$3^{+}$-\cluster}
\newcommand{\OTwoPMerge}{\OMerge{}{{2^+}}}
\newcommand{\TwoPMerge}{\Merge{}{{2^+}}}
\newcommand{\TwoPClusters}{\mathcal{M}_{\TwoPMerge}}
\newcommand{\TwoPCluster}{$\TwoPMerge$-\cluster}
\newcommand{\TwoPClaim}{$2^+$-claim}
\newcommand{\OThreePMerge}{\OMerge{}{{3^{+}}}}
\newcommand{\ThreePMerge}{\Merge{}{{3^{+}}}}

\newcommand{\CC}[4]{{\mathcal C}_{#1,#2;#3,#4}}
\newcommand{\CS}[3]{{\mathcal S}_{#1,#2}\ifthenelse{\equal{#3}{}}{}{(#3)}}
\newcommand{\CCS}[1]{\C S_{4,3}'\ifthenelse{\equal{#1}{}}{}{(#1)}}
\newcommand{\CCP}{\CC{4}{3}{3}{3}'}

\def\ex{\mathrm{ex}}

\articletype{RESEARCH ARTICLE}
\jname{Canadian Mathematical Society}
\jyear{2025}

\begin{document}

\renewcommand{\mid}{:}

\begin{Frontmatter}

\title[On the $(k+2,k)$-problem of Brown, Erd\H{o}s and S\'os for $k=5,6,7$]{On the $(k+2,k)$-problem of Brown, Erd\H{o}s and S\'os for $k=5,6,7$}

\author{Stefan Glock\thanks{Funded by the Deutsche Forschungsgemeinschaft (DFG, German Research Foundation) -- 542321564.}}
\author{Jaehoon Kim\thanks{Supported by the National Research Foundation of Korea (NRF) grant funded by the
Korea government(MSIT) No. RS-2023-00210430.}}
\author{Lyuben Lichev}
\author{Oleg Pikhurko\thanks{Supported by ERC Advanced Grant 101020255.}}
\author{Shumin Sun\thanks{Supported by ERC Advanced Grant 101020255.}}

\authormark{S. Glock and J. Kim and L. Lichev and O. Pikhurko and S. Sun}

\address{\orgname{Fakult\"at f\"ur Informatik und Mathematik, Universit\"at Passau}, \orgaddress{\city{Passau}, \country{Germany}}}

\address{\orgname{Department of Mathematical Sciences, KAIST}, \orgaddress{\city{Daejeon}, \country{South Korea}}}

\address{\orgname{Institute of Science and Technology Austria (ISTA)}, \orgaddress{\city{Klosterneuburg}, \country{Austria}}}

\address{\orgname{Mathematics Institute and DIMAP, University of Warwick}, \orgaddress{\city{Coventry}, \country{UK}}}

\address{\orgname{Mathematics Institute,
University of Warwick}, \orgaddress{\city{Coventry}, \country{UK}}}

\keywords[AMS subject classification]{05C65, 05C70, 05C35}


\abstract{Let $f^{(r)}(n;s,k)$ denote the maximum number of edges in an $n$-vertex $r$-uniform hypergraph containing no subgraph with $k$ edges and at most $s$ vertices.  Brown, Erd\H{o}s and S\'os [\emph{New directions in the theory of graphs (Proc.\ Third Ann Arbor Conf., Univ.\ Michigan 1971)}, pp.\ 53--63, Academic Press 1973] conjectured that the limit $\lim_{n\rightarrow \infty}n^{-2}f^{(3)}(n;k+2,k)$ exists for all $k$. The value of the limit was previously determined for $k=2$ in the original paper of Brown, Erd\H{o}s and S\'os, for $k=3$ by Glock [\emph{Bull.\ Lond.\ Math.\ Soc.,} 51 (2019) 230--236] and for $k=4$ by Glock, Joos, Kim, K\"{u}hn, Lichev and Pikhurko [\emph{Proc.\ Amer.\ Math.\ Soc.,\ Series B}, 11 (2024) 173--186] while Delcourt and Postle [\emph{Proc.\ Amer.\ Math.\ Soc.}, 152 (2024), 1881--1891] proved the conjecture (without determining the limiting value).

In this paper, we determine the value of the limit in the Brown--Erd\H{o}s--S\'os Problem for $k\in \{5,6,7\}$. More generally, we obtain the value of $\lim_{n\rightarrow \infty}n^{-2}f^{(r)}(n;rk-2k+2,k)$ for all $r\geq 3$ and $k\in \{5,6,7\}$. In addition, by combining these new values with recent results of Bennett, Cushman and Dudek [arxiv:2309.00182, 2023] we obtain new asymptotic values for several generalised Ramsey numbers.}

\end{Frontmatter}

\section{Introduction}

Given a family $\C F$ of $r$-uniform hypergraphs (in short, \emph{$r$-graphs}), denote by $\ex(n;\C F)$ the \emph{Tur\'an number} of $\C F$, i.e.~the maximum number of edges in an $n$-vertex $r$-graph containing no element of $\C F$ as a subgraph. Tur\'an problems for hypergraphs are notoriously difficult and we still lack an understanding of even seemingly simple instances such as when $\C F$ forbids the complete $3$-graph on $4$ vertices.
We refer the reader to the surveys~\cites{keevash,sido} for more background.
In this paper, we focus on the family $\C F^{(r)}(s,k)$ of all $r$-graphs with $k$ edges and at most $s$ vertices. 

Brown, Erd\H{o}s and S\'os~\cite{BES} launched the systematic study of the function
$$f^{(r)}(n;s,k):=\ex(n;\C F^{(r)}(s,k)).$$
 The case $r=2$ (resp.\ $r=3$) of the problem was previously studied by Erd\H{o}s~\cite{Erdos64a} (resp.\ Brown, Erd\H{o}s and S\'os~\cite{BrownErdosSos73}). 
Since then, the asymptotics of $f^{(r)}(n;s,k)$ as $n\to\infty$ have been intensively investigated for various natural choices of parameters $r,s,k$  (see e.g.~\cites{ASh,CGLS,EFR,G,GJKKLP,KeevashLong20x,NRS,RS,ShTamo,SidoBES}). For instance, it includes the celebrated $(6,3)$-theorem of Ruzsa and Szemer\'edi~\cite{RS} (namely, when $(r,s,k)=(3,6,3)$), as well as the notoriously difficult $(7,4)$-problem (namely, when $(r,s,k)=(3,7,4)$). Beyond its significance of being a fundamental Tur\'an problem, the Brown--Erd\H{o}s--S\'os function is closely related to problems from other areas such as  additive combinatorics (see e.g.~\cite{RS}), coding theory (e.g.\ the case $k=2$), hypergraph packing and designs (see below). 

Brown, Erd\H{o}s and S\'os~\cite{BES} proved  that 
$$
\Omega(n^{(rk-s)/(k-1)})=f^{(r)}(n;s,k)=O(n^{\lceil (rk-s)/(k-1)\rceil}).
$$
In this paper, we are interested in the case when the exponent in both the lower and the upper  bound is equal to $2$, i.e.\ $s=rk-2k+2$. 
In this setting, the natural question is whether $n^{-2} f^{(r)}(n;rk-2k+2,k)$ converges to a limit as $n\to \infty$; in fact, already Brown, Erd\H{o}s and S\'os~\cite{BES} considered this question and conjectured that the limit exists for $r=3$.
They verified their conjecture for $k=2$ by showing that the limit is $1/6$. Glock~\cite{G} proved that, when $k=3$, the limit exists and is equal to $1/5$. Recently, Glock, Joos, Kim, K\"{u}hn, Lichev and Pikhurko~\cite{GJKKLP} solved the case $k=4$ by showing that the limit equals to $7/36$. 

Already the original work of Brown, Erd\H{o}s and S\'os~\cites{BES,BrownErdosSos73} pointed connections to (approximate) designs: in particular, it was observed by them that Steiner triple systems (when they exist) give extremal examples for $k=2$. 
More generally, the celebrated theorem of R\"odl~\cite{Rodl85} which solved the Erd\H{o}s--Hanani problem from 1965 on asymptotically optimal clique coverings of complete hypergraphs can be phrased as
\begin{equation}\label{eq:Rodl85}
 \lim_{n\rightarrow \infty}n^{-t}f^{(r)}(n;2r-t,2)=\frac{(r-t)!}{r!}.
\end{equation} 
 Furthermore, the more recent results in~\cites{G,GJKKLP} linked the Brown--Erd\H{o}s--S\'os Problem to almost optimal graph packings. Namely,
in~\cite{G}, $F$-packings of complete graphs with some special graph $F$ are used, and in~\cite{GJKKLP}, 
a significant strengthening was needed to find ``high-girth'' packings. These structures are related to another famous problem of Erd\H{o}s in design theory, namely the existence of high-girth Steiner triple systems, which was recently resolved by Kwan, Sah, Sawhney and Simkin~\cite{KSSS}.

In a recent breakthrough, Delcourt and Postle~\cite{DP} proved the Brown--Erd\H{o}s--S\'os conjecture, namely, that for $r=3$ and any $k\geq 2$ the limit exists,  without determining its value. Moreover, as observed by Shangguan~\cite{sh}, their approach generalises to every uniformity $r\ge 4$. Thus the limit $\lim_{n\rightarrow \infty}n^{-2}f^{(r)}(n;rk-2k+2,k)$ exists for all $r\geq 3$ and $k\geq 2$.

While the existence of the limits is an important step forward, it would be very interesting to actually determine the limiting values, in particular in view of the fact that only few asymptotic results on degenerate hypergraph Tur\'an problems of quadratic growth are currently  known.

In this paper, we determine the limit for $k=5,6,7$ and arbitrary uniformity $r\geq 3$, as given by the following four theorems.
(Recall that the limit for $k=2$ is given in~\eqref{eq:Rodl85} while the cases $k=3,4$ were settled in~\cites{G,GJKKLP,ShTamo}.)

The following two results show that, for $k=5,7$, the limiting value is the same as for $k=3$.

\begin{theorem}\label{5edge}
For every $r\geq 3$, we have $\lim_{n\rightarrow \infty}n^{-2}f^{(r)}(n;5r-8,5)=\frac{1}{r^2-r-1}$.
\end{theorem}

\begin{theorem}\label{7edge}
    For every $r\geq 3$, we have $\lim_{n\rightarrow \infty}n^{-2}f^{(r)}(n;7r-12,7)=\frac{1}{r^2-r-1}$.
\end{theorem}

However, the case $k=6$ exhibits different behaviour when $r = 3$ and $r\geq 4$ (which parallels the situation for $k=4$), as established by the following two theorems.

\begin{theorem}\label{6edge3uniform}
    $\lim_{n\rightarrow \infty}n^{-2}f^{(3)}(n;8,6)=\frac{61}{330}$.
\end{theorem}

\begin{theorem}\label{6edgehigh}
    For every $r\geq 4$, we have $\lim_{n\rightarrow \infty}n^{-2}f^{(r)}(n;6r-10,6)=\frac{1}{r^2-r}$.
\end{theorem}

Very recently, Letzter and Sgueglia~\cite{letzter2023problem} proved various results on the existence and the value of the limit $\lim_{n\to\infty} n^{-t} f^{(r)}(n,k(r-t)+t,k)$. In particular, for $t=2$, they independently re-proved our upper bounds in Theorems~\ref{5edge} and \ref{7edge} when $r$ is sufficiently large,
and showed that $f^{(r)}(n; kr-2k+2,k)=(\frac{1}{r^2-r}+o(1))n^2$ when $k$ is even and $r\ge r_0(k)$ is large enough.\medskip

\noindent \textbf{An application to generalised Ramsey numbers.}
The following generalisation of Ramsey numbers was introduced by Erd\H{o}s and Shelah~\cite{ES75}, and  its systematic study was initiated by Erd\H{o}s and Gy\'arf\'as~\cite{EG97}.
Fix integers $p, q$ such that $p \geq 3$ and $2\le q\le\tbinom{p}{2}$. 
A \emph{$(p, q)$-colouring} of $K_n$
is a colouring of the edges of $K_n$ such that every $p$-clique has at least $q$ distinct colours among
its edges. 
The \emph{generalised Ramsey number} $\GR{n,p,q}$ is the minimum number of colours such that $K_n$ has a $(p, q)$-colouring. One relation to the classical Ramsey numbers is that
$\GR{n,p,2}>t$ if and only if every $t$-colouring of the edges of $K_n$ yields a monochromatic clique of order~$p$. 

In their work, Erd\H{o}s and Gy\'arf\'as~\cite{EG97} showed that, for every $p\geq 3$ and $q_{\rm{lin}} := 
\tbinom{p}{2}-p+3$, 
\[\GR{n,p,q_{\rm{lin}}} = \Omega(n)\quad \text{and}\quad \GR{n,p,q_{\rm{lin}}-1} = o(n),\]
while for every $p\geq 3$ and $q_{\rm{quad}} := 
\tbinom{p}{2}-\lfloor p/2\rfloor+2$,
\[\GR{n,p,q_{\rm{quad}}} = \Omega(n^2)\quad \text{and}\quad \GR{n,p,q_{\rm{quad}}-1} = o(n^2).\]
 Thus, $q_{\rm{lin}}$ and $q_{\rm{quad}}$ are the thresholds (that is, the smallest values of $q$) for $\GR{n,p,q}$ to be respectively linear and quadratic in~$n$.

Very recently, Bennett, Cushman and Dudek~\cite{BCD23} found the following connection between generalised Ramsey numbers and the Brown--Erd\H{o}s--S\'os function.

\begin{theorem}[\cite{BCD23}*{Theorem 3}]
\label{th:quadratic}
For all even $p\geq 6$, we have
\[\lim_{n\to \infty} \frac{\GR{n,p,q_{\rm{quad}}}}{n^2} = \frac{1}{2} - \lim_{n\to \infty} \frac{f^{(4)}(n;p,p/2-1)}{n^2}.\]
In particular, the limit on the left exists by~\cite{sh}.
\end{theorem}

By combining this with our results, we obtain the following new asymptotic values for the generalised Ramsey numbers at the quadratic threshold.

\begin{theorem} The following equalities hold:
\begin{eqnarray*}
\lim_{n\to \infty} \frac{\GR{n,12,62}}{n^2} = \frac{9}{22},\ \lim_{n\to \infty} \frac{\GR{n,14,86}}{n^2} = \frac{5}{12}, \ \lim_{n\to \infty} \frac{\GR{n,16,114}}{n^2} = \frac{9}{22}. \qed
\end{eqnarray*} 

\end{theorem}

Bennett, Cushman and Dudek~\cite{BCD23}*{Theorem 4} also proved that, for all $p\geq 3$, it holds that
\begin{equation}\label{eq:BCD4}
\liminf_{n\to \infty} \frac{\GR{n,p,q_{\rm{lin}}}}{n} \geq 1 - \lim_{n\to \infty} \frac{f^{(3)}(n;p,p-2)}{n^2}.
\end{equation}
Using our above results in the cases when $r=3$ and $k=p-2$ is in $\{5,6,7\}$, we get the following lower bounds at the linear threshold:
 \begin{equation}\label{eq:Lin} 
\liminf_{n\to \infty} \frac{\GR{n,7,17}}{n}\geq \frac{4}{5}, 
\ \ \liminf_{n\to \infty} \frac{\GR{n,8,23}}{n} \geq \frac{269}{330},\ \
\liminf_{n\to \infty} \frac{\GR{n,9,30}}{n}\geq \frac{4}{5}.
\end{equation}
Note that~\eqref{eq:BCD4} gives only a one-sided inequality. It happens to be tight for $p=3$ (trivially) and for $p=4$ by the result of Bennett, Cushman, Dudek and Pra\l at~\cite{BennettCushmanDudekPralat} that $\GR{n,4,5}=(\frac56+o(1))n$. However, \eqref{eq:BCD4} is not tight for $p=5$:
Gomez-Leos, Heath, Parker, Schwieder and Zerbib~\cite{GHPSZ} showed that $\GR{n,5,8}\ge \frac67(n-1)$ 
while $f^{(3)}(n;5,3)=(\frac15+o(1))n^2$, as proved in~\cite{G}. We do not know if the bounds in~\eqref{eq:Lin} are sharp.\medskip

\noindent \textbf{Organisation of the paper.} The remainder of this paper is organised as follows. Section~\ref{se:notation} introduces  some notation. An overview of our proofs can be found in Section~\ref{se:overview}.
The lower bounds are proved in Section~\ref{se:lower}, and the upper bounds are proved in Section~\ref{se:upper}. 
The proof of the upper bound of Theorem~\ref{6edge3uniform}, while using the same general proof strategy, is rather different from the other proofs in detail, so it is postponed until the end. The final section is dedicated to some concluding remarks.

\section{Notation}\label{se:notation} 
Throughout the paper, we use the following notation and definitions. Let $\I N$ denote the set of positive integers. For $m, n\in \I N$, we denote by $[n]$ the set $\{1,\dots,n\}$ and by $[m,n]$ the set $[n]\setminus [m-1]=\{m,\dots,n\}$. For a set $X$, we let 
${X\choose s}:=\{Y\subseteq X\mid |Y|=s\}$ 
be the family of all $s$-subsets of~$X$. We will often write an unordered pair $\{x,y\}$ (resp.\ triple $\{x,y,z\}$) as~$xy$ (resp.\ as~$xyz$).  
Moreover, for three real numbers $a$, $b$ and $c\geq 0$, we write $a = b\pm c$ to say that $a\in [b-c, b+c]$. Also, we write $a\gg b>0$ to mean that $b$ is a sufficiently small positive real depending on $a$.

Given an $r$-graph $G$, we denote by $V(G)$ the vertex set of $G$ and by $E(G)$ its edge set. 
Moreover, we define $|G|$ as the number of edges of $G$ and $v(G)$ as the number of vertices of~$G$. When  it is notationally convenient, we may identify an $r$-graph with its set of edges. If we specify only the edge set $E(G)$, then the vertex set is assumed to be the union of these edges, that is, $V(G):=\bigcup_{X\in E(G)} X$. 
For $r$-graphs $F$ and $H$, their \emph{union} $F\cup H$ and \emph{difference} $F\setminus H$ have edge sets respectively $E(F)\cup E(H)$ and $E(F)\setminus E(H)$ (with their vertex sets being the unions of these edges).  We reserve the lowercase letter $r$ to denote the uniformity of our hypergraphs.

For positive integers $s$ and $k$, an \emph{$(s,k)$-configuration} is an $r$-graph with $k$ edges and at most $s$ vertices, that is, an element of $\C F^{(r)}(s,k)$. An $r$-graph is called \emph{$(s,k)$-free} if it contains no $(s,k)$-configuration. Let us define another $r$-graph family
\begin{equation}\label{eq:cG}
\cG{r}{k}:=\C F^{(r)}(rk-2k+2,k) \cup\left(\bigcup_{\ell=2}^{k-1}\C F^{(r)}(r\ell-2\ell+1,\ell)\right).
\end{equation}

Thus, $\cG{r}{k}$ includes the family $\C F^{(r)}(rk-2k+2,k)$, whose Tur\'an function is the main object of study of this paper, as well as all analogous $r$-graphs for smaller sizes that are ``denser'' (that is, are subject to a stronger restriction on the number of vertices). Note that the family $\C F^{(r)}(r\ell-2\ell+1,\ell)$ that appears in the right-hand side of~\eqref{eq:cG}  for $\ell=2$ happens to be empty when $r=3$ (however, we include it to have a single formula that works for all pairs $(r,k)$). 
The family $\cG{r}{k}$ is of relevance for both the lower and the upper bounds (see Theorem~\ref{highgirth} and Lemma~\ref{chong}).

For an $r$-graph $G$, a pair $xy$ of distinct vertices (not necessarily in ${V(G)\choose 2}$) and $A\subseteq\I N\cup\{0\}$, we say that $G$ \emph{$A$-claims} the pair $xy$ if, for every $i\in A$, there are $i$ distinct edges $X_1,\dots,X_i\in E(G)$ such that $|\{x,y\}\cup (\bigcup_{j=1}^i X_j)|\leq ri-2i+2$. If $xy\in {V(G)\choose 2}$, this is the same as the existence of an $(ri-2i+2,i)$-configuration $J\subseteq G$ with $\{x,y\}\subseteq V(J)$ for every $i\in A$. When $A=\{i\}$ is a singleton, we just say \emph{$i$-claims} instead of $\{i\}$-claims.
By definition, any $r$-graph  $0$-claims any pair  (which will be notationally convenient, see e.g.\ Lemma~\ref{lm: 2}). For $i\geq 1$, let $\p{i}{G}$ be the set of all pairs in ${V(G)\choose 2}$ that are $i$-claimed by~$G$.  
For example, if $i=1$, then $\p{1}{G}$ is the usual \emph{$2$-shadow} of $G$ consisting of all pairs of vertices $uv$ such that there exists some edge $X\in E(G)$ with $u,v\in X$.
Also, let $\CI_{G}(xy)$ be the set of those $i\geq 0$ such that the pair $xy$ is $i$-claimed by $G$, that is,
\beq{eq:I}
 \textstyle
 \CI_{G}(xy):=\left\{ i\geq 0: \exists\mbox{ distinct } X_1,\dots,X_i\in E(G)\ \ \big|\{x,y\}\cup (\bigcup_{j=1}^i X_j)\big|\leq ri-2i+2\right\}.
\eeq

More generally, for disjoint subsets $A,B\subseteq \I N$, we say that $G$ \emph{\OutIn{A}{B}-claims} a pair $xy$ if $A\cap \CI_{G}(xy)=\emptyset$ and $B\subseteq \CI_{G}(xy)$. In the special case when $A=\{1\}$ and $B=\{i\}$ we just say \emph{\OutIn{1}{i}-claims}; also, we  let $\pp{1}{i}{G}:=\p{i}{G}\setminus \p{1}{G}$ denote the set of pairs in ${V(G)\choose 2}$ that are \OutIn{1}{i}-claimed by~$G$.

A {\emph{diamond}} is an $r$-graph consisting of two edges that share exactly $2$ vertices. Thus, a $(2r-3,2)$-free $r$-graph $G$ \OutIn{1}{2}-claims a pair of vertices $xy$ if and only if $xy\notin \p{1}{G}$ and there is a diamond $\{X_1,X_2\}\subseteq G$ such that $x,y\in X_1\cup X_2$.

\section{Overview of the proofs}\label{se:overview}

For the lower bounds, we combine a result from~\cite{GJKKLP} that allows us to build relatively dense $\C F^{(r)}(rk-2k+2,k)$-free $r$-graphs $G$ from a fixed $\cG{r}{k}$-free $r$-graph $F$. Namely, $G$ will be the union of many edge-disjoint copies of $F$ and, of course, the main issue is to avoid forbidden subgraphs coming from different copies of~$F$. In order to attain the desired lower bound on $|G|$, the packed copies of $F$ will be allowed to share pairs (but not triples) of vertices. Pairs inside $V(F)$ that will be allowed to be shared will be limited to those $uv$ for which $\CI_{F}(uv)$ does not contain any $i$ with $1\leq i\leq k/2$. 
This will automatically exclude forbidden subgraphs in $G$ coming from at most 2 copies of~$F$. Then, a result from~\cite{GJKKLP} will be used to eliminate any forbidden configurations whose edges come from at least 3 different copies of $F$.

If $r=3$ and $k\in\{5,7\}$, then we take for $F$ the union of many diamonds $\{x_iy_ia, x_iy_ib\}$ sharing only the pair $ab$ of vertices. This is a straightforward generalisation of the construction for $k=3$ by Glock~\cite{G}. However, if $r\ge 4$ and $k\in \{5,7\}$, then finding a suitable $F$ is a new difficult challenge, not present in~\cite{GJKKLP}. 
The initial idea that eventually led to its resolution was to take two sufficiently sparse $(r-2)$-graphs with edge sets $\{K_1^1,\dots,K_1^t\}$ and $\{K_2^1,\dots,K_2^t\}$, and let $F$ be the union of diamonds $\{\{x_i,y_i\}\cup K_1^i,\{x_i,y_i\}\cup K_2^i\}$, $1\le i\le t$, for new vertices $x_i,y_i$. Some further ideas are needed to fix the two big issues of this construction: namely, avoiding any subgraph in $\cG{r}{k}$ and having ``overhead" (like the pair $\{a,b\}$ in the construction for $r=3$) of size neglibile compared to $|F|$. We refer the reader to Section~\ref{se:7edge} for details.

For the case $(r,k)=(3,6)$, we provide an explicit construction of a 3-graph $F_{63}$ on $63$ vertices and $61$ edges, while for $r\geq 4, k=6$, the lower bound comes from the trivial construction when $F$ is a single edge.

Concerning the upper bounds, we will need the following result (proved for $r=3$ in~\cite{DP}*{Theorem~1.7} and then extended to any $r$ in
\cite{sh}*{Lemma~5}) which allows us to get rid of smaller ``denser" structures.

\begin{lemma}[\cite{sh}*{Lemma 5}]\label{chong}
    For all fixed $r\geq 3$ and $k\geq 3$,
\beq{eq:chong}
\limsup_{n\rightarrow \infty}\frac{f^{(r)}(n;rk-2k+2,k)}{n^2}\leq \limsup_{n\rightarrow \infty}\frac{\ex(n,\cG{r}{k})}{n^2}.
\eeq
\end{lemma}

Since $\C F^{(r)}(rk-2k+2,k)\subseteq \cG{r}{k}$, the opposite inequality in~\eqref{eq:chong} trivially holds. Also, the main results of \cites{DP,sh} show that both ratios in~\eqref{eq:chong} tend to a limit (which is the same for both) as $n\to\infty$. By Lemma~\ref{chong}, in order to obtain an upper bound on $f^{(r)}(n;rk-2k+2,k)$, it is enough to consider only those $r$-graphs $G$ on $[n]$ which are $\cG{r}{k}$-free.
For any such $r$-graph $G$, we define a partition of the edge set $E(G)$
by starting with the trivial partition into single edges and iteratively merging parts as long as possible using some  merging rules (that depend on $k$ and $r$). Then, we specify a set of weights that each final part (which is a subgraph $F\subseteq G$) attributes to some of the pairs in ${V(F)\choose 2}$ and use combinatorial arguments to show that every vertex pair receives total weight at most 1. Thus, the total weight assigned by the parts is at most ${n\choose 2}$ which translates into an upper bound on~$|G|$. The main difficulty lies in designing the merging and weighting rules, which have to be fine enough to detect even the extremal cases (which are quite intricate constructions) but coarse enough to be still analysable. 

The most challenging case here is $(r,k)=(3,6)$, where our solution uses a rather complicated weighting rule with values in $\{0,\frac6{61},\frac{11}{61},\frac{25}{61},\frac12,\frac{36}{61},\frac{55}{61},1\}$. Note that any weighting rule that gives the correct limit value of $
\frac{61}{330}$ has to be tight on optimal packings of the $63$-vertex configuration $F_{63}$ from the lower bound. Unfortunately, this seems to force any such rule to be rather complicated.

\section{Lower Bounds}\label{se:lower}

To prove our lower bounds, 
we use the following result,  
which is derived from~\cite{GJKKLP}.

\begin{theorem}[\cite{GJKKLP}*{Theorem~3.1}]
\label{highgirth}
Let $k\geq 2$, $r\geq 3$ and let $F$ be a $\cG{r}{k}$-free $r$-graph. Then,
$$
\liminf_{n\rightarrow \infty}n^{-2}f^{(r)}(n;rk-2k+2,k)\geq \frac{|F|}{2\,|\p{\le \lfloor k/2\rfloor}{F}|},
$$
where we define $P_{\le t}(F):=\{xy\in {V(F)\choose 2}\mid  \CI_{F}(xy)\cap [t]\not=\emptyset\}$ to consist of all pairs $xy$ of vertices of $F$ such that $\CI_{F}(xy)$ contains some $i$ with $1\leq i\leq t$.
\end{theorem}

We remark that the $\cG{r}{k}$-freeness captures the two conditions needed for $F$ in~\cite{GJKKLP}*{Theorem~3.1} and that the choice of $J:=\p{\le \lfloor k/2\rfloor}{F}$ there guarantees that $J$ contains the $2$-shadow of $F$ and that, using the notation from~\cite{GJKKLP}, the pair $(F,J)$ has non-edge girth greater than $k/2$.

To give the reader a little bit of motivation for this theorem, we briefly sketch where the ratio $\frac{|F|}{2\,|J|}$ comes from, where $J=\p{\le \lfloor k/2\rfloor}{F}$. The proof goes via  packing many edge-disjoint copies of the graph $J$ and then putting a copy of $F$ ``on top'' of each~$J$. Note that the total number of edge-disjoint copies of $J$ that we can find in $K_n$ is roughly $\binom{n}{2}/|J|$, and each copy of $F$ adds new $|F|$ edges to our $r$-graph. Hence, in total, we will have approximately $\frac{|F|}{2\,|J|}n^2$ edges, as desired. To ensure that the resulting $r$-graph remains \mbox{$(rk-2k+2,k)$}-free, the recently developed theory of conflict-free hypergraph matchings~\cites{dpfinding,GJKKL} is used, and the $\cG{r}{k}$-freeness condition of Theorem~\ref{highgirth} is necessary to apply this method.
We also remark that Theorem~\ref{highgirth} was used in~\cite{GJKKLP} to settle the case $k=4$, and by Delcourt and Postle~\cite{DP} and by Shangguan~\cite{sh} to prove the existence of the limits.

\subsection{Lower bounds in Theorems~\ref{5edge} and~\ref{7edge}}\label{se:7edge}

We apply Theorem~\ref{highgirth} to derive first the lower bounds in Theorems~\ref{5edge} and~\ref{7edge} (that is, when $k\in\{5,7\}$).
Note that if $F$ is a diamond then  
$$
\frac{|F|}{2\,|\p{1}{F}|}=\frac{2}{2(2\binom{r}{2}-1)}=\frac{1}{r^2-r-1}$$ 
 is exactly the bound we are aiming for.
The problem is that if $k\ge 4$ and we apply Theorem~\ref{highgirth} for this $F$ then $\p{\le \lfloor k/2\rfloor}{F}$ includes all pairs inside $V(F)$ and the theorem gives a weaker bound. Thus, we essentially want $F$ to consist of many edge-disjoint diamonds in such a way that the \OutIn{1}{2}-claimed pairs are ``reused'' by many different diamonds.
To illustrate our approach, we start with the simpler $3$-uniform case.

\vspace{12pt}
\bpf[Proof of the lower bounds in Theorems~\ref{5edge} and~\ref{7edge} with $r=3$]
Recall that we forbid $\C F^{(3)}(k+2,k)$ for $k=5,7$ here. 
Fix a positive integer $t$ and consider the $3$-graph $F$ consisting of $t$ diamonds $\{x_iy_ia,x_iy_ib\}$ where the $2t$ vertices $x_1, \ldots, x_t,y_1,\ldots, y_t$ are all distinct. 

Let us show that $F$ is $\cG{3}{5}$-free and $\cG{3}{7}$-free. Take any set $X\subseteq V(F)$ of size~$\ell$. If $\{a,b\}\subseteq X$ then $X$ can contain at most $\lfloor (\ell-2)/2\rfloor$ of the pairs $x_iy_i$ and thus spans at most twice as many edges in~$F$. If $X$ is disjoint from $\{a,b\}$ then $X$ spans no edges. In the remaining case $|X\cap\{a,b\}|=1$, the set $X$ contains at most $\lfloor (\ell-1)/2\rfloor$ of the pairs $x_iy_i$ and thus spans at most this many edges in $F$. Thus, for $\ell=4,5,6,7,9$, we see that $X$ spans at most $2,2,4,4,6$ edges, respectively.
Thus, $F$ is $\cG{3}{5}$-free and $\cG{3}{7}$-free, as claimed.

The above argument gives that $F$ is $(5,3)$-free and that every $(4,2)$-configuration in $F$ is $\{x_iy_ia,x_iy_ib\}$ for some $i\in [t]$. Thus, $\p{\le3}{F}\setminus \p{1}{F}$ consists only of the pair~$ab$.

As a result, Theorem~\ref{highgirth} implies that, for $k=5,7$,
$$
 \liminf_{n\rightarrow \infty}n^{-2}f^{(3)}(n;k+2,k)\geq \frac{|F|}{2\,|\p{\le3}{F}|}=\frac{2t}{2(5t+1)}.
$$
 By taking $t\to\infty$, we conclude that the lim-inf is at least $1/5$, as desired. 
\epf

Let us now informally describe how one can generalise the above
construction to higher uniformity~$r\geq 4$. We will often use the following definition.
Given an $r$-graph $G$, the \emph{girth} of $G$ is the smallest integer $\ell\geq 2$ such that there exist edges $X_1, \ldots, X_\ell$ spanning at most $(r-2)\ell+2$ vertices. For example, the girth is strictly larger than $2$ if and only if $G$ is \emph{linear} (that is, every two edges intersect in at most one vertex).
In informal discussions, we use the phrase ``high girth'' to assume that the girth is at least an appropriate constant.

Similar to the above $r=3$ case, we would like to find a
suitable $r$-graph $F$ as the union of diamonds such that the
set of $2$- or $3$-claimed pairs not in $\p{1}{F}$ is much smaller than~$|F|$. The difficulty
here is that, for example, if a pair is \OutIn{1}{2}-claimed by two diamonds, then, by the
$(4r-7,4)$-freeness requirement of Theorem~\ref{highgirth}, these two diamonds
cannot share any other vertices. In particular, we cannot simply replace $a$ and $b$ by two $(r-2)$-sets in the above construction for $r=3$.
One approach would be to consider two linear $(r-2)$-graphs $\C {K}_1$ and $\C {K}_2$
on disjoint vertex sets $A_1$ and $A_2$ 
with $|\C {K}_1|=|\C {K}_2|$ and $m:=|A_1|=|A_2|$. Let us pick a matching $M$  between some edges of $\C K_1$ and $\C K_2$, say consisting of the pairs $\{K_1^j,K_2^j\}$ for $j=1,\dots,|M|$. We add new vertices $x_j$ and $y_j$ for each $j=1,\dots,|M|$ and define $F=F(M)$ as the $r$-graph with edge set
 \beq{eq:F}
 E(F):=\left\{K_1^j\cup \{x_j,y_j\}\mid j=1,\dots,|M|\right\}\bigcup \left\{K_2^j\cup \{x_j,y_j\}\mid j=1,\dots,|M|\right\}.
\eeq 
Thus, $F$ is a
union of $|M|$ diamonds. It is easy to show that $F$ is necessarily
$(4r-7,4)$-free (see Claim~\ref{cl:small}) and that $\pp{1}{2} F$ is the union of
$\{ab: a\in K_1^j, b\in K_2^j\}$ for $1\le j\le |M|$. 
Moreover, we can additionally ensure that $\C K_1$ and $\C K_2$ have large girth, which follows from the recent results on conflict-free hypergraph matchings~\cites{dpfinding,GJKKL}.
As a direct consequence of this high-girth assumption, we can see that we do not get any forbidden configurations in $F$ when we use edges only from one ``side'' of the construction, say~$\C K_i$. Indeed, for any $\ell\geq 2$ such edges, by the high-girth assumption, the union of the corresponding $(r-2)$-sets in $\C K_i$ has more than $(r-2-2)\ell+2$ vertices, and when adding the $2\ell$ new  vertices $x_j$ and $y_j$ as above, we get more than $r\ell-2\ell+2$ vertices, that is, there is no $(r\ell-2\ell+2,\ell)$-configuration.

There are still two serious issues even for $k=5$. First, we have not guaranteed that the number of \OutIn{1}{2}-claimed pairs is much smaller than the number of edges in~$F$.
Indeed, even if
$|M|=\Theta(m^2)$ (which is the largest possible order of magnitude by  $|M|\le
{m\choose 2}/{r-2\choose 2}$), the set $\pp{1}{2} F$ may have size comparable
to $|F|=2\,|M|$ since potentially a positive fraction of pairs between $A_1$ and $A_2$ could be \OutIn{1}{2}-claimed. To ensure that $|\pp{1}{2} F|$ is much smaller than $|F|$, we form a random bipartite graph $G_3$ with parts $A_1$ and $A_2$ where every edge is included with small probability~$\alpha$.
Then, we allow $K_1^j\in \C {K}_1$ to be matched to $K_2^j\in \C {K}_2$ only if all pairs
in $K_j^1\times K_j^2$ are edges of~$G_3$. This ensures that $\pp{1}{2} F$ is a subgraph of $G_3$ and hence $|\pp{1}{2} F|\leq |G_3|\approx \alpha m^2$.

The second (much more complicated) problem is that we have to avoid dense configurations when using edges from both sides of the construction, which could overlap significantly in the ``middle layer'' formed by the vertices $x_j,y_j$.
Hence, roughly speaking, when we have two collections of $(r-2)$-sets in $\C {K}_1$ and $\C {K}_2$ that contain few vertices, we want to avoid matching many of these $(r-2)$-sets with each other to form diamonds. 
Formally, we construct an auxiliary bipartite graph $H$ where $K_1\in \C {K}_1$ and $K_2\in \C {K}_2$ are adjacent if $\{\{a,b\}: a\in K_1, b\in K_2\}$ is a subset of $E(G_3)$ (so a diamond could be attached to $K_1, K_2$), and we define a family $\C {C}$ of sets of disjoint edges of $H$ such that if the matching $M$ avoids $\C C$, then $F(M)$ avoids all forbidden configurations.
Then, the goal is to find a large matching in $H$ which avoids each of the problematic configurations in~$\C {C}$. 
However, this would still not be possible with the current construction.
Roughly speaking, the problem is that, for every $u_1\in K_1^i\in \C {K}_1$ and $u_2\in K_2^i\in \C {K}_2$, 
being able to attach a diamond to $K_1^i,K_2^i$ implies that $u_1u_2\in E(G_3)$. 
However, the presence of $u_1u_2$ in $G_3$ increases significantly the probability that, 
for any fixed pair $K_1^j\in \C {K}_1$ and $K_2^j\in \C {K}_2$ such that $K_1^i\cap K_1^j = \{u_1\}$ and $K_2^i\cap K_2^j = \{u_2\}$,
a diamond can be attached to the pair $K_1^i, K_2^i$.
As it turns out, the number of such pairs $(K_1^j, K_2^j)$ happens to be too large.
To fix this, we first randomly sparsify the complete graphs on $A_1$ and $A_2$ with a well-chosen probability, and then restrict our attention to $(r-2)$-sets in $A_1, A_2$ which form cliques in the underlying random graphs. This allows better control on the number of edges $K_1^j\in \C {K}_1$ and $K_2^j\in \C {K}_2$ containing a pair $u_1u_2\in E(G_3)$.

We will use the following  concentration inequality for functions of independent coordinates satisfying a Lipschitz condition, which is known as \emph{the Bounded Difference Inequality} or  \emph{McDiarmid's inequality}; it can be also derived from the Azuma-Hoeffding Inequality.

\begin{lemma}[\cite{mcdiarmid}*{Lemma~1.2}]\label{lem:BDI}
Let $X_1,X_2,\dots,X_n$ be independent random variables with $X_i$ taking values in $\Lambda_i$, and let $f:\Lambda_1\times \dots \times \Lambda_n \to \mathbb R$ be a function that satisfies the following Lipschitz condition for some numbers $(c_i)_{i=1}^n$: for every $i\in [n]$ and every two vectors $x, \tilde{x}\in \Lambda_1\times \dots \times \Lambda_n$ that differ only in the $i$-th coordinate, it holds that $|f(x)-f(\tilde{x})|\leq c_i$.

Then, the random variable $Z := f(X_1, \dots, X_n)$ satisfies
\begin{align*}
    \mathbb P\left[\,|Z-\E{Z}|>s\,\right]\leq 2\exp\left(-\dfrac{2s^2}{\sum_{i=1}^n c_i^2}\right).
\end{align*}

\end{lemma}

Let us also recall the classical \emph{Chernoff bound} stating that, for every binomial random variable $X$ and every $t\ge 0$, 
\[\mathbb P\left[\,|X-\mathbb E[X]|\ge t\,\right]\le 2\exp\left(-\tfrac{t^2}{2(\mathbb E[X]+t/3)}\right),\]
see e.g.\ \cite{JLR00}*{Theorem 2.1}.

Furthermore, we will need a simplified version of a result from~\cite{GJKKL} on the existence of approximate clique packings of high girth.

\begin{theorem}[\cite{GJKKL}*{Theorem 1.4}]\label{packing}
    For all $c_0>0$, $\ell\geq 2$ and $r\geq 3$, there exists $\varepsilon_0>0$ such that, for all $\varepsilon\in (0,\varepsilon_0)$, there exists $m_0$ such that the following holds for all $m\geq m_0$ and $c\geq c_0$. Let $G$ be a graph on $m$ vertices such that every edge of $G$ is contained in $(1\pm m^{-\varepsilon})cm^{r-2}$ cliques of order $r$.
    Then, there exists a $K_r$-packing $\C K $ in $G$ of size $|\C K|\ge (1-m^{-\varepsilon^3})|G|/\binom{r}{2}$ such that, for every $j\in [2,\ell]$, any set of $j$ elements in $\C K$ spans more than $(r-2)j+2$ vertices.
\end{theorem}

Note that if we consider the packing $\C K$ returned by Theorem~\ref{packing} as an $r$-graph, then the last requirement is precisely that the girth of $\C K$ is larger than $\ell$.   

Now we are ready to provide the construction which establishes the lower bounds in Theorems~\ref{5edge} and~\ref{7edge} for $r\ge 4$. 

\begin{lemma}\label{57cons}
Fix any integer $r\geq 4$. Then, for a sufficiently small real $\alpha>0$ and a sufficiently large integer~$m$, that is, for $1/r\gg \alpha\gg 1/m$,
there exists an $r$-graph $F$ satisfying each of the following properties:
\begin{enumerate}[(a)]   
\item\label{it:57Free} $F$ is $(5r-8,5)$-free and $(7r-12,7)$-free,
\item\label{it:57Free'} $F$ is $(2r-3,2)$-free, $(3r-5,3)$-free, $(4r-7,4)$-free and $(6r-11,6)$-free,
\item\label{it:|F|} $|F| = \Omega(\alpha^{3/4} m^2)$,
\item\label{it:J} $|\p{\le3}{F}|\le  \frac{r^2-r-1}{2}\, |F|+ 2\alpha m^2$.
\end{enumerate}
\end{lemma}
\bpf
Let $A_1$ be a set of size $m$. Sample every edge of the complete graph on $A_1$ independently with probability $\beta:=\alpha^{3/4}$ to get a random graph $G_1$ on~$A_1$. 
In the sequel, implicit constants in the $O,\Omega,\Theta$-notation may depend on $r$ but not on $\alpha$ and $m$ unless the dependence is explicitly indicated in a lower index such as $O_{\alpha}$. 
We say that an event holds \emph{with high probability} if its probability tends to $1$ as $m\to \infty$.

\begin{claim}
\label{cl:cliques} 
With high probability, $G_1$ satisfies the following properties.
\begin{enumerate}[(i)]
\item\label{it:degree} For every vertex $v\in V(G_1)$, we have $d(v)=\Theta(\beta m).$
\item\label{it:codegree} For any pair of vertices $u,v\in V(G_1)$, we have $|N(u)\cap N(v)|=\Theta(\beta^2m).$
\item\label{it:CliquesOnEdge} If $r\ge 5$, then every edge in $G_1$ is contained in $(1\pm m^{-1/3}) c m^{r-4}$ cliques of size $r-2$, where $c := \beta^{{r-2\choose 2}-1}/(r-4)!$.
\item\label{it:cliques} There is an $(r-2)$-graph $\C K_1$ of girth at least $8$, with vertex set $A_1$ and edge set being a collection of edge-disjoint $(r-2)$-cliques in $G_1$ such that all but $o(m^2)$ edges of $G_1$ belong to a clique in $\C K_1$.
\end{enumerate}
\proofthmfalse
\end{claim}
\bpf[Proof of Claim~\ref{cl:cliques}]
The first two properties follow easily by noting that for a vertex $v$ (resp.\ a pair $u,v$) the probability of failure is $\me^{-\Omega_\alpha(m)}$ by Chernoff's bound and 
then taking the union bound  over all (polynomially many in $m$) choices. 

Let us turn to the third claim. Fix an edge $uv\in G_1$ (that is, we condition on $uv$ being sampled). Let $X$ be the number of $(r-2)$-cliques in $G_1$ containing~$uv$. Since each potential clique containing $uv$ has $\tbinom{r-2}{2}-1$ edges other than $uv$, each of which appears independently with probability $\beta$, we have 
$$
\E{X}=\tbinom{m-2}{r-4} \beta^{{r-2\choose 2}-1} = c m^{r-4}+O_\alpha(m^{r-5}).
$$ 

Let us show that $X$ is concentrated. 
We use the Bounded Difference Inequality (Lemma~\ref{lem:BDI}). Altering the state of a pair with one endvertex among $u,v$ may change the value of $X$ by at most $m^{r-5}$, and every other edge may change the value of $X$ by at most $m^{r-6}$. Thus, by Lemma~\ref{lem:BDI},
we have
\begin{align*}
\mathbb P\left[\,X\neq (1\pm m^{-1/3}) c m^{r-4}\,\right]&\le 
\mathbb P\left[\,|X - \E{X}|\geq  cm^{r-4-1/3}/2\,\right]\\
&\le\; 2\exp\left(-\frac{c^2m^{2(r-4-1/3)}/2}{2m\cdot m^{2(r-5)} + m^2\cdot m^{2(m-6)}}\right)\ =\ o(m^{-2}).
\end{align*}
(If $r=5$, then $X$ is the number of triangles, which is binomially distributed, and Chernoff's bound can be applied instead of Lemma~\ref{lem:BDI}.)
The union bound over all ${m\choose 2}$ choices of $uv$ finishes the proof.

Let us turn to the existence of~$\C K_1$. Note first that the case $r=4$ is trivial since we can take each edge of $G_1$ as a clique of order~$2$ (and any $2$-graph has infinite girth according to our definition of girth), so we assume that $r\geq 5$.
Theorem~\ref{packing} for $c_0:=c$, $\ell:=8$ and $r-2$ returns some $\varepsilon_0$. Let $m_0$  be the value returned by the theorem for $\varepsilon:= \min(1/3,\varepsilon_0/2)$. Since $m_0$ depends only on $r$ and $\alpha$, we can assume that $m>m_0$. Thus, Theorem~\ref{packing} applies to  any graph $G_1$ satisfying Property~\ref{it:CliquesOnEdge} and produces $\C K_1$ with all the stated properties.
\epf

We now fix a graph $G_1$ on the set $A_1$ and an $(r-2)$-graph $\C K_1$ satisfying Claim~\ref{cl:cliques}. Let $(A_2,G_2, \C {K}_2)$ be a disjoint copy of $(A_1,G_1,\C K_1)$.  We identify $\C K_1$ and $\C K_2$ with their edge sets. 
Let $G_3$ be a random bipartite graph with parts $A_1$ and $A_2$ where every edge between $A_1$ and $A_2$ is sampled independently with probability $\alpha$.

We set $t:=|\C {K}_1|=|\C {K}_2|$ and note that $t=\Theta(\beta m^2)$.
We also define an auxiliary bipartite graph $H$ with parts $\C {K}_1$ and $\C {K}_2$ where  $K_1\in \C {K}_1$ and $K_2\in \C {K}_2$ are adjacent in $H$ if each of the $(r-2)^2$ pairs between $K_1\subseteq A_1$ and $K_2\subseteq A_2$ is an edge of~$G_3$. 
In particular, a pair in $\C {K}_1\times \C {K}_2$ is an edge of $H$ with probability $\alpha^{(r-2)^2}$. Define 
$$
d := \alpha^{(r-2)^2}t = \Theta(\beta \alpha^{(r-2)^2} m^2).
$$

\begin{claim}\label{cl:reg}
With high probability, all vertices in the graph $H$ have degree $(1\pm m^{-1/3}) d$.
\end{claim}
\bpf[Proof of Claim~\ref{cl:reg}]
Take any vertex $K$ of $H$, that is, an edge of $\C K_i$ for $i=1$ or $2$. Since $\C K_1$ and $\C K_2$ are fixed, the degree of $K$ in $H$ is a function of the $(r-2)m$ independent Bernoulli variables that encode the edges of $G_3$ between $K\subseteq A_i$ and the opposite side $A_{3-i}$. 
Furthermore, one edge in $G_3$ can influence the appearance of at most $\tfrac{m-1}{r-3}\leq m$ edges in $H$ containing~$K$, since every pair of vertices in $A_{3-i}$ is contained in at most one clique.

Thus, denoting the degree of $K$ in $H$ by $\mathrm{deg}_H(K)$ and using that its expectation is $d$, the Bounded Difference Inequality (Lemma~\ref{lem:BDI}) implies that
\[
\mathbb P\left[|\mathrm{deg}_H(K) - d| \geq \frac{d}{m^{1/3}}\right]\leq 2\exp\left(-\frac{2\alpha^{2(r-2)^2} t^2/m^{2/3}}{(r-2)m\cdot m^2}\right) = \exp(-\Omega_{\alpha}(m^{1/3})).\]
A union bound over all $O(m^2)$ edges in $\C K_1\cup \C K_2$ proves the claim.
\epf

Our goal will be to find a matching $M$ in $H$ with size $\Omega(t)=\Omega(\beta m^2)$ and certain additional properties. Given $M$, we define the $r$-graph $F=F(M)$ as in~\eqref{eq:F}. Namely, we start with the edgeless hypergraph on $A_1\cup A_2$ and, for every edge $K_1^jK_2^j$ in the matching $M$, add to $F$ new vertices $x_j,y_j$ and new edges $K_1^j\cup \{x_j, y_j\}$ and $K_2^j\cup \{x_j, y_j\}$ (forming a diamond). 

We remark that, at this point of the proof, we do not fix the choice of $G_3$ yet but view $G_3$ (and hence $H$) as a random graph since we still want to bound the number of certain problematic subconfigurations.

We can rule out some small configurations in $F$ without any additional assumptions on $M$. 

\begin{claim}\label{cl:small}
The $r$-graph $F=F(M)$ is $(4r-7,4)$-free, $(3r-4,3)$-free  
and $(2r-3,2)$-free.
\end{claim}
\bpf[Proof of Claim~\ref{cl:small}] Every triple of edges $X_1, X_2, X_3\in F$ satisfies that each of $|X_1\cap X_2|$, $|X_2\cap X_3|$ and $|X_3\cap X_1|$ is at most 1, or contains a diamond $(X_i, X_j)$ for some $i\neq j$ with the third edge intersecting $X_i\cup X_j$ in at most one vertex; in either case, $|X_1\cup X_2\cup X_2|\ge 3r-3$.
In particular, this gives that $F$ is $(2r-3, 2)$-free and $(3r-4,3)$-free. 

Now, suppose on the contrary that $F$ has a $(4r-7,4)$-configuration, say coming from some $(r-2)$-graphs $F_1\subseteq \C K_1$ and $F_2\subseteq \C K_2$. For $i=1,2$, let $e_i:=|F_i|$ and let 
$$
d_i:=(r-2)|F_i|-|\cup_{X\in F_i} X|,
$$ 
calling it the \emph{defect} of $F_i$.
Thus, $e_1+e_2=4$. By symmetry, assume that $e_1\ge e_2$. By $e_1\in [2,4]$ and the high-girth assumption, we obtain that 
\beq{eq:d1}
 d_1\le e_1(r-2)-(e_1(r-2-2)+3)=2e_1-3.
\eeq 

Let $d'\le e_2$ be the number of pairs in $M$ between $F_1$ and $F_2$  (which is exactly the number of diamonds in the hypothetical $(4r-7,4)$-configuration in $F$ that we started with). Thus we have
\beq{eq:d'}
 4r-7 \ge \big(e_1(r-2)-d_1\big)+\big(e_2(r-2)-d_2\big)+2(4-d')=4r-d_1-d_2-2d',
 \eeq
 that is, $d_1+d_2\ge 7-2d'$. 

Now, it is routine to derive a contradiction. Although some of the following cases can be combined together, we prefer (here and later) to treat each possible value of $e_1$ as a separate case for clarity. 
If $e_1=4$, then $e_2=0$, $d_2=0$, $d'=0$, and thus $d_1\ge 7$. If $e_1=3$, then $e_2=1$, $d_2=0$, $d'\le 1$ and thus $d_1\ge 5$. If $e_1=2$, then $e_2=2$, $d_2\le 1$ and $d'\le 2$ giving $d_1\ge 2$. However, the obtained lower bound on $d_1$ contradicts~\eqref{eq:d1}  in each of the three possible cases. Thus, $F$ contains no $(4r-7, 4)$-configuration, as desired.\epf

To ensure that $F$ is $(5r-8,5)$-free, $(6r-11,6)$-free and $(7r-12,7)$-free, we have to construct the matching $M$ a bit more carefully. In the following, we will define some problematic configurations and show that there are only few of them with high probability. 
We will then be able to use a probabilistic argument to construct a matching $M$ that avoids all problematic configurations.

Let us introduce some further terminology. For an $(r-2)$-graph $\C K$ and integers $e',d'$, let $\CS{e'}{d'}{\C K}$ be the family of all subgraphs of $\C K$ with $e'$ edges and defect $d'$. (Recall that this means that the union of these $e'$ edges has exactly $e'(r-2)-d'$ vertices.)
We refer the reader to Figure~\ref{fg:CS} for some special cases of this definition that will play important role in our proof.
\begin{figure}
\centering
    \begin{scriptsize}
    \begin{tabular}[b]{|>{\centering}m{3cm}|>{\centering}m{3cm}|>{\centering\arraybackslash}m{3cm}|}    
\hline
    {
    \begin{tikzpicture}[scale=0.6,line cap=round,line join=round,x=1cm,y=1cm,baseline={(0,-1)}]

    \draw [rotate around={33.690067525979785:(-3,2)},line width=0.5pt] (-3,2) ellipse (1.5cm and 0.7cm);
    \draw [rotate around={-33.690067525979785:(-3,0.7)},line width=0.5pt] (-3,0.7) ellipse (1.5cm and 0.7cm);
    \draw[color=black] (-3,-2.2) node {\large{$\CS{2}{1}{\C K_i}$}};
    \end{tikzpicture}
    }  
    & 
    {
    \begin{tikzpicture}[scale=0.6,line cap=round,line join=round,x=1cm,y=1cm,baseline={(0,-0.5)}]
    \draw [rotate around={33.690067525979785:(-6,3)},line width=0.5pt] (-6,3) ellipse (1.5cm and 0.7cm);
    \draw [rotate around={-33.690067525979785:(-6,1.7)},line width=0.5pt] (-6,1.7) ellipse (1.5cm and 0.7cm);
    \draw [rotate around={90:(-3.7,2.3)},line width=0.5pt] (-3.7,2.3) ellipse (1.5cm and 0.7cm);
    \draw[color=black] (-5,-1.3) node {\large{$\CS{3}{1}{\C K_i}$}};
    \end{tikzpicture}
    } 
    &
    {
    \begin{tikzpicture}[scale=0.6,line cap=round,line join=round,x=1cm,y=1cm,baseline={(0,3)}]
    \draw [rotate around={90:(-5.5,2.2)},line width=0.5pt] (-5.5,2.2) ellipse (1.5cm and 0.7cm);
    \draw [rotate around={45:(-4.7,1.9)},line width=0.5pt] (-4.7,1.9) ellipse (1.5cm and 0.7cm);
    \draw [rotate around={-45:(-6.2,1.9)},line width=0.5pt] (-6.2,1.9) ellipse (1.5cm and 0.7cm);
    \draw[color=black] (-5.3,-1.4) node {\large{$\CS{3}{2}{\C K_i}$}};
    \end{tikzpicture}
    }
    
\\
\hline
    {
    \begin{tikzpicture}[scale=0.6,line cap=round,line join=round,x=1cm,y=1cm,baseline={(0,-1)}]
    \draw [rotate around={33.690067525979785:(-3,2.5)},line width=0.5pt] (-3,2.5) ellipse (1.5cm and 0.7cm);
    \draw [rotate around={-33.690067525979785:(-3,1.2)},line width=0.5pt] (-3,1.2) ellipse (1.5cm and 0.7cm);
    \draw [rotate around={33.690067525979785:(-3,0)},line width=0.5pt] (-3,0) ellipse (1.5cm and 0.7cm);
    \draw[color=black] (-3,-2) node {\large{$\CS{3}{2}{\C K_i}$}};
    \end{tikzpicture}
    }   
    & 
   {
   \begin{tikzpicture}[scale=0.6,line cap=round,line join=round,x=1cm,y=1cm,baseline={(0,0)}]
    \draw [rotate around={33.690067525979785:(-3,3.5)},line width=0.5pt] (-3,3.5) ellipse (1.5cm and 0.7cm);
    \draw [rotate around={-33.690067525979785:(-3,2.2)},line width=0.5pt] (-3,2.2) ellipse (1.5cm and 0.7cm);
    \draw [rotate around={90:(-2,2.8)},line width=0.5pt] (-2,2.8) ellipse (1.6cm and 0.7cm);
    \draw[color=black] (-2.7,-0.7) node {\large{$\CS{3}{3}{\C K_i}$}};
    \end{tikzpicture}
    } 
    &
    {
    \begin{tikzpicture}[scale=0.7,line cap=round,line join=round,x=1cm,y=1cm,baseline={(0,2.5)}]
    \draw [rotate around={90:(-7,0.6)},line width=0.5pt] (-7,0.6) ellipse (0.8cm and 0.6cm);
    \draw [rotate around={90:(-7,1.83)},line width=0.5pt] (-7,1.83) ellipse (0.8cm and 0.6cm);
    \draw [rotate around={-25:(-8,0.45)},line width=0.5pt] (-8,0.45) ellipse (1.3cm and 0.4cm);
    \draw [rotate around={25:(-6,0.45)},line width=0.5pt] (-6,0.45) ellipse (1.3cm and 0.4cm);
    \draw[color=black] (-6.8,-2) node {\large{$\CS{4}{3}{\C K_i}$}};
    \end{tikzpicture}
    }

\\
\hline
    {
    \begin{tikzpicture}[scale=0.6,line cap=round,line join=round,x=1cm,y=1cm,baseline={(0,-1)}]
    \draw [line width=0.5pt] (-7.8,2.8) circle (1cm);
    \draw [line width=0.5pt] (-7.8,0.2) circle (1cm);
    \draw [rotate around={90:(-6,1.5)},line width=0.5pt] (-6,1.5) ellipse (1cm and 1cm);
    \draw [rotate around={90:(-7.2,1.5)},line width=0.5pt] (-7.2,1.5) ellipse (1cm and 1cm);
    \draw[color=black] (-6.85,-2.1) node {\large{$\CS{4}{3}{\C K_i}$}};
    \end{tikzpicture}
    }   
&
   {
   \begin{tikzpicture}[scale=0.6,line cap=round,line join=round,x=1cm,y=1cm,baseline={(0,-1)}]
    \draw [rotate around={33.690067525979785:(-7.1,2.3)},line width=0.5pt] (-7.1,2.3) ellipse (1.5cm and 0.7cm);
    \draw [rotate around={-33.690067525979785:(-7.1,1)},line width=0.5pt] (-7.1,1) ellipse (1.5cm and 0.7cm);
    \draw [rotate around={90:(-6.1,1.6)},line width=0.5pt] (-6.1,1.6) ellipse (1.6cm and 0.7cm);
    \draw [rotate around={90:(-4.5,1.6)},line width=0.5pt] (-4.5,1.6) ellipse (1.6cm and 0.7cm);
    \draw[color=black] (-6.1,-2.4) node {\large{$\CS{4}{3}{\C K_i}$}};
    \end{tikzpicture}
    } 
    &
    {
    \begin{tikzpicture}[scale=0.6,line cap=round,line join=round,x=1cm,y=1cm,baseline={(0,3.2)}]
    \draw [rotate around={0:(-7,2)},line width=0.5pt] (-7,2) ellipse (1.4142135623730974cm and 1cm);
    \draw [rotate around={-45:(-7.75,0.75)},line width=0.5pt] (-7.75,0.75) ellipse (1.2905694150420877cm and 0.7352342586156395cm);
    \draw [rotate around={45:(-6.25,0.75)},line width=0.5pt] (-6.25,0.75) ellipse (1.2905694150420928cm and 0.7352342586156425cm);
    \draw [line width=0.5pt] (-7,3.4) circle (0.9cm);
    \draw[color=black] (-6.8,-1.8) node {\large{$\CS{4}{4}{\C K_i}$}};
    \end{tikzpicture}
    }

\\
\hline
    {
    \begin{tikzpicture}[scale=0.6,line cap=round,line join=round,x=1cm,y=1cm,baseline={(0,-0.7)}]
    \draw [rotate around={90:(-8,1)},line width=0.5pt] (-8,1) ellipse (1.2807764064044187cm and 0.8002425902201228cm);
    \draw [rotate around={90:(-6,1)},line width=0.5pt] (-6,1) ellipse (1.2807764064044058cm and 0.8002425902201146cm);
    \draw [rotate around={0:(-7,0)},line width=0.5pt] (-7,0) ellipse (1.2807764064044145cm and 0.80024259022012cm);
    \draw [rotate around={0:(-7,2)},line width=0.5pt] (-7,2) ellipse (1.2807764064044145cm and 0.80024259022012cm);
    \draw[color=black] (-7,-2.3)node {\large{$\CS{4}{4}{\C K_i}$}};
    \end{tikzpicture}
    }   
    & 
   {\begin{tikzpicture}[scale=0.6,line cap=round,line join=round,x=1cm,y=1cm,baseline={(0,-0.5)}]
     \draw [rotate around={0:(-6.5,0.3)},line width=0.5pt] (-6.5,0.3) ellipse (1.5cm and 0.3cm);
    \draw [rotate around={-60:(-5.8,1.5)},line width=0.5pt] (-5.8,1.5) ellipse (1.5cm and 0.3cm);
    \draw [rotate around={60:(-7.2,1.5)},line width=0.5pt] (-7.2,1.5) ellipse (1.5cm and 0.3cm);
    \draw [rotate around={90:(-6.5,1.4)},line width=0.5pt] (-6.5,1.4) ellipse (1.36cm and 0.3cm);
    \draw[color=black] (-6.4,-2.1) node {\large{$\CS{4}{5}{\C K_i}$}};
    \end{tikzpicture}
    } 
    &
    {
    \begin{tikzpicture}[scale=0.6,line cap=round,line join=round,x=1cm,y=1cm,baseline={(0,2)}]
    \draw [rotate around={0:(-6.5,-1.1)},line width=0.5pt] (-6.5,-1.1) ellipse (1.5cm and 0.3cm);
    \draw [rotate around={-60:(-5.8,0.1)},line width=0.5pt] (-5.8,0.1) ellipse (1.5cm and 0.3cm);
    \draw [rotate around={60:(-7.2,0.1)},line width=0.5pt] (-7.2,0.1) ellipse (1.5cm and 0.3cm);
    \draw [rotate around={0:(-6.5,-0.1)},line width=0.5pt] (-6.5,-0.1) ellipse (1.5cm and 0.3cm);
    \draw[color=black] (-6.4,-3.45) node {\large{$\CS{4}{5}{\C K_i}$}};
    \end{tikzpicture}
    }

\\
\hline

\end{tabular}
\end{scriptsize}
\caption{
Examples of $\mathcal{S}_{e',d'}$-subgraphs (that is, having size $e'$ and defect $d'$) in the high-girth $(r-2)$-graph $\C K_i$ for some pairs $(e',d')$.\label{fg:CS} Since the hypergraph $\C K_i$ is linear, each drawn  intersection has size~1. For $(e',d')$ in $\{(2,1),(3,1),(3,3)\}$, the family $\mathcal{S}_{e',d'}(\C K_i)$ consists of a unique $(r-2)$-graph up to isomorphism. For $(e',d')$ in $\{(3,2),(4,5)\}$, there are exactly two non-isomorphic examples. For the remaining pairs $(e',d')$, we provide a non-exhaustive list.
}
\end{figure}
For integers $e_1,e_2,d_1,d_2$ with $e_1\ge e_2$, let the set $\CC{e_1}{d_1}{e_2}{d_2}$ consist of all matchings 
\beq{eq:Ks}
\left\{\{K_i^j,K_{3-i}^j\}\mid  j\in [e_2]\right\}
\eeq
 in $H$ for some $i=1$ or $2$ such that the $(r-2)$-graph $\{K_{3-i}^{1},\dots,K_{3-i}^{e_2}\}$ is in $\CS{e_2}{d_2}{\C K_{3-i}}$ while $\{K_i^1,\dots,K_i^{e_2}\}$ extends to an element of $\CS{e_1}{d_1}{\C K_i}$. (Recall that the set in~\eqref{eq:Ks} is a matching in $H$ if, for $i=1,2$, the $(r-2)$-sets $K_i^1,\dots,K_i^{e_2}$ are pairwise distinct
and, for every $j\in [e_2]$, all pairs in $K_1^j\times K_2^j$ are edges of~$G_3$.)
 In other words, the set $\CC{e_1}{d_1}{e_2}{d_2}$ can be constructed as follows. 
 Pick a subgraph of size $e_1$ and defect $d_1$ in one of $\C K_1$ and $\C K_2$, then a subgraph of size $e_2$ and defect $d_2$ in the other one; viewing these as two sets of vertices in the bipartite graph $H$, add to $\CC{e_1}{d_1}{e_2}{d_2}$ all matchings in $H$ of size $e_2$ (that is, fully pairing the smaller subgraph). Note that the definition of $\CC{e_1}{d_1}{e_2}{d_2}$ is symmetric in $\C K_1$ and $\C K_2$.

The final family of conflicts that our matching $M$ will have to avoid will consists of four families of the type $\CC{e_1}{d_1}{e_2}{d_2}$ plus a related family. In order to illustrate the general proof, we will show first that forbidding $\CC{3}{3}{2}{1}$ alone takes care of all undesired configurations in $F=F(M)$ with at most $6$ edges. 
Note that, by the high-girth assumption, every pair in $\CC{3}{3}{2}{1}$ is a matching of a ``loose $2$-path"  in one of $\C K_1$ and $\C K_2$ into a ``loose $3$-cycle" in the other. 

\begin{claim}\label{cl:C3321}
 If $M$ is a matching in $H$ that avoids $\CC{3}{3}{2}{1}$, then the $r$-graph $F=F(M)$, as defined in~\eqref{eq:F}, is $(5r-8,5)$-free and $(6r-11,6)$-free.
\end{claim} 
\bpf[Proof of Claim~\ref{cl:C3321}]
Suppose on the contrary that we have a $(5r-8,5)$-configuration in $F$. For $i=1,2$, let $F_i$ be the $(r-2)$-graph consisting of the edges of $\C K_i$ involved in the configuration; also, let $e_i:=|F_i|$ be the size and $d_i:=e_i(r-2)-v(F_i)$ be the  defect of $F_i$. Thus, $e_1+e_2=5$.

Let $d'$ be the number of edges in $M$ between $F_1$ and $F_2$. By a calculation analogous to~\eqref{eq:d'}, we have that
 \beq{eq:C3321}
 d_1+d_2\ge 8-2d'.
 \eeq

Without loss of generality, assume that $e_1\ge e_2$. Of course, $d'\le e_2$. Since $3\le e_1\le 5$, we have by the large girth assumption that~\eqref{eq:d1} holds, that is, $d_1\le 2e_1-3$.

If, $e_1=5$ then $e_2=0$ and so $d_2=d'=0$, but then~\eqref{eq:d1} and \eqref{eq:C3321}  give a contradiction $8\le d_1\le 7$.  If, $e_1=4$ then $e_2=1$ and so $d_2=0$ and $d'\le 1$; however, then~\eqref{eq:d1} and \eqref{eq:C3321} give $6\le d_1\le 5$, a contradiction again. Finally, suppose that $e_1=3$. We have that $e_2=2$ and thus, $d_2\le 1$ and $d'\le 2$. Now,~\eqref{eq:d1} and \eqref{eq:C3321}  give that $3\le d_1\le 3$. Thus, we have equalities everywhere; in particular, $d_1=3$, $d_2=1$ and $d'=2$. We see that the pair of edges of $M$ coming from the two diamonds in $F$ belongs to $\CC{3}{3}{2}{1}$, a contradiction.

Now, suppose on the contrary that $F$ contains a $(6r-11,6)$ configuration. For $i=1,2$, let $F_i$ be the $(r-2)$-subgraph of $\C K_i$ involved in it, with $e_i$ and $d_i$ denoting its size and defect.
Without loss of generality, assume that $e_1\ge e_2$. Let $d'$ be the number of $M$-edges between $F_1$ and~$F_2$. By a version of~\eqref{eq:d'}, we get that
\beq{eq:C33212a}
 d_1+d_2\ge 11-2d',
 \eeq
while the inequality in~\eqref{eq:d1} remains unchanged.

If $e_1=6$, then $e_2=d_2=d'=0$, so~\eqref{eq:d1} and~\eqref{eq:C33212a} give that $11\le d_1\le 9$, a contradiction. If $e_1=5$, then $e_2=1$, $d_2=0$, $d'\le 1$, and we get
that $9\le d_1\le 7$, a contradiction. If $e_1=4$, then $e_2=2$, $d_2\le 1$, $d'\le 2$, and
we get that $6\le d_1\le 5$, a contradiction. 

Finally, it remains to check the (slightly less straightforward) case when $e_1=3$. Here, we have $e_2=3$, $d_2\le 3$
and $d'\le 3$. Thus,~\eqref{eq:d1} and~\eqref{eq:C33212a} give that $11-2d'\le d_1+d_2\le 6$. We conclude that $d'=3$ and thus, $M$ fully matches $F_1$ and~$F_2$.
By symmetry, suppose that $d_1\ge d_2$. Thus, $d_1=3$. Since $d_2\ge 5-d_1\ge 2$ is positive, we can find $K_2^1,K_2^2\in F_2$ that intersect. The $(r-2)$-graphs $F_1\in\CS{3}{3}{{\C K_1}}$ and $\{K_2^1,K_2^2\}\in \CS{2}{1}{{\C K_2}}$ show that the edges of $M$ containing $K_2^1$ and $K_2^2$ form a pair in $\CC{3}{3}{2}{1}$, a contradiction.
\epf

Next, we define the ``exceptional'' family. First, for $i=1,2$, let $\CCS{\C K_i}$ consist of those $(r-2)$-graphs in $\CS{4}{3}{\C K_i}$ that do not contain an isolated edge (equivalently, not containing an $\CS{3}{3}{}$-subgraph). It is easy to show (see Claim~\ref{cl:Struc}) that $\CCS{\C K_i}$ consists precisely of subgraphs of $\C K_i$ whose edge set can be ordered as $\{X_1,\dots,X_4\}$ such that $|X_i\cap (\bigcup_{j=1}^{i-1} X_j)|=1$ for each $i\in[2,4]$, that is,  it consists of ``loose subtrees" with 4 edges.
Also, let the conflict family $\CCP$ be obtained by taking full $H$-matchings of some element of $\CS{3}{3}{\C K_i}$ into $\CCS{\C K_{3-i}}$ for $i=1,2$. Clearly, $\CCP\subseteq \CC{4}{3}{3}{3}$. 

The final family $\C C$ of the conflicts that we are going to use is
 \beq{eq:CC}
 \C C:=\CC{3}{3}{2}{1}\cup \CC{4}{4}{3}{2}\cup \CC{4}{5}{2}{1}\cup \CC{5}{7}{2}{1}\cup \CCP.
\eeq

Basically, our definition of the conflict family $\C C$ is motivated by the proof of Claim~\ref{cl:7Free} below: for each still possible way of having a $(7r-12,7)$-configuration in $F(M)$, we add further conflicts  that rule it out. We do not try to take minimal possible conflict families, but rather the ones that are easy to describe.
Note that we cannot take the full family $\CC{4}{3}{3}{3}$  as the upper bound of Claim~\ref{cl:Expectations} on the expected number of conflicts fails for it; fortunately, its smaller subfamily $\CCP$ suffices.

\begin{claim}\label{cl:7Free} If a matching $M$ in $H$ does not contain any element of $\C C$ as a subset, then the $r$-graph $F=F(M)$ defined by~\eqref{eq:F} contains no $(7r-12,7)$-configuration.\end{claim}
\bpf[Proof of Claim~\ref{cl:7Free}] Suppose on the contrary that $F$ contains a $(7r-12,7)$-configuration. For $i=1,2$, let $F_i$ be the $(r-2)$-subgraph of $\C K_i$ involved in it, with $e_i$ and $d_i$ denoting its size and defect. Thus, $e_1+e_2=7$.
Without loss of generality, assume that $e_1\ge e_2$. Let $d'$ be the number of the diamonds involved.
By a calculation analogous to~\eqref{eq:d'}, we have
 \beq{eq:7Count}
 d_1+d_2\ge 12-2d'.
 \eeq

First, we rule out the easy cases when $e_1\ge 6$. If $e_1=7$, then $e_2=d_2=d'=0$ but~\eqref{eq:d1} and~\eqref{eq:7Count} give that $12\le d_1\le 11$, a contradiction. If $e_1=6$, then $e_2=1$, $d_2=0$ and $d_1'\le 1$, so we get $10\le d_1\le 9$, a contradiction again. 

Thus, $e_1\le 5$. Since $e_2\ge 2$, the large girth assumption gives similarly to~\eqref{eq:d1} that
 \begin{equation}\label{eq:d2}
 d_2\le 2e_2-3.
 \end{equation}

Suppose that $e_1=5$. Then  $e_2=2$ and $d'\le 2$. Our bounds~\eqref{eq:d1}, \eqref{eq:7Count}  and \eqref{eq:d2} give that $d_1\le 7$, $d_1+d_2\ge 8$ and $d_2\le 1$, respectively. Thus, we have equalities everywhere. In particular, $F_1\in\CS{5}{7}{\C K_1}$, $F_2\in \CS{2}{1}{\C K_2}$, $d'=2$ and the pair of the involved $M$-edges belongs to $\CC{5}{7}{2}{1}$, a contradiction.

So, assume for the rest of the proof that $e_1=4$. Then, $d'\le e_2=3$. By~\eqref{eq:d1}, \eqref{eq:7Count}  and \eqref{eq:d2}, we have that $d_1\le 5$, $d_1+d_2\ge 12-2d'\ge 6$ and $d_2\le 3$. It follows that $d'\ge 2$.

First, suppose that $d'=2$. Thus, we must have that $d_1=5$ and $d_2=3$. Let $K_2^1,K_2^2\in F_2$ be the two edges matched by~$M$. Since $F_2\in\CS{3}{3}{\C K_2}$,
we have that $\{K_2^1,K_2^2\}\in\CS{2}{1}{\C K_2}$. Thus, we have a conflict from $\CC{4}{5}{2}{1}$, a contradiction.

Thus, $d'=3$, that is, $F_2$ is fully matched into $F_1$ by $M$. 

Suppose first that $d_1+d_2=6$. The case $(d_1,d_2)=(4,2)$ is impossible as then we get a conflict from $\CC{4}{4}{3}{2}$. Next, suppose that $(d_1,d_2)=(3,3)$. To avoid a conflict from $\CCP$, it must be the case that $F_1$ contains an $\CS{3}{3}{}$-subgraph, say with edges $K_1^1,K_1^2,K_1^3$. The matching $M$ matches at least two of these three edges, say to $K_2^1,K_2^2\in F_2$. Since every two edges of $F_2$ intersect, we have that $\{K_2^1,K_2^2\}\in\CS{2}{1}{\C K_2}$ and thus $M$ contains a conflict from $\CC{3}{3}{2}{1}$, a contradiction.  The only possible remaining case $(d_1,d_2)=(5,1)$ gives that $F_2\in\CS{3}{1}{\C K_2}$ is fully matched into $F_1\in\CS{4}{5}{\C K_1}$ by~$M$. But the two intersecting edges $K_2^1,K_2^2\in F_2$ form an element of $\CS{2}{1}{\C K_2}$ fully matched into $F_1$, that is, $M$ contains a conflict from $\CC{4}{5}{2}{1}$, a contradiction.

Thus we can assume that $d_1+d_2\ge 7$.
If $d_1=4$, then $d_2=3$ and $F_2\in \CS{3}{3}{\C K_2}$. Among the three $M$-matches of $F_2$ in $F_1$, some two, say $K_1^1,K_1^2\in F_1$, must intersect. (Indeed, otherwise these three edges contribute 0 to the defect of $F_1$ while  the remaining edge can contribute at most 3, a contradiction to $F_1\in\CS{4}{4}{\C K_1}$.) But then $\{K_1^1,K_1^2\}\in\CS{2}{1}{\C K_1}$ and $F_2\in \CS{3}{3}{\C K_2}$ give a conflict from $\CC{3}{3}{2}{1}$, a contradiction.
Thus it remains to consider the case when $d_1=5$. Note that all pairs of edges in $F_1\in \CS{4}{5}{\C K_1}$ except at most one pair intersect, and that $d_2\ge 7-d_1=2$.  
Thus, out of the three $M$-edges between $F_1$ and $F_2$, we can pick two such that their two endpoints on each side  intersect. Among the two remaining edges of $F_1$, at least one intersects both of these $F_1$-endpoints in two distinct vertices. Thus, we get two edges in $M$ between some sets in $\CS{3}{3}{\C K_1}$ and $\CS{2}{1}{\C K_2}$, which is a conflict from $\CC{3}{3}{2}{1}$. This final contradiction proves that $F=F(M)$ is indeed $(7r-12,7)$-free.\epf

Next, we need to bound from above the number of choices of $F_i\in \CS{e'}{d'}{\C K_i}$ for each involved pair~$(e',d')$. For this, we would like to construct each such $F_i$ from the empty $(r-2)$-graph by iteratively adding edges or pairs of edges at each step. Let $F'\subseteq F_i$ be the currently constructed subgraph. A \emph{$j$-attachment} for $j=0,1,2$ occurs when we add one new edge that shares exactly $j$ vertices with $V(F')$. To make a \emph{$3$-attachment}, we add two new edges $K,K'$ such that each of the three intersections $K\cap K'$, $K\cap V(F')$ and $K'\cap V(F')$ consists of exactly one vertex and these three vertices are pairwise distinct, that is, 
$$ 
 |K\cap K'|=|K\cap V(F')|=|K'\cap V(F')|=1\quad\mbox{and}\quad K\cap K'\cap V(F')=\emptyset.
 $$
 If we construct $F_i$ this way, then we let $a_j$ be the number of $j$-attachments for $0\le j\le 3$; note that then $a_0+a_1+a_2+2a_3=e'$ and $a_1+2a_2+3a_3=d'$.

\newcommand{\MS}[2]{\B T_{#1,#2}}

 \begin{claim}\label{cl:Struc} Let $i=1$ or $2$. Then, the following holds for every family $\C S$ listed in Table~\ref{ta:cliques}. 
 \begin{table}[h!]
 \begin{center}
 \begin{tabular}{|c|c|c|c|c|c|}
 \hline
 $\C S$ & $a_0$ & $a_1$ & $a_2$ & $a_3$ & $|\C S|$\\
 \hline\hline
 $\CS{2}{1}{\C K_i}$ &  1 & 1& 0 &0 & $O(\beta^2m^3)$\\
 \hline
  $\CS{3}{1}{\C K_i}$ & 2 & 1 & 0 & 0&$O(\beta^3m^5)$\\
 \hline
  $\CS{3}{2}{\C K_i}$ & 1 & 2 & 0 & 0& $O(\beta^3m^4)$ \\
 \hline
  $\CS{3}{3}{\C K_i}$ & 1 & 0 & 0 &1& $O(\beta^3m^3)$\\
 \hline
  $\CS{4}{4}{\C K_i}$ & 1 & 1 & 0 & 1& $O(\beta^4m^4)$\\
 \hline
  $\CS{4}{5}{\C K_i}$ &1 &0 & 1 & 1 &$O(\beta^3m^3)$\\
 \hline
  $\CS{5}{7}{\C K_i}$ & 1 & 0 & 2 &1& $O(\beta^3m^3)$\\ 
  \hline
 $\CCS{\C K_i}$& 1&3&0&0& $O(\beta^{4}m^{5})$\\
\hline
 \end{tabular}
 \vspace{1em}
 \caption{The values for Claim~\ref{cl:Struc}}
 \label{ta:cliques}
 \end{center}
 \end{table}
  \begin{enumerate}[(i)]
\item\label{it:attachments} Every $F_i\in \C S$ can be constructed using the above attachment operations starting from the empty graph such that the corresponding vector $(a_0,a_1,a_2,a_3)$ is exactly as stated in Table~\ref{ta:cliques}.
\item\label{it:CSSizes} The size of $\C S$ is at most the expression given in the last column of the table.
\end{enumerate}
  \end{claim}
 \bpf[Proof of Claim~\ref{cl:Struc}] Let us prove the first part for $F_i$ from a ``regular'' family $\C S=\CS{e'}{d'}{\C K_i}$. The cases $e'\le 3$  are straightforward to check, so assume that $e'\ge 4$.
 
Let $e'=4$. Take any $F_i=\{K^1,\dots,K^4\}$ in $\CS{4}{d'}{\C K_i}$.

Let $d'=4$. Suppose first that some three edges have defect $3$, say $\{K^1,K^2,K^3\}\in\CS{3}{3}{\C K_i}$. We can build this subgraph using a 0-attachment and a 3-attachment. The remaining edge $K^4$ contributes exactly 1 to the defect, so its addition is a 1-attachment, giving the desired. So suppose $F_i$ has no $\CS{3}{3}{}$-subgraph. Add one by one some two intersecting edges, say $K^1$ and $K^2$, doing a 0-attachment and a 1-attachment. By the $\CS{3}{3}{}$-freeness, each of the remaining edges $K^3$ and $K^4$ intersects $K^1\cup K^2$ in at most one vertex. The above three intersections contribute at most 3 to the total defect while $|K^3\cap K^4|\le 1$ contributes at most 1, so all these intersections are non-empty (and of size exactly 1). Thus we can add $\{K^3,K^4\}$ as one 3-attachment.

Let $d'=5$. Start with an intersecting pair of edges, say $K^1$ and $K^2$. Observe that at least one of the remaining edges, say $K^3$, intersects $K^1$ and $K^2$ at two different vertices (as otherwise the defect would be at most 4 by the linearity of $F$). We can construct $\{K^1,K^2,K^3\}\in \CS{3}{3}{\C K_i}$ using a 0-attachment and a 3-attachment. The addition of the remaining edge adds $2$ to the defect and thus is a 2-attachment. This gives the vector $(1,0,1,1)$, as desired.

So suppose that $(e',d')=(5,7)$.
Take any $F_i=\{K^1,\dots,K^5\}$ in $\CS{5}{7}{\C K_i}$. 
By the linearity of $F$, there are at least $d'=7$ pairs of intersecting edges so, by Mantel's theorem,  there are three pairwise intersecting edges, say  $K^1,K^2,K^3$. These edges contribute at most 3 to the defect. The pair of the remaining two edges contributes $|K^4\cap K^5|\le 1$ to the defect, so there are at least $d'-3-1=3$ further intersections. Thus, at least one of the edges $K^4, K^5$, say $K^4$, satisfies 
 \beq{eq:K4}
  |K^4\cap (K^1\cup K^2\cup K^3)|\ge 2.
  \eeq 
  It follows that some three of $K^1,\dots,K^4$ have defect 3, say $\{K^1,K^2,K^3\}\in\CS{3}{3}{\C K_i}$. By the same argument as before, we can assume that~\eqref{eq:K4} holds. Note that we must have equality in~\eqref{eq:K4} as otherwise the union of $K^1,\dots,K^4$ has at most $4(r-4)+2$ vertices, contradicting the high-girth assumption on $\C K_i$.
Thus, the defect of $\{K^1,\dots,K^4\}$ is exactly 5 and 
$$
|K^5\cap (K^1\cup\dots\cup K^4)|=d'-5=2.
$$
We conclude that we can build $F_i$ using a 0-attachment and a  3-attachment (to construct $\{K^1,K^2,K^3\}\in\CS{3}{3}{\C K_i}$) and two 2-attachments (to add $K^4$ and then $K^5$), as desired.

Let us turn to $F_i$ from the ``exceptional" family $\C S=\CCS{\C K_i}$. We know that $F_i$ has no isolated edge. Starting with  any edge $K^1\in F_i$ (a 0-attachment), add some edge intersecting it (a 1-attachment), say $K^2$. At least one of the two remaining edges, say $K^3$, has non-empty intersection with $K^1\cup K^2$ (as otherwise $d'\le 2$, a contradiction). This intersection consists of exactly one vertex (as otherwise the defect of $\{K^1,K^2,K^3\}$ is already $3$, a contradiction). Again, by $d'=3$, the remaining edge $K^4$ shares exactly one vertex with $K^1\cup K^2\cup K^3$. So when we add $K^3$ and then $K^4$, we use two 1-attachments, giving $(a_0,a_1,a_2,a_3)=(1,3,0,0)$, as desired.

Let us turn to the second part, namely bounding the number of choices of $F_i\in \C S$. We build each such $F_i$ as in Part~\ref{it:attachments}.
Given a partially built $F'\subseteq F_i$, there are by  Claim~\ref{cl:cliques} at most $O(\beta m^2)$, $O(\beta m)$, $O(1)$ and $O(\beta^2m)$ ways to do a $j$-attachment for $j=0,1,2,3$ respectively. For example, every 3-attachment can be obtained by taking some $u,v\in V(F')$, having at most ${3(r-2)\choose 2}=O(1)$ choices, then choosing some $w\in N_{G_i}(u)\cap N_{G_i}(v)$, having $O(\beta^2m)$ choices, and then taking (if they exist) the unique two edges of $\C K_i$ that contain the pairs $uw$ and $vw$. Thus, the number of $F_i\in\C S$ that give a fixed vector $(a_0,a_1,a_2,a_3)$ is at most 
\begin{equation}\label{eq:TUpper}
O\left((\beta m^2)^{a_0}\cdot (\beta m)^{a_1}\cdot (\beta^2m)^{a_3}\right)=O\left(\beta^{a_0+a_1+2a_3}\,m^{2a_0+a_1+a_3}\right).
\end{equation}
 By the first part, this directly gives an upper bound on $|\C S|$ stated in the table. 
\epf 

Recall that we defined $d=\alpha^{(r-2)^2}t$.

\begin{claim}\label{cl:Expectations} For every quadruple $(e_1,d_1,e_2,d_2)$  such that $\CC{e_1}{d_1}{e_2}{d_2}$ appears in the right-hand side of~\eqref{eq:CC}, the expected value of $|\CC{e_1}{d_1}{e_2}{d_2}|$ over the 
random graph $G_3$ is at most $O(d^{e_2}t)$. Also, the expectation of $|\CCP|$ is~$O(d^3t)$.
\end{claim} 
\bpf[Proof of Claim~\ref{cl:Expectations}] We can bound the expectation of $|\CC{e_1}{d_1}{e_2}{d_2}|$ from above by 
 $$
 \sum_{i=1}^2 |\CS{e_1}{d_1}{\C K_i}|\cdot |\CS{e_2}{d_2}{\C K_{3-i}}|\cdot \alpha^\ell,
 $$
 where $\ell$ is the smallest number of pairs in $A_1\times A_2$ that are covered by a full matching of some element of $\CS{e_2}{d_2}{\C K_{3-i}}$ into some element of $\CS{e_1}{d_1}{\C K_i}$. Using the upper bounds coming from Claim~\ref{cl:cliques} and recalling that $\beta=\alpha^{3/4}$, $t=\Theta(\beta m^2)$ and $d=\alpha^{(r-2)^2}t$, we obtain the following estimates:
 \begin{eqnarray*}
 \E{|\CC{3}{3}{2}{1}|}&=& O\left(\beta^3m^3\cdot \beta^2m^3\cdot \alpha^{2(r-2)^2-1}\right)\ =\ O(\beta^2\alpha^{-1}\cdot d^2t)\ =\ O(d^2t),\\
 \E{|\CC{4}{4}{3}{2}|}&=& O\left(\beta^4m^4\cdot \beta^3m^4\cdot \alpha^{3(r-2)^2-2}\right)\ =\ O(\beta^3\alpha^{-2}\cdot d^3t)\ =\ O(d^3t),\\
\E{|\CC{4}{5}{2}{1}|}&=& O\left(\beta^3m^3\cdot \beta^2m^3\cdot \alpha^{2(r-2)^2-1}\right)\ =\ O(\beta^2\alpha^{-1}\cdot d^2t)\ =\ O(d^2t),\\
 \E{|\CC{5}{7}{2}{1}|}&=& O\left(\beta^3m^3\cdot \beta^2m^3\cdot \alpha^{2(r-2)^2-1}\right)\ =\ O(\beta^2\alpha^{-1}\cdot d^2t)\ =\ O(d^2t).
 \end{eqnarray*}
 Finally, the obvious adaptation of this argument to the ``exceptional'' family $\CCP$ gives that
 $$
  \E{|\CCP|}\ =\ O\left(\beta^4m^5\cdot \beta^3m^3\cdot \alpha^{3(r-2)^2-2}\right)\ =\ O(\beta^3\alpha^{-2}\cdot d^3t)\ =\ O(d^3t).
$$ 
This finishes the proof   of Claim~\ref{cl:Expectations}.
\epf

Fix $G_3$ such that for each of the five families from Claim~\ref{cl:Expectations} its size is at most, say, 10 times the expected value, $|H|=\Theta(td)$ and  $|G_3|\le 2\alpha m^2$. This is possible since the last two properties hold with high probability by the Bounded Difference Inequality (Lemma~\ref{lem:BDI}) and the union bound.

Finally, let us describe how to construct a matching in $H$ avoiding all conflicts listed in Claim~\ref{cl:7Free}. We shall use the probabilistic deletion method. 
Pick every edge of $H$ randomly with probability $\mu/d$, where  $\mu$ is a sufficiently small constant which depends on the implicit constants in the asymptotic notations but not on $\alpha$, that is,  $1/r\gg \mu\gg \alpha$.
Clearly, the expected number of chosen edges is $(\mu/d) |H|=\Theta (\mu t)$ and the expected number of pairs of edges which overlap is $O(d^2t\cdot (\mu/d)^2)= O(\mu^2 t)$. By Claim~\ref{cl:Expectations}, the expected number of elements of $\C C$ all of whose edges have been chosen is at most
\[
 O\left(d^2t\cdot (\mu/d)^2 +  d^3t\cdot (\mu/d)^3\right) = O(\mu^2 t).
\] 

Let $M$ be obtained from the $\mu$-random subset of $E(H)$ by removing edges which overlap with some other chosen edge or participate in a conflict from~$\C C$ with all of its edges being chosen. By construction, $M$ is a matching in $H$ that avoids all conflicts. Also,
$$
\E{|M|}\ge \Omega(\mu t)-O(\mu^2 t)=\Omega(\mu t).
$$ 
Take an outcome such that $|M|$ is at least its expectation and define $F=F(M)$ by~\eqref{eq:F}.

Let us check that the obtained $r$-graph $F$ satisfies all stated properties listed in Lemma~\ref{57cons}. The first two, namely Parts~\ref{it:57Free} and~\ref{it:57Free'}, follow from Claims~\ref{cl:small}, \ref{cl:C3321} and~\ref{cl:7Free}.
Also, the size of $F$ is $2\,|M| =\Omega_{\mu}(t)
=\Omega_{\mu}(\alpha^{3/4} m^2)$, proving Part~\ref{it:|F|}.

Let us turn to Part~\ref{it:J}. Since $F$ consists of $|F|/2$ diamonds that do not share any pairs, we have that $|\p{1}{F}|= (2{r\choose 2}-1)|F|/2$. Since $|G_3|\le 2\alpha m^2$,
it is enough to  show that every pair  $xy\in \p{\le3}{F}\setminus \p{1}{F}$ is an edge of~$G_3$. 
Since $F$ has no $(3r-4,3)$-configuration, we have that $xy\in\pp{1}{2}{F}$. Since every pair of cliques in $\C K_1\cup \C K_2$ shares at most one vertex, some diamond coming from $\{K_1,K_2\}\in M$ 2-claims the pair $xy$. Thus $xy\in \{ab\mid a\in K_1,\ b\in K_2\}\subseteq E(G_3)$, as desired. This finishes the proof of Lemma~\ref{57cons}.\epf

\bpf[Proof of the lower bounds in Theorem~\ref{5edge} and~\ref{7edge} with $r\geq 4$]
Let $F$ be the $r$-graph given by Lemma~\ref{57cons}.
Thus $F$ is $\cG{r}{k}$-free for $k\in \{5,7\}$, $|F|= \Omega_r(\alpha^{3/4} m^2)$ and 
$$
 |\p{\le3}{F}|\leq 
 \frac{r^2-r-1}{2}\,|F|+2\alpha m^2.
 $$
 By Theorem~\ref{highgirth}, for each $k\in \{5,7\}$, we have that
$$
\liminf_{n\rightarrow \infty}n^{-2}f^{(r)}(n;k(r-2)+2,k)\geq \frac{|F|}{2\,|\p{\le3}{F}|}\geq \frac{1}{r^2-r+1+O(\alpha^{1/4})}.
$$
The lower bound $\frac{1}{r^2-r+1}$ in Theorems~\ref{5edge} and~\ref{7edge} follows by taking $\alpha \rightarrow 0$.
\epf

\subsection{Lower bounds in Theorems~\ref{6edge3uniform} and~\ref{6edgehigh}}

\begin{figure}[ht]
\centering
\definecolor{qqffqq}{rgb}{0,1,0}
\definecolor{ffqqqq}{rgb}{1,0,0}
\definecolor{rvwvcq}{rgb}{0.08235294117647059,0.396078431372549,0.7529411764705882}
\begin{tikzpicture}[line cap=round,line join=round,x=1cm,y=1cm]
\clip(-5.5,-3) rectangle (5.5,3);
\fill[line width=0.8pt,color=black,fill=black,fill opacity=1] (-5,0) -- (5,0) -- (0,0) -- cycle;
\fill[line width=0.8pt,color=qqffqq,fill=qqffqq,fill opacity=0.1] (-1,-2) -- (1,-2) -- (-5,0) -- cycle;
\fill[line width=0.8pt,color=qqffqq,fill=qqffqq,fill opacity=0.1] (1,-2) -- (5,0) -- (-1,-2) -- cycle;
\fill[line width=0.8pt,color=qqffqq,fill=qqffqq,fill opacity=0.1] (-1,2) -- (-5,0) -- (1,2) -- cycle;
\fill[line width=0.8pt,color=qqffqq,fill=qqffqq,fill opacity=0.1] (1,2) -- (5,0) -- (-1,2) -- cycle;
\fill[line width=0.8pt,color=rvwvcq,fill=rvwvcq,fill opacity=0.10000000149011612] (2,0.5) -- (1.4,0.8) -- (1,2) -- cycle;
\fill[line width=0.8pt,color=rvwvcq,fill=rvwvcq,fill opacity=0.10000000149011612] (2,0.5) -- (0,0) -- (1.4,0.8) -- cycle;
\fill[line width=0.8pt,color=rvwvcq,fill=rvwvcq,fill opacity=0.10000000149011612] (-1.4,0.8) -- (-2,0.5) -- (0,0) -- cycle;
\fill[line width=0.8pt,color=rvwvcq,fill=rvwvcq,fill opacity=0.10000000149011612] (-1.4,0.8) -- (-1,2) -- (-2,0.5) -- cycle;
\fill[line width=0.8pt,color=ffqqqq,fill=ffqqqq,fill opacity=0.1] (-1,2.5) -- (-0.5,2.5) -- (0,0) -- cycle;
\fill[line width=0.8pt,color=ffqqqq,fill=ffqqqq,fill opacity=0.1] (-1,2.5) -- (-1,2) -- (-0.5,2.5) -- cycle;
\fill[line width=0.8pt,color=ffqqqq,fill=ffqqqq,fill opacity=0.1] (0.5,2.5) -- (1,2.5) -- (1,2) -- cycle;
\fill[line width=0.8pt,color=ffqqqq,fill=ffqqqq,fill opacity=0.1] (1,2.5) -- (0,0) -- (0.5,2.5) -- cycle;
\fill[line width=0.8pt,color=rvwvcq,fill=rvwvcq,fill opacity=0.10000000149011612] (-2,-0.5) -- (-1.4,-0.8) -- (0,0) -- cycle;
\fill[line width=0.8pt,color=rvwvcq,fill=rvwvcq,fill opacity=0.10000000149011612] (1.4,-0.8) -- (2,-0.5) -- (0,0) -- cycle;
\fill[line width=0.8pt,color=rvwvcq,fill=rvwvcq,fill opacity=0.10000000149011612] (-2,-0.5) -- (-1,-2) -- (-1.4,-0.8) -- cycle;
\fill[line width=0.8pt,color=rvwvcq,fill=rvwvcq,fill opacity=0.10000000149011612] (1.4,-0.8) -- (1,-2) -- (2,-0.5) -- cycle;
\fill[line width=0.8pt,color=ffqqqq,fill=ffqqqq,fill opacity=0.10000000149011612] (-1,-2.5) -- (-0.5,-2.5) -- (0,0) -- cycle;
\fill[line width=0.8pt,color=ffqqqq,fill=ffqqqq,fill opacity=0.10000000149011612] (-1,-2.5) -- (-1,-2) -- (-0.5,-2.5) -- cycle;
\fill[line width=0.8pt,color=ffqqqq,fill=ffqqqq,fill opacity=0.10000000149011612] (1,-2.5) -- (1,-2) -- (0.5,-2.5) -- cycle;
\fill[line width=0.8pt,color=ffqqqq,fill=ffqqqq,fill opacity=0.10000000149011612] (1,-2.5) -- (0,0) -- (0.5,-2.5) -- cycle;
\draw [line width=0.8pt,color=rvwvcq] (-5,0)-- (5,0);
\draw [line width=0.8pt,color=rvwvcq] (5,0)-- (0,0);
\draw [line width=0.8pt,color=rvwvcq] (0,0)-- (-5,0);
\draw [line width=0.8pt,color=qqffqq] (-1,-2)-- (1,-2);
\draw [line width=0.8pt,color=qqffqq] (1,-2)-- (-5,0);
\draw [line width=0.8pt,color=qqffqq] (-5,0)-- (-1,-2);
\draw [line width=0.8pt,color=qqffqq] (1,-2)-- (5,0);
\draw [line width=0.8pt,color=qqffqq] (5,0)-- (-1,-2);
\draw [line width=0.8pt,color=qqffqq] (-1,-2)-- (1,-2);
\draw [line width=0.8pt,color=qqffqq] (-1,2)-- (-5,0);
\draw [line width=0.8pt,color=qqffqq] (-5,0)-- (1,2);
\draw [line width=0.8pt,color=qqffqq] (1,2)-- (-1,2);
\draw [line width=0.8pt,color=qqffqq] (1,2)-- (5,0);
\draw [line width=0.8pt,color=qqffqq] (5,0)-- (-1,2);
\draw [line width=0.8pt,color=qqffqq] (-1,2)-- (1,2);
\draw [line width=0.8pt,color=rvwvcq] (2,0.5)-- (1.4,0.8);
\draw [line width=0.8pt,color=rvwvcq] (1.4,0.8)-- (1,2);
\draw [line width=0.8pt,color=rvwvcq] (1,2)-- (2,0.5);
\draw [line width=0.8pt,color=rvwvcq] (2,0.5)-- (0,0);
\draw [line width=0.8pt,color=rvwvcq] (0,0)-- (1.4,0.8);
\draw [line width=0.8pt,color=rvwvcq] (1.4,0.8)-- (2,0.5);
\draw [line width=0.8pt,color=rvwvcq] (-1.4,0.8)-- (-2,0.5);
\draw [line width=0.8pt,color=rvwvcq] (-2,0.5)-- (0,0);
\draw [line width=0.8pt,color=rvwvcq] (0,0)-- (-1.4,0.8);
\draw [line width=0.8pt,color=rvwvcq] (-1.4,0.8)-- (-1,2);
\draw [line width=0.8pt,color=rvwvcq] (-1,2)-- (-2,0.5);
\draw [line width=0.8pt,color=rvwvcq] (-2,0.5)-- (-1.4,0.8);
\draw [line width=0.8pt,color=ffqqqq] (-1,2.5)-- (-0.5,2.5);
\draw [line width=0.8pt,color=ffqqqq] (-0.5,2.5)-- (0,0);
\draw [line width=0.8pt,color=ffqqqq] (0,0)-- (-1,2.5);
\draw [line width=0.8pt,color=ffqqqq] (-1,2.5)-- (-1,2);
\draw [line width=0.8pt,color=ffqqqq] (-1,2)-- (-0.5,2.5);
\draw [line width=0.8pt,color=ffqqqq] (-0.5,2.5)-- (-1,2.5);
\draw [line width=0.8pt,color=ffqqqq] (0.5,2.5)-- (1,2.5);
\draw [line width=0.8pt,color=ffqqqq] (1,2.5)-- (1,2);
\draw [line width=0.8pt,color=ffqqqq] (1,2)-- (0.5,2.5);
\draw [line width=0.8pt,color=ffqqqq] (1,2.5)-- (0,0);
\draw [line width=0.8pt,color=ffqqqq] (0,0)-- (0.5,2.5);
\draw [line width=0.8pt,color=ffqqqq] (0.5,2.5)-- (1,2.5);
\draw [line width=0.8pt,color=rvwvcq] (-2,-0.5)-- (-1.4,-0.8);
\draw [line width=0.8pt,color=rvwvcq] (-1.4,-0.8)-- (0,0);
\draw [line width=0.8pt,color=rvwvcq] (0,0)-- (-2,-0.5);
\draw [line width=0.8pt,color=rvwvcq] (1.4,-0.8)-- (2,-0.5);
\draw [line width=0.8pt,color=rvwvcq] (2,-0.5)-- (0,0);
\draw [line width=0.8pt,color=rvwvcq] (0,0)-- (1.4,-0.8);
\draw [line width=0.8pt,color=rvwvcq] (-2,-0.5)-- (-1,-2);
\draw [line width=0.8pt,color=rvwvcq] (-1,-2)-- (-1.4,-0.8);
\draw [line width=0.8pt,color=rvwvcq] (-1.4,-0.8)-- (-2,-0.5);
\draw [line width=0.8pt,color=rvwvcq] (1.4,-0.8)-- (1,-2);
\draw [line width=0.8pt,color=rvwvcq] (1,-2)-- (2,-0.5);
\draw [line width=0.8pt,color=rvwvcq] (2,-0.5)-- (1.4,-0.8);
\draw [line width=0.8pt,color=ffqqqq] (-1,-2.5)-- (-0.5,-2.5);
\draw [line width=0.8pt,color=ffqqqq] (-0.5,-2.5)-- (0,0);
\draw [line width=0.8pt,color=ffqqqq] (0,0)-- (-1,-2.5);
\draw [line width=0.8pt,color=ffqqqq] (-1,-2.5)-- (-1,-2);
\draw [line width=0.8pt,color=ffqqqq] (-1,-2)-- (-0.5,-2.5);
\draw [line width=0.8pt,color=ffqqqq] (-0.5,-2.5)-- (-1,-2.5);
\draw [line width=0.8pt,color=ffqqqq] (1,-2.5)-- (1,-2);
\draw [line width=0.8pt,color=ffqqqq] (1,-2)-- (0.5,-2.5);
\draw [line width=0.8pt,color=ffqqqq] (0.5,-2.5)-- (1,-2.5);
\draw [line width=0.8pt,color=ffqqqq] (1,-2.5)-- (0,0);
\draw [line width=0.8pt,color=ffqqqq] (0,0)-- (0.5,-2.5);
\draw [line width=0.8pt,color=ffqqqq] (0.5,-2.5)-- (1,-2.5);
\begin{scriptsize}
\draw [fill=black] (-5,0) circle (1.5pt);
\draw [fill=black] (0,0) circle (1.5pt);
\draw [fill=black] (5,0) circle (1.5pt);
\draw [fill=black] (-1,-2) circle (1.5pt);
\draw [fill=black] (1,-2) circle (1.5pt);
\draw [fill=black] (-1,2) circle (1.5pt);
\draw [fill=black] (1,2) circle (1.5pt);
\draw [fill=black] (2,0.5) circle (1.5pt);
\draw [fill=black] (1.4,0.8) circle (1.5pt);
\draw [fill=black] (-1.4,0.8) circle (1.5pt);
\draw [fill=black] (-2,0.5) circle (1.5pt);
\draw [fill=black] (-1,2.5) circle (1.5pt);
\draw [fill=black] (-0.5,2.5) circle (1.5pt);
\draw [fill=black] (0.5,2.5) circle (1.5pt);
\draw [fill=black] (1,2.5) circle (1.5pt);
\draw [fill=black] (-1.4,-0.8) circle (1.5pt);
\draw [fill=black] (1.4,-0.8) circle (1.5pt);
\draw [fill=black] (-2,-0.5) circle (1.5pt);
\draw [fill=black] (2,-0.5) circle (1.5pt);
\draw [fill=black] (-1,-2.5) circle (1.5pt);
\draw [fill=black] (-0.5,-2.5) circle (1.5pt);
\draw [fill=black] (1,-2.5) circle (1.5pt);
\draw [fill=black] (0.5,-2.5) circle (1.5pt);

\draw[color=black] (-5.3,0) node {$x_2$};
\draw[color=black] (5.3,0) node {$x_3$};

\draw[color=black] (-2.4,0.27) node {$D'(x_1, a_1)$};

\draw[color=black] (2.4,0.27) node {$D'(x_1, b_1)$};

\draw[color=black] (-2.4,-0.27) node {$D'(x_1, a_1')$};

\draw[color=black] (2.4,-0.27) node {$D'(x_1, b_1')$};

\draw[color=black] (-0.8,2.75) node {$D(x_1, a_1)$};

\draw[color=black] (0.8,2.75) node {$D(x_1, b_1)$};

\draw[color=black] (-0.8,-2.75) node {$D(x_1, a_1')$};

\draw[color=black] (0.8,-2.75) node {$D(x_1, b_1')$};

\draw[color=black] (-1.2,2.1) node {$a_1$};
\draw[color=black] (1.2,2.1) node {$b_1$};

\draw[color=black] (-1.2,-2.1) node {$a'_1$};
\draw[color=black] (1.2,-2.1) node {$b'_1$};
\end{scriptsize}
\end{tikzpicture}
\caption{The figure depicts the subgraph of the 3-graph $F_{63}$ ``lying'' on the pair $x_2x_3$. 
The central vertex in the figure is $x_1$, and the green diamonds correspond to $D_1$ and $D_1'$.
Copies of the same construction ``lie'' on the pairs $x_1x_2$ and $x_2x_3$ in $F_{63}$.}
\label{fig (8,6)}
\end{figure}

\bpf[Proof of the lower bound in Theorem~\ref{6edge3uniform}]
We define a $3$-graph $F_{63}$ on $63$ vertices with $61$ edges as follows. 
Let $T$ be the $3$-graph which is obtained from an edge $x_1x_2x_3$ by adding, for every $i\in [3]$, two  diamonds $D_i=\{a_ib_ix_s,a_ib_ix_t\}$ and $D'_i=\{a'_ib'_ix_s,a'_ib'_ix_t\}$ where $\{s,t\} = [3]\setminus\{i\}$. We say that these 6 diamonds are of  \emph{level~1}.
Then, consider the following $12$ pairs, which do  not belong to $\p{1} T$, namely \beq{eq:ClaimedPairs}
 x_ia_i,x_ib_i,x_ia'_i,x_ib'_i,\quad \mbox{for $i\in [3]$}.
 \eeq
Let $F_{63}$ be obtained from $T$ by adding, for every such pair $xy$, two vertex-disjoint diamonds $D(x, y)$ and $D'(x, y)$ \OutIn{1}{2}-claiming $xy$, calling these 24 diamonds \emph{of level~2}. 
Thus, $F_{63}$ has $3+6\cdot 2+24\cdot 2=63$ vertices and $1+6\cdot 2+24\cdot 2=61$ edges. 

To show that $F_{63}$ is $\cG{3}{6}$-free, we first prove the following claim. 

\begin{claim}\label{cl:procedure}
For every subgraph $G\subseteq F_{63}$, there is an integer $t \geq 1$ and a sequence $\emptyset = G_0\subseteq G_1\subseteq G_2\subseteq\ldots \subseteq G_t = G$ where, for every $i\in [0,t-1]$, $G_{i+1}\setminus G_i$ is 
either a single edge that shares at most one vertex with $G_i$, 
or a diamond  that shares at most two vertices with~$G_i$ so that if a pair is shared, then this pair is \OutIn{1}{2}-claimed by $G_{i+1}\setminus G_i$.
\end{claim}
\bpf[Proof of Claim~\ref{cl:procedure}]
We start from $G$ and consecutively delete single edges or diamonds with the required property, 
thus implicitly defining the sequence $G_1, G_2,\ldots, G_t$ in the reverse order. 
Suppose that $G$ contains an edge $X$ in some diamond $\{X,Y\}$ of level~2.
If $Y\notin G$, then $G\setminus \{X\}$ and $X$ share at most one vertex. 
Otherwise, $\{X,Y\}$ attaches to $G\setminus \{X,Y\}$ via at most two vertices and if there are exactly two shared vertices, then this pair is \OutIn{1}{2}-claimed by $\{X,Y\}$. 
Thus, the edge $X$ (in the first case) and the diamond $\{X,Y\}$ (in the second case) can be deleted from $G = G_t$ to obtain $G_{t-1}$.
Repeating this step, one can  delete all edges of $G$ from diamonds of level~2. 
After this, all edges in $G$ from diamonds of level~1 can be deleted in a similar way,  finishing the proof.
\epf

Note that, for any sequence returned by
Claim~\ref{cl:procedure}, the first non-empty $r$-graph $G_1$ contains two more vertices than edges and each new attachment cannot decrease this difference. It immediately follows that  $F_{63}$ is $(j+1,j)$-free for every $j$ and, in particular, for $j\in [3,5]$. Furthermore, if there exists some subgraph $G\subseteq F_{63}$ with $6$ edges and at most $8$ vertices, then the sequence $(G_i)_{i=0}^r$ returned by
Claim~\ref{cl:procedure} satisfies by the parity of $|G|$ that $t=3$ and that each of $G_1$, $G_2\setminus G_1$ and $G_3\setminus G_2$ is a diamond. Let us argue that this is impossible. 
If the diamond $G_1$ is of level 1, say $D_1$, then the only possibility for  the diamond $G_2\setminus G_1$ is $D_1'$, but then no other diamond $D$ of $F_{63}$ connects to $G_2=D_1\cup D_1'$ via a pair \OutIn{1}{2}-claimed by~$D$, a contradiction. Similarly, if $G_1$ is of level 2, say $D(x_1,a_1)$, then $G_2\setminus G_1$ must be $D'(x_1,a_1)$ and, again, there is no suitable choice for the third diamond~$G_3\setminus G_2$. We conclude that $F_{63}$ is $(8,6)$-free and thus $\cG{3}{6}$-free.

Note 
that the only $(4,2)$-configurations in $F$ are the 6 diamonds of level~1 (\OutIn{1}{2}-claiming only the pairs already 1-claimed by the central edge $x_1x_2x_3$) and the 24 diamonds of level~2 (\OutIn{1}{2}-claiming the 12 pairs in~\eqref{eq:ClaimedPairs}). Also, 
every $(5,3)$-configuration in $F_{63}$ consists of the central edge $x_1x_2x_3$ plus a diamond of level~1 (with both of its \OutIn{1}{3}-claimed pairs being already \OutIn{}{2}-claimed by diamonds of level~2). It follows that
$\p{\le3}{F_{63}}$ consists of $\p{1} {F_{63}}$ and the $12$ pairs in~\eqref{eq:ClaimedPairs}.
Thus $|\p{\le3}{F_{63}}|=3+(6+24)\cdot 5+12=165$. Using Theorem~\ref{highgirth}, we derive that
$$\liminf_{n\rightarrow \infty}n^{-2}f^{(3)}(n;8,6)\geq \frac{|F_{63}|}{2\,|\p{\le3}{F_{63}}|}=\frac{61}{330},$$
as desired.
\epf

\bpf[Proof of the lower bound in Theorem~\ref{6edgehigh}]
For $r\geq 4$, the lower bound on $f^{(r)}(n;6r-10,6)$ follows by Theorem~\ref{highgirth} with the $r$-graph $F$ being a single edge, for which $\p{\le3}{F}=\p{1}{F}$ has size~${r\choose 2}$. 
\epf

\section{Upper Bounds}\label{se:upper}

\subsection{Some common definitions and results}\label{se:UpperCommon}

Recall that for an $r$-graph $F$ and a pair $uv$, the set $\CI_F(uv)$ 
consists of all integers $j\ge 0$ such that $F$ has $j$ edges that together with $uv$ include at most $rj-2j+2$ vertices. 
Note that, by definition, 0 always belongs to $\CI_{G}(uv)$, which is notationally convenient in the statement of the following easy but very useful observation.

\begin{lemma}\label{lm: 2} For any $\C F^{(r)}(rk-2k+2,k)$-free $r$-graph $G$, any $uv\in {V(G)\choose 2}$ and any edge-disjoint subgraphs $F_1,\dots,F_s\subseteq G$, the  sum-set $\sum_{i=1}^s \CI_{F_i}(uv)=\left\{\sum_{i=1}^s m_i\mid m_i\in \CI_{F_i}(uv)\right\}$ does not contain~$k$.
\end{lemma}
\bpf If some sequence of $m_i$'s sums up to exactly $k$, then the corresponding $k$ edges of $G$ are all distinct and span at most $2+\sum_{i=1}^s (r-2)m_i = rk-2k+2$ vertices, a contradiction. 
\epf

As we mentioned in the introduction, our proof strategy to bound the size of an $(rk-2k+2,k)$-free $r$-graph $G$ from above is to analyse possible isomorphism types of the parts of some partition of $E(G)$ which is obtained from the trivial partition into single edges by iteratively applying some merging rules. 
Unfortunately, we did not find a single rule that works in all cases that are studied here.
In fact, for every new pair $(r,k)$ resolved in this paper except for $(3,5)$,
we build the final partition in stages (with each stage having a different merging rule) as the intermediate families are also needed in our analysis. 
Let us now develop some general notation and prove some basic results related to merging. 

Let $G$ be an arbitrary $r$-graph. When dealing with a partition $\C P$ of $E(G)$, we will view each element $F\in\C P$ as an $r$-graph whose vertex set is the union of the edges in~$F$.
Let $A,B\subseteq \I N$ be any (not necessarily disjoint) sets of positive integers. For two subgraphs $F,H\subseteq G$, if they are edge disjoint and there is a pair $uv$ such that $A\subseteq \CI_F(uv)$ and $B\subseteq \CI_H(uv)$, then we say that $F$ and $H$ are \emph{$\OMerge{A}{B}$-mergeable (via $uv$)}. Note that this relation is not symmetric in $F$ and $H$: the first (resp.\ second) $r$-graph $A$-claims (resp.\ $B$-claims) the pair~$uv$.
When the ordering of the two $r$-graphs does not matter, we use the shorthand \emph{$\Merge{A}{B}$-mergeable} to mean $\OMerge{A}{B}$-mergeable or $\OMerge{B}{A}$-mergeable.
For a partition $\C P$ of $E(G)$, its \emph{$\Merge{A}{B}$-merging} is the partition $\C M_{\Merge{A}{B}}(\C P)$ of $E(G)$ obtained from $\C P$ by iteratively and as long as possible taking a pair of distinct $\Merge{A}{B}$-mergeable parts in the current partition
and replacing them by their union.
Note that, if $F$ and $H$ are $\Merge{A}{B}$-mergeable via $uv$, then
\begin{equation}\label{eq:AB}
A\cup B \subseteq \CI_{F\cup H}(uv)\quad \text{and}\quad A+B\subseteq \CI_{F\cup H}(uv).    
\end{equation}
The first inclusion implies that, in particular, the order of merging operations does not affect the partition $\C M_{\Merge{A}{B}}(\C P)$.
Note that the final partition $\C M_{\Merge{A}{B}}(\C P)$ is a coarsening of $\C P$ and contains no $\Merge{A}{B}$-mergeable pairs of $r$-graphs. 
When $\C P$ is clear from the context, we may refer to the elements of $\C M_{\Merge{A}{B}}(\C P)$ as \emph{$\Merge{A}{B}$-\cluster{}s}. 
Likewise, a subgraph $F$ of $G$ that can appear as a part in some intermediate stage of the $\Merge{A}{B}$-merging process starting with $\C P$ is called a \emph{partial $\Merge{A}{B}$-\cluster{}} and we let $\C M'_{\Merge{A}{B}}(\C P)$ denote the set of all partial $\Merge{A}{B}$-\cluster{}s. In other words, $\C M'_{\Merge{A}{B}}(\C P)$ is the smallest family of $r$-graphs which contains $\C P$ as a subfamily and is closed under taking the union of $\Merge{A}{B}$-mergeable elements. The monotonicity of the merging rule implies that $\C M_{\Merge{A}{B}}(\C P)$ is exactly the set of maximal (by inclusion) elements of $\C M'_{\Merge{A}{B}}(\C P)$ and that the final partition $\C M_{\Merge{A}{B}}(\C P)$ does not depend on the order in which we merge parts.

In the frequently occurring case when $A=\{1\}$ and $B=\{j\}$, we abbreviate $\OMerge{\{1\}}{\{j\}}$ to $\OMerge{}{j}$ and $\Merge{\{1\}}{\{j\}}$ to $\Merge{}{j}$ in the above nomenclature.
Thus, $\OMerge{}{j}$-mergeable (resp.\ $\Merge{}{j}$-mergeable) means $\OMerge{\{1\}}{\{j\}}$-mergeable (resp.\ $\Merge{\{1\}}{\{j\}}$-mergeable). 

As an example, let us look at the following merging rule that is actually used as the first step in each of our proofs of the upper bounds. Namely, given $G$, let 
$$
 \OneClusters:=\C M_{\Merge{\{1\}}{\{1\}}}(\C P_{\mathrm{trivial}})
$$ 
 be the $\Merge{}{1}$-merging  of the trivial partition $\C P_{\mathrm{trivial}}$ of $G$ into single edges. We call the elements of $\OneClusters$ \emph{\OneCluster{}s}. Here is an alternative description of~$\OneClusters$. 
Call a subgraph $F\subseteq G$ \emph{connected} if for any two edges $X,Y\in F$ there is a sequence of edges $X_1=X,X_2,\ldots,X_m=Y$ in $F$ such that, for every $i\in [m-1]$, we have $|X_i\cap X_{i+1}|\geq 2$.   
 Then, \OneCluster{}s are exactly maximal connected subgraphs of $G$  (and partial \OneCluster{}s are exactly connected subgraphs).

We will also use (often without explicit mention) the following result, which is a generalisation of the well-known fact that we can remove edges from any connected 2-graph one by one, down to any given connected subgraph, while keeping the edge set connected. 
The assumption~\eqref{eq:Trim} states that, roughly speaking, the merging process cannot create any new mergeable pairs.

\begin{lemma}[Trimming Lemma]\label{lm:Trim} 
Fix an $r$-graph $G$, a partition $\C P$ of $E(G)$ and sets \mbox{$A,B\subseteq \I N$}.
\begin{equation}\label{eq:Trim}
\begin{gathered}
\mbox{Suppose that, for all $\OMerge{A}{B}$-mergeable (and thus edge-disjoint) $F,H\in\C M'_{\Merge{A}{B}}(\C P)$,}\\
\mbox{there exist $\OMerge{A}{B}$-mergeable  $F',H'\in\C P$ such that $F'\subseteq F$ and $H'\subseteq H$.}
\end{gathered}
\end{equation}
Then, for every partial $\Merge{A}{B}$-\cluster{}
$F_0\subseteq F$, there is an ordering $F_1,\dots,F_s$ of the elements of $\C P$ that lie inside $F\setminus F_0$ such that, for every $i\in [s]$, $\bigcup_{j=0}^{i-1} F_j$ and $F_i$ are $\Merge{A}{B}$-mergeable (and, in particular, $\bigcup_{j=0}^{i} F_j$ is a partial $\Merge{A}{B}$-\cluster{} for every $i\in [s]$).
\end{lemma}
\bpf
Suppose that the lemma is false. 
For every $H\in\C M'_{\Merge{A}{B}}(\C P)$, denote by $|H|_{\C P}$ the number of elements of $\C P$ contained in $H$.
Choose a counterexample $(F_0,F)$ with $m:=|F|_{\C P}-|F_0|_{\C P}$ smallest possible.
Note that $m>0$ as otherwise $F_0=F$ and the conclusion of the lemma vacuously holds. Consider some $\Merge{A}{B}$-merging process which leads to $F$; let $H$ be the first occurring partial $\Merge{A}{B}$-\cluster{} that shares at least one edge with each of $F_0$ and $F\setminus F_0$. Since the final partial $\Merge{A}{B}$-\cluster{} $F$ satisfies both these properties, $H$ exists. As $F_0$ and $F\setminus F_0$ are unions of some elements of $\C P$, we have $H\notin\C P$. So, let the merging process for $F$ build $H$ as the union of $\OMerge{A}{B}$-mergeable  partial $\Merge{A}{B}$-\cluster{}s $H_A,H_B\subsetneq H$. By the definition of $H$, one of $H_A$ and $H_B$, say $H_A$, shares an edge with $F_0$ but not with $F\setminus F_0$ while the opposite holds for $H_B$. 
Then, $H_A\subseteq F_0$ and $H_B\subseteq F\setminus F_0$.

Since the partial $\Merge{A}{B}$-\cluster{}s $H_A$ and $H_B$ are $\OMerge{A}{B}$-mergeable, the assumption of the lemma implies that there are $\OMerge{A}{B}$-mergeable $H'_A\subseteq H_A$ and $H'_B\subseteq H_B$ with $H'_A,H'_B\in\C P$. Then, $H'_B$, which is edge-disjoint from $F_0$, is $\Merge{A}{B}$-mergeable with $F_0\supseteq H_A$ as well. 
The minimality of $m$ guarantees an ordering $F_2,\dots, F_m$ of the elements of $\mathcal{P}$ that lie inside $F\setminus (F_0\cup H'_B)$ satisfying the statement of the claim for the pair $(F_0\cup H'_B,F)$. However, then the ordering $H'_B, F_2,\dots, F_m$ satisfies the statement of the claim for $(F_0,F)$, so this pair cannot be a counterexample, on the contrary to our assumption. 
\epf

In the special case $A=B=\{1\}$ (when partial \cluster{}s are just connected subgraphs), the assumption of Lemma~\ref{lm:Trim} is vacuously true. Since we are going to use its conclusion quite often,
we state it separately.

\begin{corollary}\label{cr:Trim} For every pair $F_0\subseteq F$ of connected $r$-graphs, there is an ordering $X_1,\dots,X_s$ of the edges in $F\setminus F_0$ such that, for every $i\in [s]$, the $r$-graph $F_0\cup \{X_1,\dots,X_i\}$ is connected.\qed\end{corollary}

We say that an $r$-graph is a \emph{$1$-tree} if it contains only one edge. 
For $i\geq 2$, we recursively define an \emph{$i$-tree} as any $r$-graph that can be obtained from an $(i-1)$-tree $T$ by adding a new edge that consists of a pair $ab$ in the $2$-shadow $\p{1} T$ of $T$ and $r-2$ new vertices (not present in~$T$). Clearly, every $i$-tree is connected. Like the usual 2-graph trees, $i$-trees are the ``sparsest" connected $r$-graphs of given size. Any $i$-tree $T$ satisfies
\beq{eq:Sizes}
 |\p{1} T|=i{r\choose 2}-i+1\quad\mbox{and}\quad |\pp{1}{2} T|\geq (i-1) (r-2)^2. 
\eeq
 (Recall that $\pp{1}{2} T$ is the set of pairs which are $2$-claimed but not \OneClaim{}ed by~$T$.) Note that the second inequality in~\eqref{eq:Sizes} is equality if, for example, $T$ is an \emph{$i$-path}, that is, we can order the edges of $T$ as $X_1,\dots,X_i$ so that, for each $j\in [i-1]$, the intersection $X_{j+1}\cap (\cup_{s=1}^j X_s)$ consists of exactly one pair of vertices and this pair belongs to $\p{1}{X_j}\setminus \p{1}{\cup_{s=1}^{j-1} X_s}$.

The following result shows that the \OneCluster{}s of any $\cG{r}{k}$-free graph have a very simple structure: namely, they are all small trees.

\begin{lemma}\label{lm:Trees} 
With the above notation, if $G$ is $\cG{r}{k}$-free, then every $F\in\OneClusters$ is an $m$-tree for some $m\in [k-1]$.
\end{lemma} 
\bpf
 Since $F$ is connected, Corollary~\ref{cr:Trim} gives an ordering $X_1,\dots,X_m$ of its edge set such that, for each $i\in [2,m]$, the $i$-th edge $X_{i}$ shares at least two vertices with some earlier edge.
 By induction, the number of vertices spanned by $\{X_1,\dots,X_i\}$ is at most $i(r-2)+2$ for each $i\in [m]$. Since $F$ is $(k(r-2)+2,k)$-free, we have $m<k$. Furthermore, for each $i\in [2,m]$, the edge $X_i$ has at least $r-2$ vertices not present in the previous edges by the $(i(r-2)+1,i)$-freeness of~$G$. It follows that $F$ is an $m$-tree, as desired.
\epf

\subsection{\texorpdfstring{Upper bound for $(r,k)= (3,5)$}{}}\label{se:35}

We begin with the case $(r,k)= (3,5)$, which is simpler but still embodies some key ideas that also apply to higher uniformities. 

\vspace{12pt}
\bpf[Proof of the upper bound of Theorem~\ref{5edge} for $r=3$]
By Lemma~\ref{chong}, it is enough to prove that $|G|\leq n^2/5$ for every $3$-graph on $n$ vertices which is $\cG{3}{5}$-free, that is, contains no $(7,5)$, $(5,4)$ and $(4,3)$-configurations.
Recall that $\OneClusters$ denotes the partition of $E(G)$ into \OneCluster{}s.
By Lemma~\ref{lm:Trees}, each \OneCluster\ is an $i$-tree with $i\leq 4$ edges.

For an $i$-tree $F$, consider the difference $2\,|\p{1} F|-5\,|F|=2(2i+1)-5i=2-i$, which is non-negative when $i\in \{1,2\}$. 
For $i\in \{3,4\}$, we would like to take an extra pair (in addition to those in $\p{1} F$) into account to make the difference non-negative. 
The following combinatorial lemma suffices for this.

\begin{claim}\label{cl:Deficit}
For every $i$-tree $F\in\OneClusters$ with $i\in\{3,4\}$, there is a pair which is \OutIn{1}{2}-claimed by $F$ but neither \OneClaim{}ed nor \OutIn{}{2}-claimed by $G\setminus F$.
\end{claim}
\bpf[Proof of Claim~\ref{cl:Deficit}] 
Note that no pair $xy$ can be \OutIn{}{2}-claimed by both $F$ and $G\setminus F$, for otherwise we can find a $3$-subtree in $F$ (by Corollary~\ref{cr:Trim}) \OutIn{}{3}-claiming $xy$. This would mean that $5\in \CI_{G}(xy)$ by~\eqref{eq:AB}, which contradicts Lemma~\ref{lm: 2}.

If $i=3$, then $F$ \OutIn{1}{2}-claims at least 2 pairs and at least one of them is not \OneClaim{}ed by $G\setminus F$: indeed,
if they are \OneClaim{}ed by different edges, then we get a 
$(7,5)$-configuration, and if they are \OneClaim{}ed by the same edge, then we get a $(5,4)$-configuration, a contradiction in either case. If $i=4$, then $F$ \OutIn{1}{2}-claims at least 3 pairs and, in fact, none can be \OneClaim{}ed by $G\setminus F$ (as otherwise we would have a $(7,5)$-configuration).
\epf

Now, for each $i$-tree $F\in \OneClusters$, define $P'_1(F)$ to be $\p{1} F$ if $i\in \{1,2\}$, and to be  $\p{1} F$ plus the pair returned by Claim~\ref{cl:Deficit} if $i\in\{3,4\}$. 

The sets $P'_1(F)$ for $F\in\OneClusters$ are pairwise disjoint and satisfy  $2\,|P'_1(F)|\geq 5\,|F|$. Thus,
$$
|G|=\sum_{F\in\OneClusters} |F|\leq \frac25 \sum_{F\in \OneClusters} |P'_1(F)| \leq \frac25 {n\choose 2}< \frac{n^2}5,
$$
giving the desired.
\epf

\subsection{Upper bounds of Theorems~\ref{5edge} and~\ref{7edge}}

In this section, we shall prove the remaining upper bounds of Theorems~\ref{5edge} and~\ref{7edge}, that is, the cases when $k=5$ and $r\geq 4$, or $k=7$ and $r\geq 3$. First, let us present some definitions and auxiliary results that are common to the proofs of both theorems.

Let $k\in \{5,7\}$ and let $G$ be an arbitrary $\cG{r}{k}$-free $r$-graph. Recall that $\C M_{\Merge{}{1}}$ denotes the partition of $E(G)$ into \OneCluster{}s. 

We say that a subgraph $H\subseteq G$ \emph{\TwoPClaim{}s} a pair $uv\in {V(G)\choose 2}$ if $H$ has a subtree $T$ with 3 edges that $\{2,3\}$-claims $uv$ (which, by Corollary~\ref{cr:Trim}, is equivalent to $T$ $2$-claiming the pair~$uv$).
Let us say that two edge-disjoint subgraphs $F,H\subseteq G$ are \emph{$\OTwoPMerge$-mergeable (via a pair $uv$)} if $uv$ is \OneClaim{}ed by $F$ and is \TwoPClaim{}ed by~$H$. 
If the order of $F$ and $H$ does not matter, we just say \emph{$\TwoPMerge$-mergeable}. 
Let $\TwoPClusters$ be the partition of $E(G)$ obtained from $\OneClusters$ by iteratively and as long as possible taking two $\TwoPMerge$-mergeable elements $F$ and $H$ and replacing them by $F\cup H$.
Let $\TwoPClusters'$ be the smallest family of subgraphs of $G$ that contains all \OneCluster{}s and is closed under $\TwoPMerge$-merging. We call the elements of $\TwoPClusters'$
\emph{partial \TwoPCluster{}s}. By the monotonicity of the merging rule, the family $\TwoPClusters$ consists of the maximal elements of $\TwoPClusters'$ and does not depend on the order of the merging steps.

Note that the relations 
of 
being $\OTwoPMerge$-mergeable and $\OMerge{\{1\}}{\{2,3\}}$-mergeable, while having many similarities, differ in general (e.g.\ a subgraph with 3 edges $3$-claiming a pair need not be a tree). 
As far as we see, the latter relation can also be used to prove the upper bounds for $k\in\{5,7\}$ but the former is more convenient for us to work with. 

If a partial \TwoPCluster{} $F$ \TwoPClaim{}s a pair $uv$ as witnessed by a 3-tree $T\subseteq F$, then trivially $T$ is a subgraph of one of the \OneCluster{}s in~$F$ and this \OneCluster{} \TwoPClaim{}s the pair~$uv$. Since the analogous statement for \OneClaim{}ed pairs is trivial, we have the analogue of Assumption~\eqref{eq:Trim} for $\OTwoPMerge$-mergeability (with $\C P = \OneClusters$). The proof of Lemma~\ref{lm:Trim} with the obvious modifications works for partial \TwoPCluster{}s. 
We will need only the following special case.

\begin{lemma}\label{lm:Cluster} 
For every $F\in\TwoPClusters'$ and $T\in \OneClusters$ such that $T\subseteq F$, 
there is an ordering $T_1,\dots,T_t$ of all \OneCluster{}s in $F$ such that $T_1=T$ and, for every $i\in [t-1]$, 
$T_{i+1}$ and $\bigcup_{j=1}^{i} T_j$ are $\TwoPMerge$-mergeable. (In particular, $\bigcup_{j=1}^{i} T_j\in\TwoPClusters'$ for every $i\in [t]$.)\qed
\end{lemma}

We say that an $r$-graph $F$ \emph{\OutIn{1}{2^+}-claims} a pair $uv$ if it \TwoPClaim{}s  but not $1$-claims $uv$, and denote the set of all such pairs by $\pp{1}{2^+}{F}$. Note that, if we build \TwoPCluster{}s by starting with \OneCluster{}s and merge some $\OTwoPMerge$-mergeable $F$ and $H$ via $uv$, then the pair $uv$ is \OutIn{1}{2^+}-claimed by~$H$. 
An integer sequence $(e_1,\dots,e_t)$ is called a \emph{composition} of a partial \TwoPCluster\ $F\in \TwoPClusters'$ if there is a sequence $(T_1,\dots,T_t)$ as in Lemma~\ref{lm:Cluster} with $|T_i|=e_i$ for each $i\in [t]$. Of course, every two compositions of the same \TwoPCluster\ are permutations of each other. \emph{The (non-increasing) composition} of $F$ is the non-increasing reordering of $(e_1,\dots,e_t)$; in general, there need not be a sequence of iterative legal merges realising it.

\vspace{12pt}
\bpf[Proof of the upper bound of Theorem~\ref{5edge} for $r\geq 4$]
By Lemma~\ref{chong}, it is enough to bound the size of any $\cG{r}{5}$-free $r$-graph $G$ on $[n]$ from above.
As before, $\OneClusters$ (resp.\ $\TwoPClusters$) denotes the set of  \OneCluster{}s (resp.\ \TwoPCluster{}s) of~$G$. 
By Lemma~\ref{lm:Trees}, each  \OneCluster\ is an $i$-tree with $i\in [4]$.  

Let every \TwoPCluster\ $F\in\TwoPClusters$ assign weight $1$ to every pair in $\p{1} F$ and in $\pp{1}{2^+} F$. Note that every pair $xy$ receives weight at most $1$. 
Indeed, suppose for contradiction that $xy$ receives weight 1 from two different \TwoPCluster{}s $F$ and $H$. 
Then, $xy$ must be \OutIn{1}{2^+}-claimed by both (as otherwise $F$ and $H$ would be merged together). 
However, this implies that $\{2,3\}\subseteq \CI_{F}(xy)\cap \CI_{H}(xy)$, which contradicts Lemma~\ref{lm: 2}.

Thus, in order to prove the theorem, it is enough to show that, for every $F\in\TwoPClusters$, we have $\lambda(F)\geq 0$ where
$$
 \lambda(F):=2\left(|\p{1} F|+|\pp{1}{2^+} F|\right)-(r^2-r-1)|F|.
$$
Indeed, we will then have that
\beq{eq:End1}
|G|=\sum_{F\in \TwoPClusters} |F|\leq \frac{2}{r^2-r-1}\,(|\p{1} F|+|\pp{1}{2^+} F|)\leq \frac{2}{r^2-r-1}{n\choose 2}.
\eeq

\begin{claim}\label{cl:35}
Every $F\in\TwoPClusters$ has at most four edges.
\end{claim}
\bpf[Proof of Claim~\ref{cl:35}]
Suppose for the sake of contradiction that $|F|\geq 5$ (and thus $|F|\geq 6$ since $F$ contains at most $(r-2)|F|+2$ vertices and $G$ is $(5r-8,5)$-free). Let $F$ be obtained by merging \OneCluster{}s $F_1,\dots,F_m\in\OneClusters$ in this order as in Lemma~\ref{lm:Cluster}.
Let us stop the merging process when we reach a partial \TwoPCluster\ containing at least 5 (and thus at least 6) edges. 
Suppose that  we have merged $F_1,\dots,F_s$ until this point. By  Lemma~\ref{lm:Trees}, we have $|F_s|\le 4$ and thus $s\ge 2$.

We claim that the last tree $F_s$ has exactly three edges.
As $G$ is $(5r-8,5)$-free, it holds that $|\bigcup_{i=1}^{s-1} F_i|\neq 5$. Thus, $F_s$ has at least two edges. 
In fact, $F_s$ cannot have exactly two edges as otherwise the subgraph  $\bigcup_{i=1}^{s-1} F_i$ of size 4 would \OutIn{}{4}-claim a pair in $\p{1}{F_s}$, contradicting Lemma~\ref{lm: 2}. 
Also, the tree $F_s$ cannot have 4 edges as otherwise any $\TwoPMerge$-merging involving it would lead to a $(5r-8,5)$-configuration and we would have $s=1$, a contradiction. 

Next, we can build the partial \TwoPCluster{} $\bigcup_{i=1}^s F_i$ by starting with $H_1:=F_s$ as in Lemma~\ref{lm:Cluster}; let us stop here at the first moment when the current partial \TwoPCluster\ $H$, say composed of $H_1,\dots,H_t\in\{F_1,\dots,F_{s}\}$ in this order, has at least five edges. 
As before, we have that $|H_t|= 3$. 
By $(5r-8,5)$-freeness, the sizes of the \OneCluster{}s in $H$ in the order of merging are either $(3,3)$ or $(3,1,3)$.
In the former (resp.\ latter) case, 
we can remove an edge from one (resp.\ each) of the 3-trees $H_1$ and $H_t$ so that the remaining subgraph is a diamond \OutIn{1}{2}-claiming the pair along which this tree attaches to the rest of~$H$.
This way, we can find a sequence of trees of sizes $(3,2)$ or $(2,1,2)$ with each one containing some 2 previously used vertices. 
This gives a forbidden $5$-edge configuration, proving that $|F|\leq 4$ for every $F\in\TwoPClusters$.\epf

It follows that, if $F\in\TwoPClusters$ is not a tree, then $F$ is made of a single edge and a 3-tree which must share exactly two vertices (for otherwise a $(4r-7,4)$-configuration appears). In this case, we have that $|\p{1} F|=4{r\choose 2}-2$ and $|\pp{1}{2^+} F|\geq 2(r-2)^2-1$. Thus,
$$
\lambda(F)\geq 2\left(4{r\choose 2}-2+2(r-2)^2-1\right)-(r^2-r+1)\cdot 4=4r^2-16r+6,
$$
which is positive  for $r\geq 4$. 
 
Let $F\in\TwoPClusters$ be an $i$-tree. If $i\geq 3$, then $|\pp{1}{2^+}{F}|\geq (i-1)(r-2)^2$ and  we have
 \begin{eqnarray*}
 \lambda(F)&\ge& 2\left(i{r\choose 2}-i+1+(i-1)(r-2)^2\right)-(r^2-r-1)i\\
 &=&i\left(2r^2-8r+7\right)-(2r^2-8r+6).
 \end{eqnarray*}
For $r\geq 4$, this expression is monotone increasing in $i$ and thus is at least its value when $i=3$, which is
$4r^2-16r+15\geq 15$. Finally, if $i=1,2$, then $\lambda(F)$ is respectively $2{r\choose 2}-(r^2-r-1)=1$ and $2\left(2{r\choose 2}-1\right)-2(r^2-r-1)=0$. 

Hence, for all $F\in\TwoPClusters$, we have that $\lambda(F)\geq 0$, which finishes the proof of the theorem by~\eqref{eq:End1}.
\epf

\bpf[Proof of the upper bound of Theorem~\ref{7edge}]
By Lemma~\ref{chong}, it is sufficient to bound the size of any $\cG{r}{7}$-free $r$-graph $G$ of order~$n$ from above. As before, $\OneClusters$ (resp.\ $\TwoPClusters$) is the set of all \OneCluster{}s (resp.\ \TwoPCluster{}s) of $G$. By Lemma~\ref{lm:Trees}, each element of $\OneClusters$ is an $i$-tree with $i\in [6]$.

Let each \TwoPCluster\ $F\in \TwoPClusters$ assign weight $1$ to each pair it \OneClaim{}s and weight $1/2$ to each pair  it \OutIn{1}{2^+}-claims. 
Let us show that every pair $xy\in \binom{V(G)}{2}$ receives weight at most~$1$. 
If a pair of vertices is \OneClaim{}ed by some \TwoPCluster{}, then it cannot be \OneClaim{}ed or \TwoPClaim{}ed by another \TwoPCluster{} as otherwise this would violate the merging rules for $\OneClusters$ or $\TwoPClusters$. 
Furthermore, a pair of vertices $xy$ cannot be $\{2,3\}$-claimed by three \TwoPCluster{}s by Lemma~\ref{lm: 2}.
So $xy$ indeed receives weight at most~1.

Thus, in order to prove the upper bound, it is sufficient to show that each \TwoPCluster\ $F$ satisfies that
\begin{equation}
 \lambda(F):=2w(F)-   (r^2-r-1)|F|\geq 0,
    \label{E7edge}
\end{equation}
where $w(F):=|\p{1} F|+\frac{1}{2}|\pp{1}{2^+} F|$ is the total weight assigned by $F$ to the vertex pairs. Indeed, we would then be done since 
 $$
 |G|=\sum _{F\in \TwoPClusters}|F|\leq \sum_{F\in \TwoPClusters}\frac{2}{r^2-r-1}\,w(F)\leq \frac{2}{r^2-r-1}\binom{n}{2}.
 $$

To show $(\ref{E7edge})$, we first prove the following claim.

\begin{claim}\label{cl:str7}
For each $F\in \TwoPClusters$, we have $|F|\leq 6$.
\end{claim}
\bpf[Proof of Claim~\ref{cl:str7}]
Suppose for contradiction that  $F\in \TwoPClusters$ has at least $7$ edges. Since $F$ is $(7r-12,7)$-free, we have $|F|\geq 8$. Let $F$ be obtained by merging $m$ distinct \OneCluster{}s $T_1,\dots,T_m\in\OneClusters$ in this order as in Lemma~\ref{lm:Cluster}.

Let $s\in [m]$ be the first index satisfying $|\bigcup_{i=1}^{s} T_i|\geq 8$. Then, $|\bigcup_{i=1}^{s-1} T_i|\leq 6$ as $F$ is $(7r-12,7)$-free. 
Hence, $|T_s|\geq 2$. It is impossible that $T_s$ and $T_1\cup \ldots\cup T_{s-1}$ are $\OTwoPMerge$-mergeable as otherwise we could trim edges from $T_s$ using Corollary~\ref{cr:Trim} to get a $(7r-12,7)$-configuration in $F$, a contradiction. Thus, $T_s$ \TwoPClaim{}s some pair $xy$ \OneClaim{}ed by $\bigcup_{i=1}^{s-1} T_i$; in particular, $|T_s|\geq 3$. 
Let $D$ be the diamond in $T_s$ \OutIn{}{2}-claiming $xy$. We note that $|(\bigcup_{i=1}^{s-1} T_i)\cup D|\geq 8$ since otherwise we could trim edges in $T_s\setminus D$ using Corollary~\ref{cr:Trim} to obtain a $(7r-12,7)$-configuration. 
Thus, $|\bigcup_{i=1}^{s-1} T_i|=6$.
Let $T_1':=T_s$ and let $T_2'$ be any \OneCluster{} in $\{T_1,\dots,T_{s-1}\}$ which is $\TwoPMerge$-mergeable with $T_1'$. By Lemma~\ref{lm:Cluster}, we can obtain the partial \TwoPCluster\ $\bigcup_{i=1}^{s-1} T_i$ by starting with 
$T_2'$ and $\TwoPMerge$-merging the remaining \OneCluster{}s one at a time, say $T_3',\dots,T_s'$ in this order.
Let $t\in [s]$ be the smallest index such that $T_1'\cup\cdots\cup T_t'$ has at least $7$ edges. By the same argument as before, we derive that
$|\bigcup_{i=1}^{t-1}T'_i|=6$ and $|T_t'|\ge3$.
Also, we can trim edges one by one from each of the ``pendant'' \OneCluster{}s $T_1'$ and $T_t'$ down to a diamond so that each intermediate subgraph is always a tree that shares 2 vertices with the partial \TwoPCluster{} $F':=\bigcup_{i=2}^{t-1} T_i'$. Since $|F'|\leq 6-|T_1'|\leq 3$, we must encounter a $(7r-12,7)$-configuration inside $\bigcup_{i=1}^t T_i'$ in this process, a contradiction.
\epf

Now, we prove $(\ref{E7edge})$ for every \TwoPCluster\ $F\in \TwoPClusters$. 

Take any $\TwoPMerge$-merging sequence $T_1,\dots ,T_m$ for $F$ as in Lemma~\ref{lm:Cluster}. For $i\in [m]$, let $e_i:=|T_i|$. Thus, $(e_1,\dots,e_m)$ is a composition of~$F$.
Then,
 $$
 |\pp{1}{2^+} F|\geq 1-m+\sum_{e_i\geq 3}(e_i-1)(r-2)^2,$$
since each $T_i$ shares exactly two vertices with $\bigcup_{j=1}^{i-1}T_j$ (as otherwise there would be an $(r\ell-2\ell+1,\ell)$-configuration in $G$ for some $\ell \in [2,6]$). Thus, by~\eqref{eq:Sizes}, we have that 
\begin{eqnarray}
\lambda(F)&\geq & 2\sum_{i=1}^m\left(e_i{r\choose 2}-e_i+1\right)
+1-m+\sum_{e_i\geq 3} (e_i-1)(r-2)^2- (r^2-r-1)\sum_{i=1}^m e_i\nonumber\\
&=&1+(r-2)^2\sum_{e_i\geq 3} (e_i-1)- \sum_{i=1}^m (e_i-1).\label{eq:temp1}
\end{eqnarray}

Our goal is to show that~\eqref{eq:temp1} is non-negative.
Let us denote by $x = x(F)$ the number of diamonds in the merging sequence of $F$, that is, the number of $i\in[m]$ with $e_i=2$. Note that $x\leq 1$: indeed, if $m\geq 2$, then $\max(e_1,e_2)\geq 3$ (in order for the first merging to occur) and Claim~\ref{cl:str7} implies that $x\leq \lfloor\tfrac{6-3}{2}\rfloor = 1$. Since the contribution of each $e_i\not=2$ to the right-hand side of~\eqref{eq:temp1} is non-negative, the expression there is at least $1-x\geq 0$, as desired.
\epf

\subsection{\texorpdfstring{Upper bounds for $k=6$}{}}

Here we set $k=6$. We will continue using the definitions of Section~\ref{se:UpperCommon} for a given $\cG{r}{6}$-free $r$-graph~$G$. In particular, recall that $\OneClusters$ denotes the set of \OneCluster{}s of~$G$.
However, unlike in the cases $k=5,7$, diamonds in the final partition would have to assign some positive weight to \OutIn{1}{2}-claimed pairs as otherwise the best we could hope for would be only $|G|\le(\frac1{r^2-r-1}+o(1))n^2$, which is strictly larger than the desired upper bound. This brings extra challenges to proving that each pair of vertices receives weight at most~1. We resolved this by using a different merging rule. Namely, in addition to the partition $\OneClusters$ of $E(G)$ into \OneCluster{}s, we will also use the partition
$$
  \TwoClusters:=\mathcal{M}_{\Merge{\{1\}}{\{2\}}}(\OneClusters),
  $$
which is obtained from the partition $\OneClusters$ by iteratively and as long as possible combining $\OMerge{}{2}$-mergeable pairs, that is, two current parts such that the first 1-claims and the second 2-claims the same pair. Also, we define $\TwoClusters'$ to consist of all subgraphs of $G$ that may appear at any stage of this process, calling  the elements of $\TwoClusters$ (resp.\ $\TwoClusters'$) \emph{\TwoCluster{}s} (resp.\ \emph{partial \TwoCluster{}s}). 
 
Let us observe some basic properties of (partial) \TwoCluster{}s. If some two partial \TwoCluster{}s $F$ and $H$ are  $\OMerge{}{2}$-mergeable via a pair $uv$ then it holds that $uv\not\in \p{1}{H}$ (as otherwise \OneCluster{}s of $F$ and $H$ \OneClaim{}ing the pair $uv$ would have been merged  when constructing~$\OneClusters$). Likewise, if a partial \TwoCluster\ $F$ $2$-claims $uv$ then the pair $uv$ is $2$-claimed by one of the \OneCluster{}s that make~$F$. Thus, Assumption~\eqref{eq:Trim} of Lemma~\ref{lm:Trim} holds and the conclusion of the lemma applies here.

For $F\in\TwoClusters'$ which is made by merging \OneCluster{}s $F_1,\dots,F_m$ in this order as in Lemma~\ref{lm:Trim}, we call the sequences of sizes $(|F_1|,\dots,|F_m|)$ a \emph{composition} of $F$. Its non-increasing reordering is called \emph{the (non-increasing) composition} of~$F$.

Also, recall that, by Lemma~\ref{lm: 2}, it holds for any edge-disjoint subgraphs $F_1,\dots,F_s\subseteq G$ and any $xy\in {V(G)\choose 2}$ that
\beq{eq:2}
6\notin \sum_{i=1}^s \CI_{F_i}(xy).
\eeq

Next, in Lemma~\ref{lm:6Comb} below, we derive some combinatorial properties that every partial \TwoCluster\ has to satisfy and that will suffice for our estimates. (An exact description of all possible partial \TwoCluster{}s is possible, with some extra work.) Since the proof of the lemma does not introduce any new ideas in addition to the ones seen before, the reader may skip it in the first reading.

\begin{lemma}\label{lm:6Comb} Let $r\geq 3$ and $G$ be an arbitrary $\cG{r}{6}$-free $r$-graph. Let a partial \TwoCluster\ $F$  be obtained by merging the elements $T_1,\dots ,T_m$ of $\OneClusters$ in this order as in Lemma~\ref{lm:Trim}. Let $(e_1,\dots,e_m)$ be the  composition of~$F$, that is, the non-increasing reordering of $(|T_1|,\dots,|T_m|)$. Then, each of the following statements holds.
 \begin{enumerate}[(a)]
 \item \label{it:>6} If $F$ has at least $7$ edges, 
 then $(e_1,\dots,e_m)$ is either $(3,2,2)$ or $(2,\dots,2,1)$ with at most $r(r-1)$ entries equal to~$2$. If, moreover, $T_1$ is the unique \OneCluster{} of size different from $2$ (that is, $|T_1|=1$ or $3$)
 then, for each $i\in [2,m]$,  $H_i:=\bigcup_{j=1}^{i-1} T_j$ and $T_i$ are $\OMerge{}{2}$-mergeable via some pair $xy\in \pp{1}{2}{T_i}$ while
 no other pair in $\pp{1}{2} {T_i}$ is \OneClaim{}ed or $2$-claimed by~$H_i$.
 \item\label{it:claim} It holds that $|\pp{1}{2} F  |\geq 1-m+\sum_{i=1}^m (e_i-1) (r-2)^2$.
 \item\label{it:2111} If $(e_1,\dots,e_m)=(2,1,1,1)$ then no pair in $\pp{1}{2} F$ is \OneClaim{}ed or $2$-claimed \mbox{by~$G\setminus F$}.
 \end{enumerate}
\end{lemma}
\bpf
Suppose first that the partial $\Merge{}{2}$-cluster $F$ has at least $7$ edges.
We prove the first two claims of the lemma for this~$F$ simultaneously.

Let $s\in [m]$ be the first index such that $|\bigcup_{i=1}^{s} T_i|\geq 7$. Then, $|H_s|\leq 5$ as $F$ is \mbox{$(6r-10,6)$}-free, and hence, $|T_{s}|\geq 2$. 
If $T_{s}$ and $H_s$ are $\OMerge{}{2}$-mergeable, then by Corollary~\ref{cr:Trim} we can remove some edges from $T_{s}$ to get a $(6r-10,6)$-configuration inside $H_s\cup T_s$, a contradiction. 
Hence, $T_{s}$ must \OutIn{1}{2}-claim some pair $xy$ \OneClaim{}ed by $H_s$. Let $D\subseteq T_s$ be the (unique) diamond \OutIn{1}{2}-claiming~$xy$.
Note that $|H_s\cup D|\geq 7$ as otherwise Corollary~\ref{cr:Trim} implies that we could remove some edges from $T_{s}\setminus D$ one by one to obtain a $(6r-10,6)$-configuration. Thus, $|H_s|=5$. 

Let $T_1':=T_s$ and let $T_2'\in \{T_1,\dots,T_{s-1}\}$ be a $\Merge{}{1}$-cluster $2$-mergeable with $T_s$. 
Let $(T_2',\dots,T_{s}')$ be the ordering of $\{T_1,\dots,T_{s-1}\}$ returned by Lemma~\ref{lm:Trim} for the partial  $\Merge{}{2}$-clusters $T_2'\subseteq \bigcup_{i=1}^{s-1} T_i$. By the choice of $T_2'$,  for each $i=2,\dots,s$ the \OneCluster\ $T_i'$ is $\Merge{}{2}$-mergeable with the partial \TwoCluster\ $\bigcup_{j=1}^{i-1}T_j'$. Let $t\in [s]$ be the first index such that $|\bigcup_{i=1}^{t} T_i'|\geq 7$. Set $H':=\bigcup_{i=1}^{t-1} T_i'$.
By the same argument as in the previous paragraph, we have that $|H'|=5$ and there is a diamond $D'\subseteq T_t'$ such that $H'$ and $D'$ are $\OMerge{}{2}$-mergeable.

Thus, we have a partial $\Merge{}{2}$-cluster $F':=\bigcup_{i=1}^{t} T_i'$  with at least $7$ edges built via the sequence $(T_1',T_2',\dots,T_t')$
so that the first $\Merge{}{1}$-cluster $T_1'=T_s$ (resp.\ the last $\Merge{}{1}$-cluster $T_t'$) can be merged with the rest through only one pair, which is \OutIn{1}{2}-claimed by the diamond $D\subseteq T_s$ (resp.\ $D'\subseteq T_t'$). Here we have the freedom to trim one or both of these two clusters, leaving any number of edges in each except exactly 1 edge. It routinely follows that $|T_1'|=|T_t'|=2$ (that is, $T_1'=D$ and $T_t'=D'$). For example, if $(|T_1'|,|T_t'|)=(3,3)$, then we can trim exactly one edge (from $T_1'$), two edges (one from each of $T_1'$ and $T_t'$) or three edges  (all of $T_1'$), with one of these operations leaving a forbidden subgraph of  $F'$  with exactly 6 edges.

Thus, the \OneCluster{}s $T_i'$ with $2\leq i\leq t-1$ contain exactly 3 edges in total. This leaves us with the following possibilities for the sequence $(e_1',\dots,e_t')$ of the encountered sizes $e_i':=|T_i'|$.

The first case that $(e_1',\dots,e_t')=(2,1,1,1,2)$ is in fact impossible because one of the 1-trees can be removed so that the remaining $r$-graph is a partial \TwoCluster\ with exactly $6$ edges, a contradiction.

If $(e_1',\dots,e_t')=(2,3,2)$ then, in order to avoid a forbidden configuration, the following statements must hold: 
$T_2'$ is a 3-path, the diamonds $T_1'$ and $T_3'$ \OutIn{1}{2}-claim pairs that are \OneClaim{}ed by the opposite end-edges of the $3$-path $T_2'$ but not by the middle edge, $|V(T_1')\cap V(T_2')|=|V(T_3')\cap V(T_2')|=2$, $|V(T_1')\cap V(T_3')|\leq 1$, 
no further \OneCluster\ can be merged with $F'=T_1'\cup T_2'\cup T_3'$ (so $F'=F$ by Lemma~\ref{lm:Trim}),
and therefore $F$ satisfies Part~\ref{it:>6} of the lemma. Part~\ref{it:claim} now easily follows.

Finally, suppose that $(e_1',\dots,e_t')$ is $(2,1,2,2)$ or $(2,2,1,2)$, with the single-edge \OneCluster{} being $\{X\}$. By two applications of Lemma~\ref{lm:Trim} (for $\{X\}\subseteq F'$ and for $F'\subseteq F$), we can additionally assume that $F'$ is made of $T_1=\{X\}$ and three diamonds $T_2$, $T_3$ and $T_4$ (and thus $F$ can be obtained from $F'$ by iteratively merging \OneCluster{}s $T_5,\dots, T_m$ in this order). 
Call a \OneCluster\ $T_i$ for $i\geq 2$ of \emph{type $ab$} if the merging chain connects it to $X$ via a pair $\{a,b\}\subseteq X$. (Note that the vertices $a,b$ are not necessarily in $T_i$: for example, $T_i$ can merge with a diamond 2-claiming $ab$.) 
By convention, we assume that the 1-tree $T_1$ is of all ${r\choose 2}$ types. 
Observe that at least two of the initial diamonds $T_2,T_3,T_4$ must be of different types (as otherwise $T_2\cup T_3\cup T_4$ would be a $(6r-10,6)$-configuration).

Let us continue denoting $H_i:=T_1\cup \dots\cup T_{i-1}$ for $i\in [m-1]$. In order to finish Part~\ref{it:>6}, it remains to prove the following.

\begin{claim}\label{cl:>6D}
For every $i\in [2,m]$, $T_i$ is a diamond that \OutIn{1}{2}-claims some previously \OneClaim{}ed pair $x_iy_i\in \p{1}{H_i}$, and no pair in $\pp{1}{2} {T_i}\setminus\{x_iy_i\}$ is \OneClaim{}ed or $2$-claimed by~$H_i$. Also, $m\leq 1 +2{r\choose 2}$.
\end{claim}
\bpf[Proof of Claim~\ref{cl:>6D}]
For the first part, we use induction on $i\in [m]$. 
To begin with, it is easy to check that the configuration on $T_1, T_2, T_3, T_4$ satisfies the statement of the claim.
Let $i\in [5,m]$ and let the \OneCluster\ $T_i$ be of type~$ab$. Note that we have at most two \OneCluster{}s of each type among the diamonds $T_2,\dots,T_{i-1}$ (otherwise, the first three diamonds of any given type would form a forbidden 6-edge configuration). If some edge $e'\in T_i$ \OneClaim{}s a pair \OutIn{1}{2}-claimed by $H_i$, then 
by keeping only the edge $e'$ in $T_i$ and removing one by one the diamonds of types different from $ab$, we can reach a partial \TwoCluster\ with exactly 6 edges, a contradiction. 
So let the diamond $D_i\subseteq T_i$ \OutIn{1}{2}-claim a pair $x_iy_i\in\p{1} {H_i}$. We know that there is at most one previous diamond $T_j$ of the same type as $T_i$ (otherwise $D_i$ with two such diamonds would form a $(6r-10,6)$-configuration). 
It follows that $D_i=T_i$ as otherwise a forbidden 6-edge configuration would be formed by $T_1$, $D_i$, some suitable edge of $T_i\setminus D_i$ plus either the diamond $T_j$ of type $ab$ (if it exists) or a diamond \OutIn{1}{2}-claiming a pair in~$\p{1} {T_1}\setminus \{ab\}$ (such diamond exists among $T_2, T_3, T_4$). If $D_i$ contains some other vertex $z_i\not\in\{x_i,y_i\}$ from an earlier \OneCluster\ of the same type~$ab$, then some edge $e'$ of $D_i$ shares at least two vertices with $H_i$, again leading to a forbidden $6$-vertex configuration in $H_i\cup \{e'\}$. 
Thus, we are done unless  $\pp{1}{2} {T_i}$ contains a pair $uv$ with both vertices in a \OneCluster{} of some different type~$a'b'$ (where $\{u,v\}$ may possibly intersect $\{x_i,y_i\}$). If each of the types $ab$ and $a'b'$ contains an earlier diamond (which then must be unique), then these two diamonds and $T_i$ form a configuration on 6 edges and at most $6r-10$ vertices; otherwise, we have in total at most two diamonds of Type $ab$ or $a'b'$ (including $T_i$) and these diamonds together with $T_1$ have $5$ edges and at most $5(r-2)+1$ vertices, a contradiction.

Finally, the inequality $m\leq 1+2{r\choose 2}$ follows from the observation made earlier that for every type $ab\in {e\choose 2}$ there are at most 2 diamonds among $T_2,\dots,T_m$ of this type.
\epf

Claim~\ref{cl:>6D} implies that the addition of each new diamond gives $(r-2)^2-1$ new \OutIn{1}{2}-claimed pairs, 
from which Part~\ref{it:claim} follows in the case $|F|\geq 7$.

Now, suppose that $|F|\leq 6$, and thus $|F|\leq 5$. When we construct $F$ by merging \OneCluster{}s one by one as in Lemma~\ref{lm:Trim}, each new \OneCluster\ shares exactly 2 vertices with the current partial 
\TwoCluster{} (since $G$ is $\cG{r}{6}$-free). Thus, Part~\ref{it:claim} follows. 

Finally, if  $(e_1,\dots,e_m)=(2,1,1,1)$, then $F$ consists of a diamond $D$ with three single edges $X_1,X_2,X_3$ attached along some pairs \OutIn{1}{2}-claimed by $D$. Take any pair $xy\in \pp{1}{2} F$; then, $xy$ is \OutIn{1}{2}-claimed by $D$ but not \OneClaim{}ed by any~$X_i$. Note that, for an edge $X\in G\setminus F$ \OneClaim{}ing $xy$ (resp.\ a diamond $D'\subseteq G\setminus F$ \OutIn{1}{2}-claiming $xy$), $F\cup\{X\}$ (resp.\ $(F\cup D')\setminus\{X_1\}$) is a $(6r-10,6)$-configuration. 
Thus, $G\setminus F$ neither \OneClaim{}s nor 2-claims $xy$, finishing the proof of Part~\ref{it:2111}. \epf

\subsubsection{Upper bound of Theorem~\ref{6edgehigh}}	

Here, we deal with the case $k=6$ and $r\geq 4$.

\vspace{12pt}
\bpf[Proof of the upper bound of Theorem~\ref{6edgehigh}]
By Lemma~\ref{chong}, it is enough to upper bound the size of a  $\cG{r}{6}$-free $r$-graph $G$ with $n$ vertices. 
Recall that, by Lemma~\ref{lm:Trees}, each element of $\OneClusters$ is an $i$-tree with $i\in [5]$.

Now, we define the weights. 
A \TwoCluster\ $F\in \TwoClusters$ assigns weight $1$ to each pair in $\p{1}{F}$
and weight $1/2$ to each pair in $\pp{1}{2}{F}$, except if the composition of $F$ is $(2,1,1,1)$ in which case every pair \OutIn{1}{2}-claimed by $F$ receives weight $1$ (instead of $1/2$).
Let us show that every pair $xy$ of vertices of $G$ receives weight at most~$1$. This is clearly true if there is a  \TwoCluster\ $F$ with the composition $(2,1,1,1)$ such that $xy\in \pp{1}{2} F$: then, $F$ gives weight $1$ to $xy$ and no other \TwoCluster\ \OneClaim{}s or 2-claims $xy$ by Part~\ref{it:2111} of Lemma~\ref{lm:6Comb}. Suppose that $xy$ is not \OutIn{1}{2}-claimed by any \TwoCluster{} with the composition $(2,1,1,1)$.
If $xy$ is \OneClaim{}ed by some \TwoCluster\ $F_1$, then it cannot be \OneClaim{}ed or 2-claimed by another \TwoCluster\ $F_2$ (as otherwise $F_1$ and $F_2$ would be merged). On the other hand, the pair $xy$ can be \OutIn{}{2}-claimed by at most two \TwoCluster{}s by \eqref{eq:2}.
In all cases, the pair $xy$ receives weight at most~1.

Let us
show that for each \TwoCluster\ $F\in \TwoClusters$, we have 
\begin{equation}
   \lambda(F):=2 w(F)- r(r-1)|F|\geq 0,
    \label{E6edge}
\end{equation}
where $w(F)$ denotes the total weight assigned by the \TwoCluster~$F$. First, consider the exceptional case when $F$ is composed of a diamond and 3 single edges. Here, $w(F)$ does not depend on how the three edges are merged with the diamond and we have
$$
\lambda(F)=
2\left(5{r\choose 2} + (r-2)^2-4\right)-r(r-1)\cdot 5=2((r-2)^2-4)\geq 0.
$$
(Note that, if $F$ gave weight of $1/2$ to each \OutIn{1}{2}-claimed pair, then~\eqref{E6edge} may be false for $r=4$, so some exceptional weight distribution is necessary.)

So let $F$ be any other (non-exceptional) \TwoCluster. For $j\in\I N$, let $n_j$ denote the number of \OneCluster{}s in $F$ with $j$ edges. Thus, $n_j=0$ for $j\geq 6$ by Lemma~\ref{lm:Trees}. 
We have
\begin{equation*}
\begin{aligned}
 &\lambda(F)= 2\,|\p{1} F|+|\pp{1}{2} F|-r(r-1)|F|\\
 &\geq  
 2\sum_{j=1}^5 \left(j{r\choose 2}-j+1\right)n_j + \left(1-\sum_{j=1}^5 n_j + \sum_{j=1}^5  (j-1) (r-2)^2n_j\right)-r(r-1) \sum_{j=1}^5 jn_j\\
 &=1+\sum_{j=1}^5  \left((r^2-4r)(j-1)+2j-3\right)n_j,
 \end{aligned}
\end{equation*}
where the inequality in the middle follows from Part~\ref{it:claim} of Lemma~\ref{lm:6Comb}. Since $r\geq 4$, the coefficient at $n_j$ is at least $2j-3$. This is negative only if $j=1$. Thus, $\lambda(F)\geq 0$ unless $n_1\geq 2$. By Part~\ref{it:>6} of Lemma~\ref{lm:6Comb} (and since we have already excluded the exceptional $(2,1,1,1)$-case of Part~\ref{it:2111}), this is only possible if $F$ has the composition $(3,1,1)$ or $(2,1,1)$. The corresponding sequences of $(n_1,n_2,n_3)$ are $(2,0,1)$ and $(2,1,0)$; thus, the corresponding values of $\lambda(F)$ are $2$ and $0$. Hence, $\lambda(F)\geq 0$ for every \TwoCluster~$F$, so the familiar double counting argument implies that $|G|\leq {r\choose 2}^{-1} {n\choose 2}$, proving the theorem.\epf

\subsubsection{Upper bound of Theorem~\ref{6edge3uniform}}

In this section, we deal with the case $(r,k)=(3,6)$.

\vspace{12pt}
\bpf[Proof of the upper bound of Theorem~\ref{6edge3uniform}]
By Lemma~\ref{chong}, it is enough to provide a uniform upper bound on the size of an arbitrary $\cG{3}{6}$-free $3$-graph $G$ on $V:=[n]$ from above.

As before, $\OneClusters$ (resp.\ $\TwoClusters$) denotes the partition of $E(G)$ into \OneCluster{}s (resp.\ \TwoCluster{}s).
Call edge-disjoint subgraphs $F,F'\subseteq G$ \emph{$\OThreePMerge$-mergeable (via $uv\in {V\choose 2}$)} if they are $\OMerge{\{1\}}{\{3,4\}}$ or $\OMerge{\{1,2\}}{\{3\}}$-mergeable via $uv$, that is, if at least one of the following two conditions holds:
\begin{itemize}
	\item $1\in \CI_{F}(uv)$ and $\{3,4\}\subseteq \CI_{F'}(uv)$, or
	\item $\{1,2\}\subseteq \CI_{F}(uv)$ and $3\in \CI_{F'}(uv)$.
\end{itemize}
 If the order of $F$ and $F'$ does not matter, then we simply say that they are \emph{$\ThreePMerge$-mergeable}.
 Let the partition $\ThreePClusters$ of $E(G)$ be obtained by starting with $\TwoClusters$ and, iteratively and as long as possible, merging any two $\ThreePMerge$-mergeable parts. 
Also, let $\ThreePClusters'$ be the set of $3$-graphs that could appear at some point of the above process.
We refer to the elements of $\ThreePClusters$ (resp.\ $\ThreePClusters'$) as \emph{\ThreePCluster{}s} (resp.\ \emph{partial \ThreePCluster{}s}). By monotonicity, the partition $\ThreePClusters$ does not depend on the order in which we perform the merging steps.

Let us observe some basic properties of $\ThreePClusters$.

\begin{lemma}\label{lm:3Merge} 
Suppose that the edge-disjoint partial \ThreePCluster{}s $F$ and $H$ are $\OThreePMerge$-mergeable via some pair $uv$. Then, there are \TwoCluster{}s $F'\subseteq F$ and $H'\subseteq H$ that are $\OThreePMerge$-mergeable via~$uv$.
\end{lemma}
\bpf Let $F'\subseteq F$ be the (unique) \TwoCluster{} that \OneClaim{}s the pair~$uv$. Let $H''$ be a $(5,3)$-configuration in $H$ that $3$-claims the pair~$uv$. Note that the pair $uv$ is not 2-claimed by $H''$ since otherwise the $2$-cluster in $H''$ claiming this pair would have been merged with $F'$.
Since $H''\subseteq G$ is $\cG{3}{6}$-free, $H''$ is either a 3-tree or the union of a single edge and a diamond that can be $\OMerge{}{2}$-merged. Thus, $H''$ lies entirely inside some \TwoCluster{} $H'$. Of course, $H'\subseteq H$. Let us show that $F'$ and $H'$ satisfy the lemma.

If $\{1,2\}\subseteq \CI_{F}(uv)$, as witnessed by an edge $e$ and a diamond $D$ in $F$, then $e$ is an edge of $D$ (as otherwise $D\cup\{e\}\cup H''$ would be a forbidden 6-edge configuration) and the lemma is satisfied (since $F'$ must contain $D$ as a subgraph).

So suppose that $4\in \CI_{H}(uv)$. We are done if $4\in \CI_{H'}(uv)$ so suppose otherwise. This assumption implies that no pair in $\pp{1}{2} {H''}$ can be \OneClaim{}ed by an edge from $G\setminus H''$. Furthermore, no pair in $\p{1}{H''}$ can be $2$-claimed by a diamond $D$ in $G\setminus H''$ as otherwise $D\cup H''$ together with an edge of $F'$ \OneClaim{}ing the pair $uv$ would form a $(8,6)$-configuration. 
Hence, $H'=H''$. Since $\CI_{H}(uv)\ni 4$ is strictly larger than $\CI_{H'}(uv)$, the \TwoCluster{} $H'$ is $\ThreePMerge$-mergeable with some other \TwoCluster{} $H'''$ in~$H\setminus H'$ via some pair~$u'v'$. 
Note that $3\in \CI_{H'}(u'v')$ since $H'$ is a $(5,3)$-configuration, so it 3-claims every pair of vertices it contains. 
By~\eqref{eq:2}, $3\notin \CI_{H'''}(u'v')$, and hence $1\in \CI_{H'''}(u'v')$. 
It follows that $u'v'\not=uv$ as otherwise $uv$ is \OneClaim{}ed by $F$ and $H'''$, contradicting the merging rule for~$\OneClusters$.
As $H'$ has only $3$ edges, the definition of $\ThreePMerge$-mergeability gives that $2\in \CI_{H'''}(u'v')$. However, in that case the union of $H'$, a diamond $D$ in $H'''$ that \OutIn{}{2}-claims the pair $u'v'$, and an edge in $F'$ that \OneClaim{}s the pair $uv$ forms a $(8,6)$-configuration,
a contradiction.\epf

Lemma~\ref{lm:3Merge} provides us with an analogue of Assumption~\eqref{eq:Trim} of Lemma~\ref{lm:Trim} and the proof of Lemma~\ref{lm:Trim} trivially adapts to $\ThreePMerge$-merging. We will need only the following special case.

\begin{claim}\label{cl:3Trim} 
For every partial \ThreePCluster{} $F$ and any \TwoCluster{} $F_0\subseteq F$, there is an ordering $F_0,\dots,F_s$ of the \TwoCluster{}s constituting $F$ such that, for each $i\in [s]$, the $3$-graphs $F_i$ and $\bigcup_{j=0}^{i-1} F_j$ are $\ThreePMerge$-mergeable.\qed
\end{claim}

Now, we consider the following two functions $f$ and $h$ from subsets of $[5]$ to the reals. Namely, for $A\subseteq [5]$, we define
$$f(A) :=\left\{\begin{array}{ll}
{55}/{61} & \text{ if } A=\{1\}, \\
1 & \text{ if } A=\{1,x\} \text{ for some } x\in \{2,3\},\\
{55}/{61} & \text{ if } A=\{1,x\} \text{ for some } x\in \{4,5\},\\
{25}/{61} & \text{ if } A=\{2\}, \\
{36}/{61} & \text{ if } A=\{2,3\}, \\
1 & \text{ if } A=\{2,3,4\},\\
{1}/{2} & \text{ if } A=\{2,4\},\\
{6}/{61} & \text{ if } A=\{3\}, \\
{11}/{61} & \text{ if } A=\{3,5\}, \\
1 & \text{ if } A=\{3,4\},\\
0 & \text{ for all other sets } A\subseteq [5],
\end{array}\right.$$
 and 
 $$
  h(A):=\max\{f(A') \mid A'\subseteq A \}.
  $$
  Then, the function $h$ is clearly non-decreasing and satisfies that
 \beq{eq:HZero} 
 h(A)>0\ \ \Leftrightarrow\ \ A\cap \{1,2,3\}\not=\emptyset.
 \eeq
In the sequel, we abbreviate $h(\{i_1,\dots,i_s\})$ to $h(i_1,\dots,i_s)$.
 
Define the weight attributed to a pair $uv\in {V\choose 2}$ by a subgraph $F\subseteq G$ to be 
$$
 w_F(uv) := h([5]\cap \CI_F(uv)).
 $$
Moreover, we set $w(uv) :=\sum_{F\in \ThreePClusters} w_F(uv)$ to be the total weight received by a pair $uv$ from all \ThreePCluster{}s.

\begin{claim}\label{cl:wuv}
For every $uv\in {V\choose 2}$, it holds that $w(uv)\leq 1$.
\end{claim}
\bpf[Proof of Claim~\ref{cl:wuv}]
Fix $uv$ and let $F_1,\dots, F_s$ be all \ThreePCluster{}s with $w_{F_i}(uv)>0$. We have to show that $\sum_{i=1}^s w_{F_i}(uv)\leq 1$.
Note that each $\CI_{F_i}(uv)$ intersects $\{1,2,3\}$ by~\eqref{eq:HZero}. 
If $s=1$, then we are done (since $h(A)\leq 1$ for every $A\subseteq [5]$), so assume that $s\geq 2$.

The following cases cover all possibilities up to a permutation of $F_1,\dots,F_s$.

\case{1}{Assume that $\CI_{F_1}(uv)$ contains $1$.}%
Then, for every $j\in [2,s]$, we have that $1,2\notin \CI_{F_j}(uv)$ (as otherwise the corresponding \OneCluster{}s of $F_1$ and $F_j$ would be merged when building $\OneClusters$ or $\TwoClusters$) and it follows from~\eqref{eq:HZero} that $3\in \CI_{F_j}(uv)$. Since the subgraphs $F_1,\dots,F_s\subseteq G$ are edge-disjoint, \eqref{eq:2} implies that $s=2$. Furthermore, since $F_1$ and $F_2$ are not $\OThreePMerge$-mergeable, it holds that $4\notin \CI_{F_2}(uv)$ and $2\notin \CI_{F_1}(uv)$. By \eqref{eq:2}, $5\notin \CI_{F_2}(uv)$ and $3\notin \CI_{F_1}(uv)$. Thus,
$$
w(uv)=w_{F_1}(uv)+w_{F_2}(uv)\leq h(1,4,5)+h(3)= \frac{55}{61}+ \frac{6}{61}=1.$$

\case{2}{Assume that no $\CI_{F_i}(uv)$ contains $1$ but $\CI_{F_1}(uv)$ contains $2$.}%
Here, it is impossible to have distinct $i,j\in [2,s]$ with $2\in \CI_{F_i}(uv)\cap \CI_{F_j}(uv)$ as otherwise the edge-disjoint subgraphs $F_1, F_i, F_j\subseteq G$ would contradict \eqref{eq:2}. We split this case into 2 subcases.
 
\case{2-1}{Assume $2\in \CI_{F_2}(uv)$.}%
Suppose first that $s\geq 3$. Then, for every $j\in [3,s]$, we have $2\notin \CI_{F_j}(uv)$ by \eqref{eq:2} and thus $3\in \CI_{F_j}(uv)$ by~\eqref{eq:HZero}. It follows from \eqref{eq:2} that $s=3$ and, moreover,  $4\notin \CI_{F_3}(uv)$ and $3,4\notin \CI_{F_1}(uv)\cup \CI_{F_2}(uv)$. Hence
$$
w(uv)=w_{F_1}(uv)+w_{F_2}(uv)+w_{F_3}(uv)\leq 2\,h(2,5)+h(3,5)= 2\cdot  \frac{25}{61} + \frac{11}{61} = 1.
$$
If $s=2$, then \eqref{eq:2} implies that $4\notin \CI_{F_1}(uv)\cup \CI_{F_2}(uv)$ and $3\notin \CI_{F_1}(uv)\cap \CI_{F_2}(uv)$. Therefore, $$w(uv)=w_{F_1}(uv)+w_{F_2}(uv)\leq h(2,5)+h(2,3,5)=\frac{25}{61} + \frac{36}{61}=1.$$

\case{2-2}{Assume that $2\notin \CI_{F_j}(uv)$ for all $j\in[2,s]$.}%
By~\eqref{eq:HZero}, $\CI_{F_j}(uv)$ contains $3$ for every $j\in[2,s]$. By \eqref{eq:2}, it holds that $s=2$ and, moreover, $3\notin \CI_{F_1}(uv)$ and  $4\notin \CI_{F_2}(uv)$. Thus, we have
$$
w(uv)=w_{F_1}(uv)+w_{F_2}(uv)\leq h(2,4,5)+h(3,5)=\frac{1}{2} + \frac{11}{61} \leq 1.
$$

\case{3}{Assume that no $\CI_{F_i}(uv)$ contains $1$ or $2$.}%
By~\eqref{eq:HZero}, we have $3\in \CI_{F_i}(uv)$ for each $i\in [s]$. However, our assumption that $s\geq 2$ contradicts~\eqref{eq:2}. 
This finishes the case analysis and the proof.
\epf

Now, let us show that, for every $F\in \ThreePClusters$, the total weight 
$$
w(F): = \sum_{uv\in \binom{V}{2}} w_{F}(uv)
$$
assigned by $F$ to different vertex pairs is at least $\frac{165}{61}\,|F|$.

First, we check this for the \ThreePCluster{}s $F$ consisting of a single \OneCluster{}.

\begin{claim}\label{cl:wF}
For all $F\in \OneClusters$, we have $w(F)\geq \frac{165}{61}\,|F|$.
\end{claim}
\bpf[Proof of Claim~\ref{cl:wF}]
Recall that $F$ is an $i$-tree with $i\leq 5$ by Lemma~\ref{lm:Trees}. Assume that $i\geq 2$ as otherwise $w(F)=3h(1)=\frac{165}{61}$ and the claim holds.
    
Every pair in $\p{1} F$ is $\{1,2\}$-claimed (in fact, $\{1,\dots,i\}$-claimed by Corollary~\ref{cr:Trim}) and, in particular, receives weight at least $h(1,2)=1$ from~$F$. 
Moreover, since $F$ is an $i$-tree, it \OneClaim{}s $2i+1$ pairs. 
Then, if $i=2$, $w(F) = 5\,h(1,2)+h(2) =5+ \frac{25}{61} 
= \frac{165}{61}\cdot 2$, as desired.

Now, assume that $i\geq 3$. Then, each pair in $\pp{1}{2} F$ is $\{2,3\}$-claimed by Corollary~\ref{cr:Trim} and receives weight at least $h(2,3)=\frac{36}{61}$. 
Moreover, $|\pp{1}{2} F|\geq i-1$, and if we exclude all pairs in $\p{1} F$ and some $i-1$ pairs in $\pp{1}{2} F$ then, regardless of the structure of $F$, there will remain at least $i-2$ pairs that are $3$-claimed. 
Indeed, easy induction shows that any $i$-tree with $i\geq 3$ contains at least $2i-3$ different sub-paths of length $2$ or $3$ such that the opposite vertices of degree 1 in these paths give distinct pairs outside of $\p{1}{F}$ that are $3$-claimed by $F$.
Thus,
\begin{eqnarray*}
 w(F)&\geq& (2i+1)\,h(1,2)+(i-1)\,h(2,3)+(i-2)\,h(3)\\
  &=&(2i+1)+(i-1)\,\frac{36}{61}+(i-2)\,\frac{6}{61}\\
  &=& \frac{13-i}{61} + \frac{165}{61}\,i\ >\ \frac{165}{61}\,|F|,
\end{eqnarray*}
 as required.
\epf

As a next step, we estimate the weight assigned by those \TwoCluster{}s that consist of more than one \OneCluster{}.

\begin{claim}\label{cl:wFLower}
For all $F\in \TwoClusters\setminus \OneClusters$, we have 	$w(F)\geq \frac{165}{61}\,|F|$.
\end{claim}
\bpf[Proof of Claim~\ref{cl:wFLower}] Note that if $F,H\subseteq G$ are $\OMerge{}{2}$-mergeable via a pair $uv$, then $\{1\}+\CI_{H}(uv)\subseteq \CI_{F\cup H}(uv)$ and $\{2\}+ \CI_{F}(uv)\subseteq \CI_{F\cup H}(uv)$ holds. In particular, we conclude by \eqref{eq:2} that $5\notin \CI_{H}(uv)$ and $4\notin \CI_{F}(uv)$.

Suppose first that $|F|\geq 7$. By Lemma~\ref{lm:6Comb}\ref{it:>6}, there are two cases to consider. 
First, let $F$ be made from a 3-tree by $\OMerge{}{2}$-merging two diamonds one by one. Note that $F$ $\{1,2\}$-claims all $17$ pairs in $\p{1} F$. 
Since each new diamond $abx, aby$ attaches to the rest via its \OutIn{1}{2}-claimed pair $xy$, which is also $\{1,2\}$-claimed by the previous edges (in particular, one of these edges is $xyz$ for some vertex $z\in V\setminus \{a,b\}$), this gives 2 further pairs $\{3,4\}$-claimed by~$F$, namely $za$ and $zb$ (so 4 such pairs in total for the two diamonds). Thus,
$$w(F)\geq 17\,h(1,2)+4\,h(3,4)=17+4
>\frac{165}{61}\cdot 7,$$
 as desired.
So, by Lemma~\ref{lm:6Comb}\ref{it:>6}, we can assume that $F$ is made from a single edge $e$ by iteratively $\OMerge{}{2}$-merging $i\in [3,6]$ diamonds.
Then, all $5i+3$ pairs in $\p{1} F$ are $\{1,3\}$-claimed. Also, $F$ $\{3,5\}$-claims further $2i$ pairs. 
Indeed, each new diamond $abx,aby$ $2$-claims a pair $xy$ \OneClaim{}ed by some previous edge $xyz$, and since $i\geq 3$, the pairs $za$ and $zb$ are $\{3,5\}$-claimed by the final \TwoCluster\ $F$. Thus, we have
\begin{eqnarray*}
 w(F)&\ge& (5i+3)\,h(1,3)+2i\cdot h(3,5)\ =\ (5i+3)+2i\cdot \frac{11}{61}\\
 &=& \frac{18-3i}{61} + \frac{165}{61}\cdot (2i+1)
 \ \geq\  \frac{165}{61}\cdot |F|.
 \end{eqnarray*}

Thus, suppose that $|F|\leq 6$. By $(8,6)$-freenees, we have that $|F|\leq 5$. First, consider the case that $F=F_1\cup F_2$ for $\OMerge{}{2}$-mergeable  $F_1,F_2\in \OneClusters$ via some pair $uv$ (thus $1\in \CI_{F_1}(uv)\setminus \CI_{F_2}(uv)$).
Let $F_1$ be an $i$-tree and $F_2$ be a $j$-tree. 
Then, $F$ has $i+j$ edges and \OneClaim{}s $2i+2j+2$ pairs. Note that $j\geq 2$ as $F_2$ has to contain a diamond. Also,  since $G$ is $\cG{3}{6}$-free, the subgraphs $F_1$ and $F_2$ do not share any further vertices in addition to $u$ and $v$.

Suppose first that $i=1$. Every pair in $\p{1} F$ is $\{1,3\}$-claimed by~$F$. If $j=2$, then $F$ $3$-claims the remaining $2$ pairs and we have 
$$
w(F) =8\,h(1,3)+2\,h(3) =8 + 2\cdot \frac{6}{61} = \frac{500}{61} > \frac{165}{61}\cdot 3.
$$
If $j\geq 3$, then $F$ $\{3,4\}$-claims at least $2$ pairs not in $\p{1} F$. Hence, 
$$
w(F) \geq (2j+4)\,h(1,3)+2\,h(3,4)=(2j+4)+2 = \frac{(201- 43j)}{61} + \frac{165}{61}\, (j+1),
$$
which is at least $\frac{165}{61}\,(j+1)$ since $j= |F|-1\leq 4$, as desired.

Suppose that $i\geq 2$. Then, every pair in $\p{1} F$ is $\{1,2\}$-claimed by $F$. Also, given an edge $uvx$ in $F_1$ and a diamond $abu,abv$ in $F_2$, the \TwoCluster{} $F$ $\{3,4\}$-claims the pairs $ax, bx\notin \p{1}{F}$.
Moreover, if $i=2$, then the (unique) pair in ${V(F_1)\choose 2}\setminus \p{1}{F_1}$ is $\{2,4\}$-claimed by $F$; combining this with the fact that $j=|F|-i \leq 3$ yields 
\begin{eqnarray*}
w(F) &\geq& (2j+6)\, h(1,2)+2\, h(3,4)+ h(2,4)\ =\ (2j+6) + 2 + \frac{1}{2}\\
 &=&
 \  \frac{377-86j}{122}+\frac{165}{61}\, (j+2) >\ \frac{165}{61}\,|F|.
 \end{eqnarray*}
If $i=3$, then $j=2$. Again, $F$ $\{1,2\}$-claims all 12 pairs in $\p{1} F$ and $\{3,4\}$-claims another 2 pairs, but it also $\{2,3\}$-claims at least $2$ other pairs, namely the pairs 2-claimed but not \OneClaim{}ed by $F_1$. Hence, we have
$$
w(F) \geq 12\,h(1,2)+2\,h(3,4)+2\,h(2,3)=12 + 2 + 2\cdot \frac{36}{61} = \frac{926}{61} > \frac{165}{61}\cdot 5.
$$


Now, note that a 
\TwoCluster{} made of at least four \OneCluster{}s has at least 6 edges: indeed, this \TwoCluster{} was obtained by doing at least 3 consecutive $\Merge{\{1\}}{\{2\}}$-mergings, so some \OneCluster{} in it contains at least three edges or some two \OneCluster{}s in it contain at least two edges each.
Thus, it remains to consider the case when $F$ is obtained by merging three trees $F_1,F_2,F_3\in \OneClusters$. 
We cannot have $|F|\leq 4$ as then, the \TwoCluster{} $F$ would be made of at least two 1-trees and at most one 2-tree, which is impossible for $3$-graphs. Hence, we obtain that $|F|=5$ with the  composition $(3,1,1)$ or $(2,2,1)$.

If $F$ has composition $(3,1,1)$, then $F$ $\{1,3\}$-claims all $13$ pairs in $\p{1} F$ and $\{3,4\}$-claims at least $4$ other pairs (namely, for each 1-tree $xyz$ and the corresponding diamond $aby,abz$ inside the 3-tree, the pairs $ax$ and $bx$ are $\{3,4\}$-claimed by $F$).  We have that
$$
w(F)\geq 13\,h(1,3)+4\,h(3,4)=13+4=\frac{1037}{61}>\frac{165}{61}\cdot 5.
$$

Suppose that $F$ has composition $(2,2,1)$. Then, $F$ $\{1,3\}$-claims  all $13$ pairs in $\p{1} F$.
Also, $F$ can be built from its $1$-tree by attaching each new diamond $abx,aby$ via a previously \OneClaim{}ed pair $xy$, say by $xyz\in F$; here each of the pairs $az$ and $bz$ is $\{3,5\}$-claimed by the final \TwoCluster{}~$F$. Thus, we have that
$$
w(F)\geq 13\,h(1,3)+4\,h(3,5)=13+4\cdot \frac{11}{61}=\frac{837}{61}>\frac{165}{61}\cdot 5,
$$
which finishes the case analysis and the proof.
\epf

Finally, we prove the following claim. 

\begin{claim}\label{cl:wFLower2}
For all $F\in \ThreePClusters\setminus \TwoClusters$, we have that 
$w(F)\geq \frac{165}{61}\,|F|$.	
\end{claim}
\bpf[Proof of Claim~\ref{cl:wFLower2}]
Recall that if $F_1,F_2\in\TwoClusters$ are $\OThreePMerge$-mergeable via $uv$, then $1\in \CI_{F_1}(uv)$ and $3 \in \CI_{F_2}(uv)$ and, in addition, either $2\in \CI_{F_1}(uv)$ or $4\in \CI_{F_2}(uv)$. In particular, we have that $1,2\notin \CI_{F_2}(uv)$ (as otherwise we would have already merged $F_1$ and $F_2$ when constructing~$\TwoClusters$); also, by~\eqref{eq:2} $5\not\in \CI_{F_2}(uv)$ and $3\not\in\CI_{F_1}(uv)$.

\begin{figure}
\begin{center}
\includegraphics[scale=0.45]{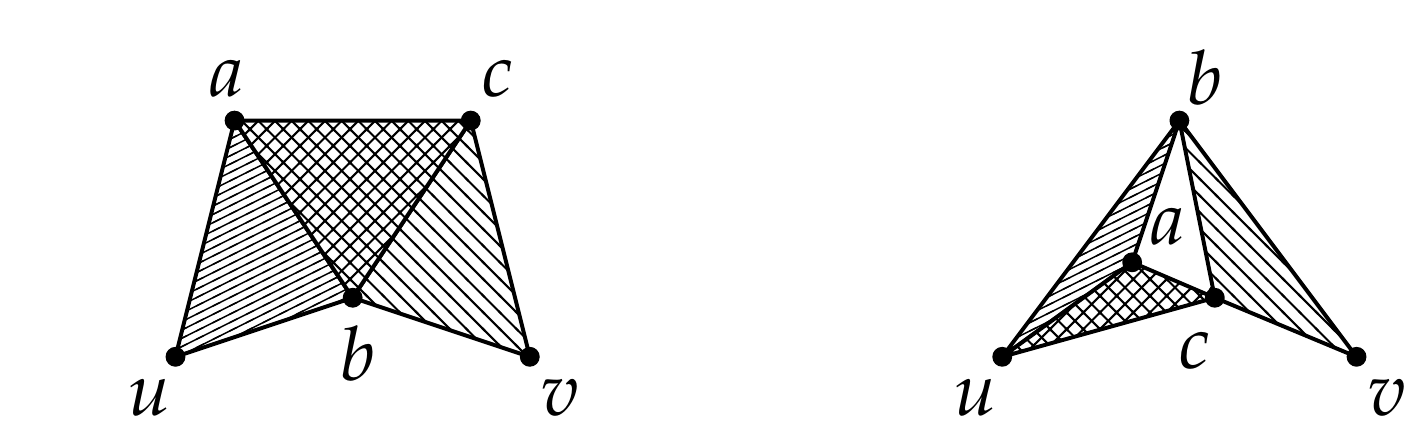}
\end{center}
\caption{Configurations $P_3$ and $C_3$.}\label{fi:P3C3}
\end{figure}

Let $F\in \ThreePClusters\setminus \TwoClusters$ be made of $F_1,\dots,F_s \in \TwoClusters$ $\ThreePMerge$-merged in this order as in Claim~\ref{cl:3Trim}.
Assume without loss of generality that $F_1$ and $F_2$ are $\OMerge{}{3^+}$-mergeable via some pair~$uv$. Let $F':=F_1\cup F_2$.
As $1,2\notin \CI_{F_2}(uv)$ and $3\in \CI_{F_2}(uv)$, we know that there is a $(5,3)$-configuration in $F_2$ containing $uv$ which is either a $3$-path $P_3=\{uab,abc,bcv\}\subseteq F_2$ or a \TwoCluster\ $C_3=\{aub,auc,bcv\}\subseteq F_2$ 
(which is the union of a $1$-tree and a $2$-tree that are $\OMerge{}{2}$-mergeable), see Figure~\ref{fi:P3C3}.
Since $3\notin \CI_{F_1}(uv)$ by \eqref{eq:2}, $F_1$ cannot be an $i$-tree for $i\geq 3$. We split the proof into $3$ cases depending on whether $F_1$ is a $2$-tree, a $1$-tree or an element of $\TwoClusters\setminus \OneClusters$.

\case{1}{Assume that $F_1$ is a $2$-tree.}%
Then, $1,2\in \CI_{F_1}(uv)$. By \eqref{eq:2}, we have $4,5\notin \CI_{F_2}(uv)$, and recall that $F_2$ contains one of $P_3$ or $C_3$ as described above. Also, the union $F'':=F_1\cup P_3$ or $F_1\cup C_3$ is a $(7,5)$-configuration, so  
\beq{eq:F1C3}
 \binom{V(F'')}{2}\cap \p{1}{G\setminus F''}=\emptyset.
 \eeq
Moreover, since $5\notin \CI_{F_2}(uv)$, no pair in $V(P_3)$ (resp.\ $V(C_3)$) can be $2$-claimed by $G\setminus F'$. It follows that the \TwoCluster{} $F_2$ is equal to $P_3$ or $C_3$, and thus $F''=F'$. 

Suppose that $F$ consists of $s\geq 3$ 2-clusters. Let the next $\ThreePMerge$-merging step (of  $F'$ and $F_3$) be via a pair~$u'v'$. It follows from~\eqref{eq:F1C3} (and from $F''=F'$) that $u'v' \in \p{1} {F'}\cap \p{3}{F_3}$. By \eqref{eq:2}, it holds that $3\notin \CI_{F'}(u'v')$. 
As $F_2$ contains 3 edges, necessarily $u'v'\in \p{1} {F_1}$. 
Also, since $u'v'$ is in $\p{1}{F_1}\subseteq \p{2}{F_1}$, we see that $F_1$ alone and $F_3$ are $\OThreePMerge$-mergeable via~$u'v'$.
Now, our previous argument about the structure of $F_2$ applies to $F_3$ as well and shows that $F_3$ is isomorphic to a copy of $P_3$ or $C_3$ 3-claiming $u'v'$.
Furthermore, the $(5,3)$-configurations $F_2$ and $F_3$ can share at most one vertex, so ${V(F_2)\choose 2}\cap {V(F_3)\choose 2}=\emptyset$. The very same reasoning applies in turn to each $F_i$ with $i\in[4,s]$ with the analogous conclusion. Indeed, by Lemma~\ref{lm:3Merge}, there is $j\in [1,i-1]$ such that $F_i$ and $F_j$ $\ThreePMerge$-merge. 
Together with the $(8,6)$-freeness and the fact that $F_1\cup F_2, \ldots, F_1\cup F_{i-1}$ form $(7,5)$-configurations, $F_i$ can only $\ThreePMerge$-merge with $F_1$.  
Let us bound from below the total weight assigned by $F$, being rather loose in our estimates. 
The five pairs in $\p{1} {F_1}$ are $\{1,2\}$-claimed (and we just ignore the remaining pair inside $V(F_1)$ which is $2$-claimed). 
For each $i\in [2,s]$, if $F_i$ is a 3-path, then it $\{1,3\}$-claims all seven pairs in its $2$-shadow $\p{1} {F_i}$ and two further pairs inside $V(F_i)$ are $\{2,3,4\}$-claimed by~$F_1\cup F_i\subseteq F$ (for example, if $i=2$, then these are the pairs $uc$ and $av$). All these pairs are unique to $F_i$ and thus $F_i$ contributes at least $7\,h(1,3)+2\,h(2,3,4)=9$ to $w(F)$. 
Also, if $F_i$ is isomorphic to $C_3$, then it $\{1,3\}$-claims all 8 pairs in $\p{1} {F_i}$ and one further pair inside $V(F_i)$  is $\{3,4\}$-claimed by $F_1\cup F_i$ (for example, if $i=2$ then it is the pair $av$). Thus, $F_i$ contributes at least $8\,h(1,3)+h(3,4)=9$ to $w(F)$. We conclude by $s\geq 2$ (which follows from $F\notin\TwoClusters$) that
\begin{eqnarray*}
 w(F) &\geq& 5\,h(1,2)+9(s-1)\ =\ 5   +9(s-1)\\
 &=& \frac{ 54s -79 }{61} + \frac{165}{61}\cdot (3s-1)\ >\ \frac{165}{61}\,|F|,
\end{eqnarray*}
 as desired.

\case{2}{Assume that $F_1$ is a $1$-tree $\{uvw\}$.}%
As $F_1$ and $F_2$ are $\OThreePMerge$-mergeable via $uv$, we have $3,4\in \CI_{F_2}(uv)$. Thus, $|F_2|\geq 4$ and $|F_2|\notin \{5,6\}$ (otherwise, we would get a $(8,6)$-configuration). Also, 
$1,2\notin\CI_{F_2}(uv)$ since $F_1,F_2\in\TwoClusters$ are distinct.
Hence, as before, $F_2$ contains a copy of $P_3$ or $C_3$ 3-claiming $uv$.

Suppose first that $|F_2|\geq 7$. Then, $F_2$ has the structure given by Lemma~\ref{lm:6Comb}\ref{it:>6}, consisting of a 1-tree or 3-tree with a number of diamonds merged in one by one. A $(5,3)$-configuration in $F_2$ that contains $uv$ involves at most two \OneCluster{}s of $F_2$, whose union $F''$ has at most $5$ edges; moreover, if $F''$ consists of two \OneCluster{}s, then these \OneCluster{}s are $\Merge{}{2}$-mergeable. The proof of Lemma~\ref{lm:6Comb}\ref{it:>6} shows that, if we build $F_2\in\TwoClusters$ by starting with the partial \TwoCluster{} $F''$ and attaching \OneCluster{}s one by one, then we reach a $(7,5)$-configuration $F'''$ just before the number of edges jumps over~$6$. 
However, $F_1$ shares at least $2$ vertices with $F'''\supseteq F''$, so $F_1\cup F'''$ is a $(8,6)$-configuration. This contradiction shows that $|F_2|=4$.

Now, denote the $4$-edge hypergraph $F_2$ by $T_4$ if $F_2$ consists of (a copy of) $P_3$ that 3-claims $uv$ together with one more edge (thus, $T_4$ is a $4$-tree or a \TwoCluster\ made of a $3$-tree and a $1$-tree as its \OneCluster{}s), and denote $F_2$ by $C_4$ if $F_2$ consists of (a copy of) $C_3$ that 3-claims $uv$ with one more edge (thus, the \OneCluster{}s of $C_4$ are two $2$-trees, or a $3$-tree and a $1$-tree).
Recall that $F'=F_1\cup F_2$.

Suppose first that $s\geq 3$. Let $F_3$ and $F'$ be merged via a pair $u'v'$.
The $3$-graphs $F_3$ and $F'$ cannot be $\OThreePMerge$-mergeable via $u'v'$ as otherwise an edge in $F_3$ containing $u'v'$ together with $F'$ would be a $(8,6)$-configuration. So $F'$ and $F_3$ are $\OThreePMerge$-mergeable, that is, $1\in\CI_{F'}(u'v')$ and $3\in \CI_{F_3}(u'v')$. Then, $3\notin \CI_{F'}(u'v')$. This greatly limits the number of possibilities for the pair $u'v'$ inside $\p{1}{F'}$.
If $F_2=T_4$, one can easily notice that $u'v'$ cannot be in $\p{1} {F_2}\cup \{uv\}$ and thus \ $u'v'\in \p{1} {F_1}\setminus \{uv\}$. Hence, $F_1$ alone and $F_3$ are $\OThreePMerge$-mergeable and the above analysis for $F_2$ applies to $F_3$ as well, showing that the \TwoCluster{} $F_3$ is isomorphic to a copy of $T_4$ or $C_4$ 3-claiming $u'v'$. If $F_2=C_4$, then either $u'v'\in \p{1} {F_1}\setminus \{uv\}$ 
(and, again, $F_3$ is isomorphic to $T_4$ or $C_4$) or, for some vertices $a,b,c,d$, we have $F_2=\{aub,auc,bcv,cdv\}$ and $u'v'$ is $cd$ or $dv$; also, the argument of Case~$1$ (with $\{bcv,cdv\}$ playing the role of the 2-tree $F_1$ from Case~$1$) shows that $F_3$ rooted at $u'v'$ is isomorphic to $P_3$ or $C_3$, no other $F_i$ with $i\geq 4$ can be $\ThreePMerge$-merged with $F_3$, and all pairs in ${V(F_3)\choose 2}$ are unique to $F_3$. Using Lemma~\ref{lm:3Merge}, the same argument applies to each new $F_i$ with $i\geq 4$.
To summarise, we obtained that each $F_i$ for $i\geq 2$ is either some instance of $T_4$ or $C_4$ $\ThreePMerge$-merged with the single edge $F_1$, or an instance of $P_3$ or $C_3$ $\ThreePMerge$-merged with a copy of $C_4$ as specified above; also, the only pairs shared between these \TwoCluster{}s are the pairs along which these $\ThreePMerge$-mergings occur.

Now, assume the final $F$ consists of one $1$-tree, $i$ copies $T_4$ or $C_4$, and $j$ copies of $P_3$ or~$C_3$. (Although we could say more about the structure of $F$, e.g.\ that $i\leq 3$ and $j\leq 2i$, these observations are not needed for our estimates.)
Then, $|F|=1+4i+3j$. Each copy of $T_4$ (which is a 4-tree or a 3-tree $\Merge{}{}{2}$-merged with a single edge) has, in addition to the pair via which it is merged with $F_1$, at least nine $\{1,3\}$-claimed pairs and other two $\{2,3,4\}$-claimed pairs. Thus, it contributes at least $9\,h(1,3)+2\,h(2,3,4)=11$ to $w(F)$. 
Likewise, each copy of $C_4$ contributes at least $8\,h(1,3)+2\,h(1,2)+h(3,4)=11$ to~$w(F)$. Also, as in Case~1, each copy of $P_3$ or $C_3$ contributes at least~$9$ to~$w(F)$. Additionally, we have $3$ pairs in $\p{1} {F_1}$ which are $\{1,4\}$-claimed by $F_1\cup F_2$.
 Hence, we get the required:
 \begin{eqnarray*}
 w(F)&\geq& 3\,h(1,4)+11i+9j \ =\ 3\cdot \frac{55}{61}+ 11i+9j\\
 & =&\frac{11i+54j}{61}+\frac{165}{61}\, (1+4i+3j)\ >\  \frac{165}{61}\, |F|.
 \end{eqnarray*}

\case{3}{Assume that $F_1\in \TwoClusters\setminus \OneClusters$.}%
Here, $F_1$ is a $2$-merging of at least two \OneCluster{}s and thus has at least 3 edges.
Let $T_1\in\OneClusters$ be the \OneCluster\ in $F_1$ 1-claiming~$uv$ (recall that $uv$ is the pair via which $F_1$ and $F_2$ are $\OThreePMerge$-mergeable). The tree $T_1$ has at most 2 edges as otherwise $T_1$ together with $F_2$ would contain an $(8,6)$-configuration, a contradiction. 
Also, $T_1$ cannot be a $1$-tree as otherwise $F_1\setminus T_1$ would contain a diamond $2$-claiming a pair in $\p{1} {T_1}$ (since $F_1\in \TwoClusters\setminus \OneClusters$), which implies that $3\in \CI_{F_1}(uv)$, a contradiction.
Therefore, $T_1$ is a $2$-tree. As $F_1\in\TwoClusters\setminus\OneClusters$,  $T_1$ has to be $2$-merged with some other \OneCluster\ $T_2\subseteq F_1$. 
It is impossible that $T_2$ and $T_1$ are $\OMerge{}{2}$-mergeable as otherwise $T_1$ plus an edge of $T_2$ would be a $(5,3)$-configuration containing $uv$ in $F_1$, which contradicts $3\notin \CI_{F_1}(uv)$. Thus, $T_1$ and $T_2$ are $\OMerge{}{2}$-mergeable.  
Also, $T_2$ has at most 3 edges since trees with more edges would form an $(8,6)$-configuration with~$T_1$. It is routine to check that we can assume (after swapping $u$ and $v$ if necessary) that, for some vertices $a,b\in V$, the 2-tree $T_1$ is $\{avu,avb\}$ and the pair $\Merge{}{2}$-claimed by $T_2$ is $ab$ or~$bv$. 
Also, any further 2-merging involving $T_1\cup T_2$ would cause an $(8,6)$-configuration. We conclude that $F_1=T_1\cup T_2$.  

The same argument as in Case 1 shows that $F_2$ is given by $P_3$ (a 3-path) or $C_3$ (a diamond and a single edge), see Figure~\ref{fi:P3C3}. 
Let $F':=F_1\cup F_2$. We have $2$ subcases depending on~$T_2$.

\case{3-1}{Assume that $|T_2|=2$.}%

Suppose first that $s\geq 3$. If $F_3$ and $F'$ are $\OThreePMerge$-mergeable via some pair $u'v'$ then, by $(8,6)$-freeness, $u'v'$ must be one of the two pairs  \OutIn{1}{3}-claimed by $F_1$, and $F_3$ is a $1$-tree; therefore, by Claim~\ref{cl:3Trim}, we can reorder \TwoCluster{}s constituting $F$ starting with $F_3$ and follow the analysis in the Case 2. Now assume that both pairs \OutIn{1}{3}-claimed by $F_1$ are not used for further $\ThreePMerge$-mergings. Thus, $F'$ and $F_3$ are $\OThreePMerge$-mergeable via $u'v'$. 
Then, since $3\in \CI_{F_3}(u'v')$, $u'v'$ must be in $\p{1} {F'}\setminus \p{3}{F'}=\{au\}$ (recall that $3\in \CI_{F_2}(uv)\subseteq \CI_{F'}(uv)$), 
$F_3$ rooted at $u'v'=au$ is isomorphic to $P_3$ or $C_3$; moreover no further $\ThreePMerge$-mergings are possible and thus $s=3$. 
Thus, for both $s=2$ and $s=3$, we can assume the final \ThreePCluster\ $F$ is made of
$F_1$ and $j$ copies of $P_3$ or $C_3$ where $j\in \{1,2\}$. Then, $|F|=4+3j$, $F_1$ contributes at least $8\,h(1,3)+2\,h(1,2)+2h(3,4)+h(2,4)=10+2+\frac{1}{2}=\frac{25}{2}$ to $w(F)$. Hence, we have 
\begin{eqnarray*}
w(F)&\geq& 
\frac{25}{2}+9j
\ =\
 \frac{108j+205}{122}+ \frac{165}{61}\, (4+3j)\ >\ \frac{165}{61}
\,|F|,
\end{eqnarray*}
as desired.

\case{3-2}{Assume that $T_2$ is a 3-tree.}%
Here $\CI_{F_1}(uv)\supseteq\{1,2,4,5\}$ so $F_2$ has exactly 3 edges and 
$F'=F_1\cup F_2$ is a $(10,8)$-configuration. Suppose first that $s\geq 3$. Since no $3$-claimed pair of $F'$ can be \OneClaim{}ed by $G\setminus F'$,
the $3$-graphs $F_3$ and $F'$ cannot be $\OThreePMerge$-mergeable. So let $F'$ and $F_3$ be $\OThreePMerge$-mergeable via some pair $u'v'$. Since $3\in \CI_{F_3}(u'v')$, $u'v'$ must be in $\p{1} {F'}\setminus \p{3}{F'}=\{au\}$; furthermore, $F_3$ rooted at $u'v'=au$ is isomorphic to $P_3$ or $C_3$, no further $\ThreePMerge$-mergings are possible and $s=3$. Thus, for both $s=2$ and $s = 3$, the final \ThreePCluster\ $F$ consists of $F_1$ and $j$ copies of $P_3$ or $C_3$ where $j\in \{1,2\}$. Here, $|F|=5+3j$. Note that $F_1$ contributes at least $10\,h(1,3)+2\,h(1,2)+4\,h(3,4)=12+4=16$ to the total weight. We have
$$
 w(F)\ \geq\ 16+9j\ =\ \frac{54j+151}{61}+\frac{165}{61}\,(5+3j)\ >\  \frac{165}{61}\, |F|.$$

This finishes the proof of the claim.
\epf

Hence, by the previous claims, we conclude that 
$$
|G|=\sum_{F\in \ThreePClusters } |F| \leq \frac{61}{165} \sum_{F\in \ThreePClusters } w(F) = \frac{61}{165} \sum_{uv\in \binom{V(H)}{2}} w(uv) \leq \frac{61}{165}  \binom{n}{2}.
$$
This proves Theorem~\ref{6edge3uniform}.
\epf

\section{Concluding remarks}
In this paper, we made progress on the Brown--Erd\H{o}s--S\'os Problem (the case of 3-uniform hypergraphs) with $k\in \{5,6,7\}$ edges and its extension to $r$-graphs. 
We note that a further extension of the Brown--Erd\H{o}s--S\'os Problem proposed by Shangguan and Tamo~\cite{ShTamo} asks to determine whether the limits
$$\lim_{n\rightarrow \infty}n^{-t}f^{(r)}(n;k(r-t)+t,k)$$
exist for all fixed $r,t,k$ and, if so, to find their values. (The case that we studied here corresponds to $t=2$.) 

Let us briefly summarise what is known. The results in~\cites{BES,Rodl85} resolve the case $k=2$.
In~\cites{G,GJKKLP,ShTamo}, this problem was completely solved  when $k\in \{3,4\}$.
In~\cites{DP,sh}, the existence of the limit was proved for $t=2$.
In
\cite{letzter2023problem}, the value of the limit (and thus its existence) was established for even $k$ when $r\gg k,t$ is sufficiently large. Also, it was proved in \cite{letzter2023problem} that if $k\in \{5,7\}$ then the limit exists for any $r$ and $t$.
Our results determine the limit values when $t=2$ and $k\in \{5,6,7\}$.

It would be interesting to study the existence of limits for the remaining sets of parameters as well as their precise values.

\section*{Acknowledgments} 
We thank Felix Joos and Marcus K\"uhn for useful discussions. In particular, some starting ideas of the present paper originated from the joint work~\cite{GJKKLP}. Also, we thank the anonymous referee for useful comments.

\begin{Backmatter}
 


\printaddress

\end{Backmatter}

\end{document}